%% file: main4arxiv.tex
\documentclass[10pt,a4paper]{article}
\usepackage[a4paper,left=2cm,right=2cm,top=2.5cm,bottom=2.5cm]{geometry}

\usepackage{framed,multirow}

\usepackage{amssymb}
\usepackage{latexsym}

\usepackage{amsmath}
\usepackage{amsthm}
\usepackage{todonotes}
\usepackage{url}
\usepackage[ruled,linesnumbered]{algorithm2e}
\usepackage{hyperref}
\usepackage{bm}
\usepackage[makeroom]{cancel}
\usepackage{subcaption}

\theoremstyle{definition}

\newcommand*{\mv}[1]{\bm{#1}}
\newcommand*{\mt}[1]{\bm{\mathcal{#1}}}
\newcommand{\mvec}[1]{\texttt{vec}\left[#1\right]}
\newcommand{\qr}[1]{\texttt{qr}\left[#1\right]}

\newcommand{\diag}[1]{\text{diag}\left(#1\right)}

\newcommand{\Oh}[1]{\mathcal{O}\left(#1\right)}

\usepackage{pgfplots}
\usepackage{tikz}

\newcommand{\N}{\mathbb{N}}

\newcommand{\C}{\mathbb{C}}
\newcommand{\T}{{\mkern-1.5mu\mathsf{T}}}
\renewcommand{\H}{\mathsf{H}}
\newcommand*{\conj}[1]{\overline{#1}}

\usepackage{url}
\definecolor{newcolor}{rgb}{.8,.349,.1}
\usepackage{xcolor}

\begin{document}

\title{Solving for the low-rank tensor components of a scattering wave function}
\author{Jacob Snoeijer and Wim Vanroose}
\maketitle

\begin{abstract}
  Atomic and molecular breakup reactions, such as multiple-ionisation,
  are described by a driven Schr\"odinger equation.  This equation is
  equivalent to a high-dimensional Helmholtz equation and it has
  solutions that are outgoing waves, emerging from the target.  We
  show that these waves can be described by a low-rank
  approximation.  For 2D problems this it a matrix product of two
  low-rank matrices, for 3D problems it is a low-rank tensor decomposition.
  We propose an iterative method that solves, in an alternating way, for
  these low-rank components of the scattered wave. We illustrate the
  method with examples in 2D and 3D.
\end{abstract}



\setcounter{secnumdepth}{4}

\section{Introduction}
An in-coincidence experiment measures simultaneously the outgoing
momenta of multiple products of a microscopic reaction
\cite{schmidt2002coltrims}. It is an instrument that can study the
correlations in reactions involving multiple particles.  In double
ionization, for example, a single photon ionizes, simultaneously, two
electrons and the outgoing momenta of both particles are captured
\cite{akoury2007simplest}. The reaction probes the correlation between
two electrons in, for example, a chemical bound at the moment of
photon impact.  The outgoing wave of the two electrons is described by
a 6D correlated wave and results in a cross section that depends on
four angles, the directions of the first and the second electron.

Free-electron lasers, and similar experiments around the world, are
expected to generate a wealth of this high-dimensional scattering
data.  This will result in high-dimensional forward and inverse wave
problems that need to be solved to interpret the data.

The experimental cross section are often smooth functions as a
function of the angles. Similarly, some parts of the scattering
solution, such as single-ionization, only probes a limited subspace
of the possible full solution space. The scattering solution can then
be described by a low-rank wave function, a product of one-particle
bound states with scattering waves in the other coordinates.

This paper introduces a low-rank representation for the
scattering solutions, not only for the single ionization but also for
double and triple ionization waves that appear in breakup reaction.

We also propose and analyze an alternating direction algorithm that
directly solves for the low-rank components that describe the
solution. This reduces a large-scale linear system to smaller,
low-dimensional, scattering problems that are solved in a iterative
sequence.  The proposed method can be generalized to high-dimensional
scattering problems where a low-rank tensor decomposition is used to
represent the full scattering wave function.

Efficient low-rank tensor representations are used in quantum physics
for quite some time already \cite{huckle2013computations,
  murg2010simulating}.  They are also used in the applied mathematics
literature to approximate high-dimension problems, for a review see
\cite{grasedyck2013literature, hackbusch2012tensor, kolda2009tensor}.
Methods such ALS \cite{holtz2012alternating}, DMRG
\cite{oseledets2011dmrg}, and AMEn \cite{dolgov2014alternating} use in
alternating directions, a small linear system to determine the
low-rank components of a tensor decomposition. These innovations have
not found their application in computational scattering theory.

To calculate cross sections, from first principles, we start
from a multi-particle Schr\"odinger equation. The equation is
reformulated into a driven Schr\"odinger equation with an unknown scattering
wave function and a right hand side that describes the
excitation, for example, a dipole operator working on the initial
state.

Since the asymptotic behaviour of a scattering function for multiple
charged particles is in many cases unknown, absorbing boundary
conditions \cite{ECS,PML} are used.  Here, an artificial layer is
added to the numerical domain that dampens outgoing waves. The
outgoing wave boundary conditions are then replaced with
homogeneous Dirichlet boundary conditions at the end of the artificial
layer. These boundary do not require any knowledge about the
asymptotic behaviour, which becomes very complicated for these
multiple charged particles.

The resulting equation is discretized on a grid and results in a
large, sparse indefinite linear system.  It is typically solved by a
preconditioned Krylov subspace method \cite{cools2016fast}.  However,
the preconditioning techniques for indefinite systems are not as
efficient as preconditioners for symmetric and positive definite
system.  And solving the resulting equation is still a
computationally expensive task, often requiring a distributed
calculation on a supercomputer.

To compare the resulting theoretical cross sections with experimental
data, a further postprocessing step is necessary. The cross section is
the farfield map and this is calculated through integrals of the
scattering wave function, which is the solution of the linear system,
and a Greens function \cite{mccurdy2004solving}.

The main result of the paper is that we show that scattering waves
that describe multiple ionization can be represented by a low-rank
tensor. We first show this for a 2D wave and then generalize the
results to 3D waves.  The methodology can be generalized to higher
dimension.

The outline of the paper is a follows. In section~\ref{sec:stateofart} we
review the methodology that solves a forward scattering problem. It
results in a driven Schr\"odinger equation with absorbing boundary
conditions.  From the solution we can extract the cross section using
an integral.  In section~\ref{sec:lowrank} we illustrate, in 2D, that a
solution can be approximated by a truncated low-rank approximation. We
also show that these low-rank components can be calculated directly
with an iterative method.  In section~\ref{sec:3dhelmholtz} we show that
this methodology generalizes to 3D and higher dimensional problems.
We use a truncated tensor decomposition and determine the components 
with a similar iterative method.
A discussion of some numerical results and a comparison of the different
presented versions of the method is given in section~\ref{sec:NumericalResults}.
In the final section, Sec.~\ref{sec:discussion}, we summarize some
conclusions and discuss some possible extensions of the presented method.

\section{State of the art} \label{sec:stateofart}
This section summarize the methodology that solves a forward
break-up problems with charged particles. The methodology is developed
in a series of papers \cite{rescigno2000numerical,mccurdy2004solving}
and applied to solve the impact-ionization problem
\cite{rescigno1999collisional} and double ionization of molecules
\cite{vanroose2005complete, vanroose2006double}.  These methods are
being extended to treat, for example, water
\cite{streeter2018dissociation}.

The helium atom, He, is the simplest system with double ionization
\cite{briggs2000differential}. It has two electrons with coordinates
$\mathbf{r}_1 \in \mathbb{R}^3$ and $\mathbf{r}_2 \in \mathbb{R}^3$
relative to the nucleus positioned at $\mathbf{0}$.  The driven
Schr\"odinger equation for $u(\mathbf{r}_1,\mathbf{r}_2) \in C^2$ then
reads
\begin{equation}\label{DrivenSchrodinger}
\left( {- \frac{1}{2} \Delta_{\mathbf{r}_1} - \frac{1}{2} \Delta_{\mathbf{r}_2} }  { - \frac{1}{\|\mathbf{r}_1\|}- \frac{1}{\|\mathbf{r}_2\|}} + {\frac{1}{\|\mathbf{r}_1-\mathbf{r}_2\|} }- { E} \right) u(\mathbf{r}_1,\mathbf{r}_2) =   \mathcal{\mu}  \phi_0(\mathbf{r}_1,\mathbf{r}_2) \quad  \forall \mathbf{r}_1, \mathbf{r}_2 \in \mathbb{R}^3,
\end{equation}
where the right hand side is
the dipole operator $\mathcal{\mu}$ working on the ground state $\phi_0$, the
eigenstate with the lowest energy $\lambda_0$.  The operators
$-\frac{1}{2} \Delta_{\mathbf{r}_1}$ and
$-\frac{1}{2}\Delta_{\mathbf{r}_2}$ are the Laplacian operators for
the first and second electron and model the kinetic energy.  The
nuclear attraction is $-1/\|\mathbf{r}_1\|$ and $-1/\|\mathbf{r}_2\|$
and the electron-electron repulsion is
$1/\|\mathbf{r}_1-\mathbf{r}_2\|$.

The total energy $E = h\nu + \lambda_0$ is the energy
deposited in the system by the photon, $h\nu$, and the energy $\lambda_0$ of
the ground state. If the $E>0$, both electrons can
escape simultaneously from the system. The solution
$u(\mathbf{r}_1,\mathbf{r}_2)$ then represents a 6D wave emerging from
the nucleus.

The equation can be interpreted as a Helmholtz equation with a
space-dependent wave number, $k^2(\mathbf{r}_1,\mathbf{r}_2)$,
\begin{equation}\label{eq:helmholtz}
 \left( - \Delta_{6D} - {  k^2(\mathbf{r}_1,\mathbf{r}_2)} \right)
 u(\mathbf{r}_1,\mathbf{r}_2) = {  f(\mathbf{r}_1,\mathbf{r}_2)} \quad \forall (\mathbf{r}_1,
 \mathbf{r}_2) \in \mathbb{R}^6.
\end{equation}

In this paper we prefer to write this Helmholtz equation as
\begin{equation}
   \left( - \Delta_{6D}   -   k_0^2  (1 + \chi(\mathbf{r}_1,\mathbf{r}_2) ) \right)  u(\mathbf{r}_1,\mathbf{r}_2)  =  f(\mathbf{r}_1,\mathbf{r}_2) \quad  \forall (\mathbf{r}_1, \mathbf{r}_2) \in \mathbb{R}^6,
\end{equation}
where $k_0^2$ is a constant wave number, in this case related to the
total energy $E$, and a space-dependent function $\chi: \mathbb{R}^6
\rightarrow \mathbb{R}$, that goes to zero if $\|\mathbf{r}_1\|
\rightarrow \infty $ or $\|\mathbf{r}_1\| \rightarrow \infty $ that
represents all the potentials.

\subsection{Expansion in spherical waves and absorbing boundary conditions}
For small atomic and molecular system, where spherical symmetry is relevant, the
system is typically written in spherical coordinates and expanded in
spherical harmonics.  With $\mathbf{r}_1(\rho_1, \theta_1,
\varphi_1)$ and $\mathbf{r}_2(\rho_2, \theta_2, \varphi_2)$ we can
write
\begin{equation}\label{eq:partialwaves}
 u(\mathbf{r}_1,\mathbf{r}_2) = \sum_{l_1}^\infty \sum_{m_1 = -l_1}^{l_1} \sum_{l_2=0}^{\infty} \sum_{m_2=-l_2}^{l_2}  u_{l_1m_1,l_2m_2}(\rho_1,\rho_2)  Y_{l_1m_1}(\theta_1,\varphi_1) Y_{l_2m_2}(\theta_2,\varphi_2),
\end{equation}
where $Y_{l_1m_1}(\theta_1, \varphi_1)$ and $Y_{l_2m_2}(\theta_2,
\varphi_2)$ are spherical Harmonics, the eigenfunctions of the angular
part of a 3D Laplacian in spherical coordinates. In practice the sum in equation~\eqref{eq:partialwaves} is
truncated. The expansion is then a low-rank, truncated, tensor decomposition
of a 6D tensor describing the solution.

For each $l_1$, $m_1$, $l_2$ and $m_2$ combination, the radial function
$u_{l_1m_1,l_2m_2}(\rho_1,\rho_2)$ describes an outgoing wave that
depends on the distances $\rho_1$ and $\rho_2$ of the two electrons to
the nucleus.

A coupled equation that simultaneously solves for all the
$u_{l_1m_1l_2m_2}(\rho_1,\rho_2)$'s is found by inserting the
truncated sum in \eqref{DrivenSchrodinger}, multiplying with $Y_{l_1m_1}^*(\theta_1,
\varphi_1)$ and $Y_{l_2m_2}^*(\theta_2,\varphi_2)$ and integrating
over all the angular coordinates,
\begin{equation}\label{eq:coupledpartialwaves}
  \begin{aligned}\left(-\frac{1}{2}\frac{d^2}{d \rho_1^2} + \frac{l_1(l_1+1)}{2\rho_1^2 }
 -\frac{1}{2}\frac{d^2}{d \rho_2^2} + \frac{l_2(l_2+1)}{2\rho_2^2 } + V_{l_1m_1l_2m_2}
(\rho_1,\rho_2) -E \right) u_{l_1m_1l_2m_2}(\rho_1,\rho_2)   \\ + \sum_{l^\prime_1m^\prime_1l^\prime_2m^\prime_2} V_{l_1m_1l_2m_2, l^\prime_1m^\prime_1l^\prime_2m^\prime_2}(\rho_1,\rho_2) u_{l^\prime_1m^\prime_1l^\prime_2m^\prime_2}(\rho_1,\rho_2) = f_{l_1m_1,l_2m_2}(\rho_1,\rho_2) \quad \forall l_1m_1 l_1m_2. \quad \forall \rho_1,\rho_2 \in [0,\infty[
\end{aligned}
\end{equation}
with boundary conditions $u(\rho_1=0,\rho_2)=0$ for all $\rho_2 \ge 0$ and $u(\rho_1,\rho_2=0)=0$ for all $\rho_1 \ge 0$.

The equation \eqref{eq:coupledpartialwaves} is typically discretized
on a spectral elements quadrature grid \cite{rescigno2000numerical}.

To reflect the physics, where electrons are emitted from the system,
outgoing wave boundary conditions need to be applied at the outer
boundaries. There are many ways to implement outgoing wave boundary
conditions. Exterior complex scaling (ECS) \cite{ECS} for
example, is frequently used in the computational atomic and molecular
physics literature. In the computational electromagnetic scattering a
perfectly matched layers (PML) \cite{PML} is used, which can also be
interpreted as a complex scaled grid. \cite{chew19943d}

\subsection{Calculation of the amplitudes}
To correctly predict the probabilities of the arriving particles at
the detector, we need the amplitudes of the solution far away from the
molecule. These are related to the asymptotic amplitudes of the
wave functions.

Let us go back to the formulation with Helmholtz equation, as given in
\eqref{eq:helmholtz}. Suppose that we have solved the following
Helmholtz equation with absorbing boundary conditions, in any
representation,
\begin{equation} \label{Helmholtz}
  \left(- \Delta -k_0^2 \left( 1+ \chi(\mathbf{x}) \right) \right) u_{\text{sc}}(\mathbf{x}) = f(\mathbf{x}),  \quad \forall \mathbf{x} \in [-L,L]^d,
\end{equation}
where $f$ is only non-zero on the real part of the grid $[-L,L]^d
\subset \mathbb{R}^d$. Similarly, $\chi(\mathbf{x})$ is only non-zero
on the box $[-L,L]^d$.

The calculation of the asymptotic amplitudes requires the solution
$u_{\text{sc}}(\mathbf{x})$ for an $\mathbf{x}$ outside of the box
$[-L, L]^d$. To that end, we reorganize equation~\eqref{Helmholtz},
after we have solved it, as follows
\begin{equation} \label{eq:reorganized}
  \left( - \Delta -k_0^2 \right) u_{\text{sc}} = f + k_0^2  \, \chi \, u_{\text{sc}}.
\end{equation}
The right hand side of \eqref{eq:reorganized} is now only non-zero on
$[-L,L]^d$, since both $f$ and $\chi$ are only non-zero
there. Furthermore, since we have solved \eqref{Helmholtz} we also
know $u_\text{sc}$ on $[-L,L]^d$.  So the full right hand side of
\eqref{eq:reorganized} is known. The remaining left hand side of
\eqref{eq:reorganized} is now a Helmholtz equation with a constant
wave number $k_0^2$.  For this equation the Greens function is known
analytically. 
\begin{equation}
  \begin{aligned}
    u_{\text{sc}}(\mathbf{x}) &= \int  G(\mathbf{x},\mathbf{y}) \left( f(\mathbf{y}) +   k_0^2 \,\chi(\mathbf{y}) u_\text{sc}(\mathbf{y})\right) d\mathbf{y} \\
         &=   \int\limits_{[-L,L]^d} G(\mathbf{x},\mathbf{y}) \left( f(\mathbf{y}) +   k_0^2\, \chi(\mathbf{y}) u_\text{sc} (\mathbf{y}) \right) d\mathbf{y} \quad \forall x \in \mathbb{R}^d,
\end{aligned}
\end{equation}
where $f$ and $\chi$ are limited to $[-L,L]^d$ thus we can truncate the
integral to the box $[-L,L]^d$.

This methodology was successfully applied to calculate challenging
break up problems, for example \cite{rescigno1999collisional}.

\subsection{Single ionization versus double ionization}
Let us discuss the qualitative behaviour of the solution for single
and double ionization. To illustrate the behaviour, we truncate the
partial wave expansion, \eqref{eq:partialwaves}, to the first
term. This is known as a $s$-wave expansion. The 6D wave function is
then approximated as
\begin{equation}
  u(\mathbf{r}_1,\mathbf{r}_2) \approx  u(\rho_1,\rho_2)  Y_{00}(\theta_1,\varphi_1) Y_{00}(\theta_2,\varphi_2).
\end{equation}
The radial wave, $u(\rho_1,\rho_2)$, then fits a 2D Helmholtz equation
\begin{equation}\label{eq:radialequation} 
 \left(- \frac{1}{2} \frac{d^2}{d\rho_1^2} - \frac{1}{2} \frac{d^2}{d\rho_2^2} + V_1(\rho_1) + V_2(\rho_2) + V_{12}(\rho_1,\rho_2) - E \right) u(\rho_1,\rho_2) = f(\rho_1,\rho_2), \quad \forall \rho_1,\rho_2 \in [0,\infty[,
\end{equation}
where $V_1(\rho_1)$ and $V_2(\rho_2)$ represents the one-particle
potentials and $V_{12}(\rho_1,\rho_2)$ the two-particle repulsion.
This model is known as a $s$-wave or Temkin-Poet model
\cite{temkin1962nonadiabatic, poet1978exact}.

Before the photo-ionization, the atom is in a two-particle ground
state. In this $s$-wave model, it is the eigenstate of
\begin{equation}
 \left(- \frac{1}{2} \frac{d^2}{d\rho_1^2} - \frac{1}{2} \frac{d^2}{d\rho_2^2} +
 V_1(\rho_1) + V_2(\rho_2) + V_{12}(\rho_1,\rho_2)
 \right)\phi_0(\rho_1,\rho_2) = \lambda_0 \phi_0(\rho_1,\rho_2).
\end{equation}
 with the lowest energy. Simultaniously,
there are one-particle states that are eigenstates of
\begin{equation}\label{eq:eigenstate_H1}
  \left(-\frac{1}{2} \frac{d^2}{d\rho_1^2} + V_1(\rho_1) \right) \phi_i(\rho_1)  =  \mu_i \phi_i(\rho_1)
\end{equation}
and
\begin{equation}
  \label{eq:eigenstate_H2}
  \left(- \frac{1}{2} \frac{d^2}{d\rho_2^2} + V_2(\rho_2) \right) \varphi_i(\rho_2)  =  \nu_i \varphi_i(\rho_2).
\end{equation}

\begin{figure}
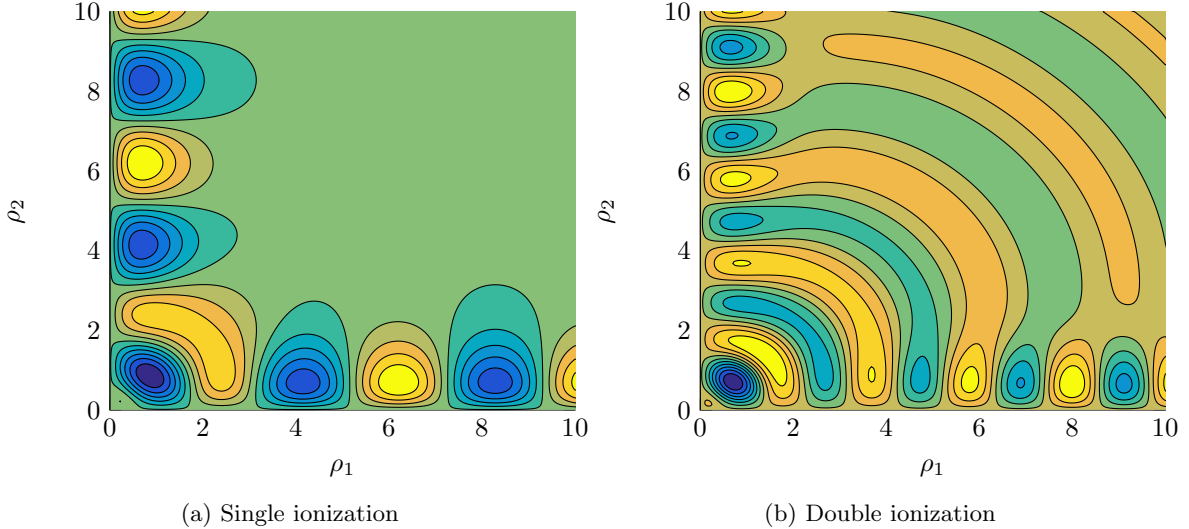

	\centering
	\begin{subfigure}{.45\textwidth}
		\centering
		\input{singleionization2d-contourf-fullrank.tikz}
		\caption{Single ionization}
		\label{fig:singleionization}
	\end{subfigure}
	\begin{subfigure}{.45\textwidth}
		\centering
		\input{doubleionization2d-contourf-fullrank.tikz}
		\caption{Double ionization}
		\label{fig:doubleionization}
	\end{subfigure}  \caption{Left: When the energy $E = \hbar \nu + \lambda_0 <0$, there is
    only single ionization. The solution is then localized along the
    edges, where the solution is a combination of an outgoing wave in
    the $\rho_1$ and a bound state in $\rho_2$, or vice-versa. Right: For
    energy $E>0$, there is,  in addition to single ionization with
    solution localized along the edges, a double ionization wave where
    both coordinates can become large.
  \label{fig:single_vs_double}}
\end{figure}

When \eqref{eq:radialequation} is solved with the energy $E = \hbar
\nu + \lambda<0$, there is only single ionization. Only one of the two
coordinates $\rho_1$ or $\rho_2$ can become large and the solution,
as can be seen in Figure~\ref{fig:singleionization}, is
localized along both axis. The solution is a product of an outgoing
wave in one coordinate and a bound state in the other
coordinate. For example, along the $\rho_2$-axis, the solution is
described by $A_i(\rho_2)\phi_i(\rho_1)$, where $A_i(\rho_2)$ is a
one-dimensional outgoing wave, with an energy $E-\mu_i$ and
$\phi_i(\rho_1)$ is a bound state of \eqref{eq:eigenstate_H1} in the
first coordinate with energy $\mu_i$.  Similary, there is a wave,
along the $\rho_1$ axis, that is an outgoing wave of the form
$B_l(\rho_1) \varphi_l(\rho_2)$, with a scattering wave in the first coordinate, $\rho_1$,
and a bound state in the coordinate $\rho_2$, solution of
\eqref{eq:eigenstate_H2}.

When \eqref{eq:radialequation} is solved with energy $E = \hbar \nu +
\lambda \ge 0$ there is also double ionization and both coordinates
$\rho_1$ and $\rho_2$ can become large. We see, in
Figure~\ref{fig:doubleionization}, a (spherical) wave in the middle of
the domain, where both coordinates can be become large. To describe
this solution the full coordinate space is necessary.  Note that
these solutions still show single ionization along the axes. Even for
$E> 0$, one particle can take away all the energy and leave the other
particle as a bound state.

\begin{figure}
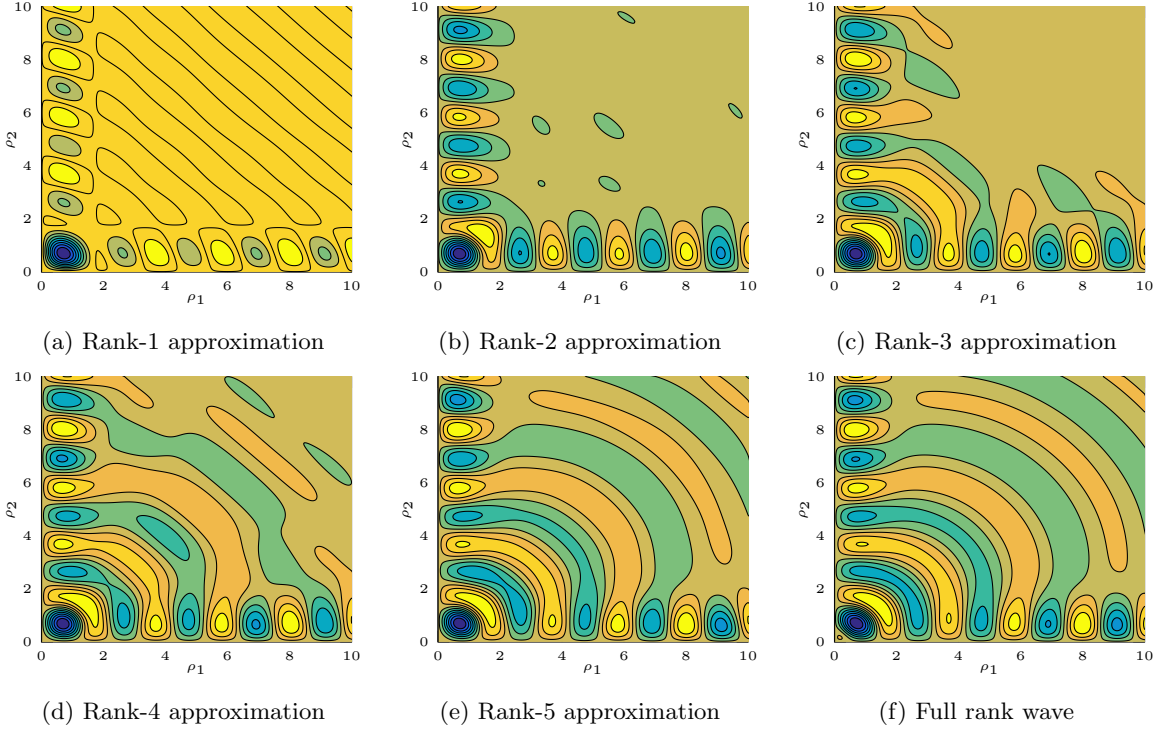

	\centering
	\begin{subfigure}{.3\textwidth}
		\centering
		\input{wave2d-contourf-rank=1.tikz}
		\caption{Rank-1 approximation}
	\end{subfigure}
	\begin{subfigure}{.3\textwidth}
		\centering
		\input{wave2d-contourf-rank=2.tikz}
		\caption{Rank-2 approximation}
	\end{subfigure}
	\begin{subfigure}{.3\textwidth}
		\centering
		\input{wave2d-contourf-rank=3.tikz}
		\caption{Rank-3 approximation}
	\end{subfigure}
	\\~\\
	\begin{subfigure}{.3\textwidth}
		\centering
		\input{wave2d-contourf-rank=4.tikz}
		\caption{Rank-4 approximation}
	\end{subfigure}
	\begin{subfigure}{.3\textwidth}
		\centering
		\input{wave2d-contourf-rank=5.tikz}
		\caption{Rank-5 approximation}
	\end{subfigure}
	\begin{subfigure}{.3\textwidth}
		\centering
		\input{wave2d-contourf-fullrank.tikz}
		\caption{Full rank wave}
	\end{subfigure}
	\caption{Contour plots of the double ionization wave function (bottom, right) and low-rank approximations for increasing rank.}
\end{figure}

\subsection{Coupled channel model for single ionization waves}
\label{sec:singleionization}
In this section, we write the single ionization solution as a low-rank
decomposition and derive the equations for the low-rank
components. When there is only single ionization, the total wave can
be written as
\begin{equation}\label{expansion}
  u(\rho_1,\rho_2) = \sum_{m=1}^{M} \phi_m(\rho_1) A_m(\rho_2)  +  \sum_{l=1}^L  B_l(\rho_1)\varphi_l(\rho_2),
\end{equation}
where $\phi_m(\rho_1)$ and $\varphi_l(\rho_2)$ are the  bound state eigenstates,
defined in \eqref{eq:eigenstate_H1} and \eqref{eq:eigenstate_H2}.  The
first term is localized along the $\rho_2$-axis, the second term is
localized along the $\rho_1$-axis with $\mu_i < 0$ and $\nu_i<0$.

As discussed in \cite{cools2016fast}, this expansion is not unique. We
can add multiples of $\gamma_i \varphi_i(\rho_2)$ to $A_i(\rho_2)$ and
simultaneously subtract $\gamma_l \phi_l(\rho_1)$ from $B_l(\rho_1)$
without contaminating the result.  Indeed, for any choice of $\gamma_i
\in \mathbb{C}$ and $L=M$ holds that
\begin{equation}
u(\rho_1,\rho_2) = \sum_{m=1}^{M} \phi_m(\rho_1) \left(A_m(\rho_2) + \gamma_m \varphi_m(\rho_2) \right)  +  \sum_{l=1}^L  \left(B_l(\rho_1)-\gamma_i \phi_i(\rho_1))\varphi_l(\rho_2) \right) = u(\rho_1,\rho_2).
\end{equation}
To make the expansion unique, \cite{cools2016fast} chooses to select $A_i
\perp \varphi_j$ when $j \ge i$ and $B_i\perp \phi_m$ when $i \ge j$.

In this paper, we choose to make the functions in the set $\{\phi_{i
  \in\{1,\ldots,m\}}, B_{l \in \{1,\ldots, L \}}\}$ orthogonal.  We
also assume that $V_{12}(\rho_1,\rho_2) \approx \sum_{i=1}^L\sum_{l}^M
\phi_i(\rho_1) \varphi(\rho_2) \int
\phi_i^*(\rho_1)\varphi^*_l(\rho_2) V_{12}(\rho_1,\rho_2) d\rho_1
d\rho_2.$

Given a function $f(\rho_1,\rho_2)$, a right hand side, we can now
derive the equations for $A_m$ and $B_l$. When we insert the low-rank
decomposition of the expansion \eqref{expansion}, in the 2D Hamiltonian
\eqref{eq:radialequation}, multiply with $\phi_i$ and integrate over
$\rho_1$ we find
\begin{equation}
  (H_2 + \mu_j -E) A_j(\rho_2) + \sum_{l=1}^M V_{ij}(\rho_2) A_l(\rho_2) =  \int_0^\infty \phi^*_j(\rho_1) f(\rho_1,\rho_2) d\rho_1, \quad \text{for} \quad j=1, \ldots, M \quad \text{and}\quad  \forall \rho_2 \in [0,L],
\end{equation}
with
\begin{equation}
 V_{ij}(\rho_2)=  \int_0^\infty \phi^*_j(\rho_1) V_{12}(\rho_1,\rho_2) \phi_l(\rho_1) d\rho_1 .
\end{equation}
We have used that $\phi_i \perp B_l$ to eliminate the second term in
the expansion \eqref{expansion}.

Similarly, for $B_l$, we find
\begin{equation}
  \begin{aligned}
  (H_1 + \nu_l -E) B_l(\rho_1) + \sum_{m=1}^L W_{l}(\rho_1) B_k(\rho_1) =  \int_0^\infty\!\!\!\! \varphi_l(\rho_2) \left(f(\rho_1,\rho_2)-\sum_i^M \phi_i \phi^* f(\rho_1,\rho_2) \right) d\rho_2 \quad\\
    \text{for} \quad i=1,\ldots, M, \quad \forall \rho_1 \in [0,\infty[
  \end{aligned}
\end{equation}
where
\begin{equation}
  W_{lk}(\rho_1):=\left(\int_0^\infty \!\!\!\!  \varphi^*_l(\rho_2) V_{12}(\rho_1,\rho_2) \varphi_k(\rho_2) d\rho_2 \right).
\end{equation}

\section{Low-rank matrix representation of a 2D wave function that includes both single and double ionization}
\label{sec:lowrank}
\subsection{Low rank of the double ionization solution}
We now discuss the main result of the paper. We will derive a
coupled channel equation that gives a low-rank approximation for the
double ionization wave function, as shown in Figure~\ref{fig:doubleionization}.
   
In section \ref{sec:singleionization} we have shown
that the single ionization wave can be represented by a low-rank
decomposition. In this section, we show that also the double
ionization wave can be written as a similar low-rank decomposition.

We first illustrate that the solution of a 2D driven Schr\"odinger equation
that contains both single and double ionization, it is a solution of
\eqref{eq:radialequation} with $E>0$, can be represented by a similar low-rank
decomposition.

\begin{figure}
	\begin{tabular}{ll}
		\includegraphics[width=0.47\textwidth]{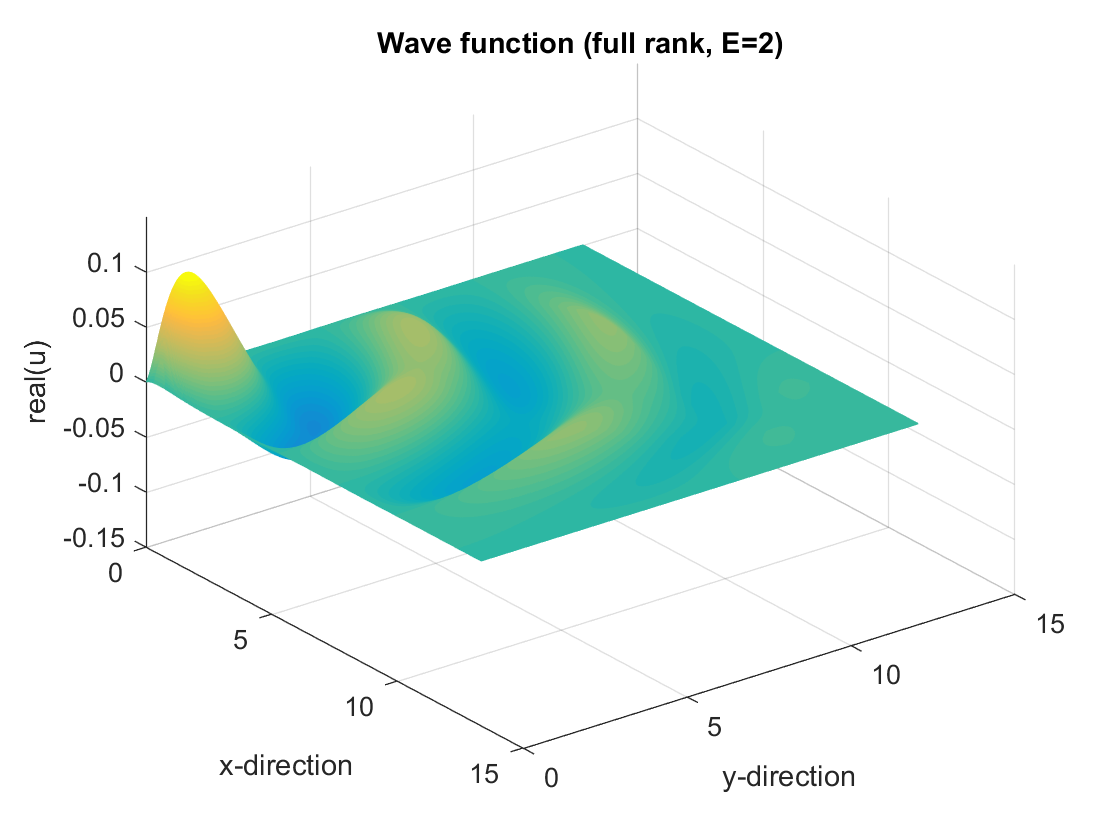} &     \includegraphics[width=0.47\textwidth]{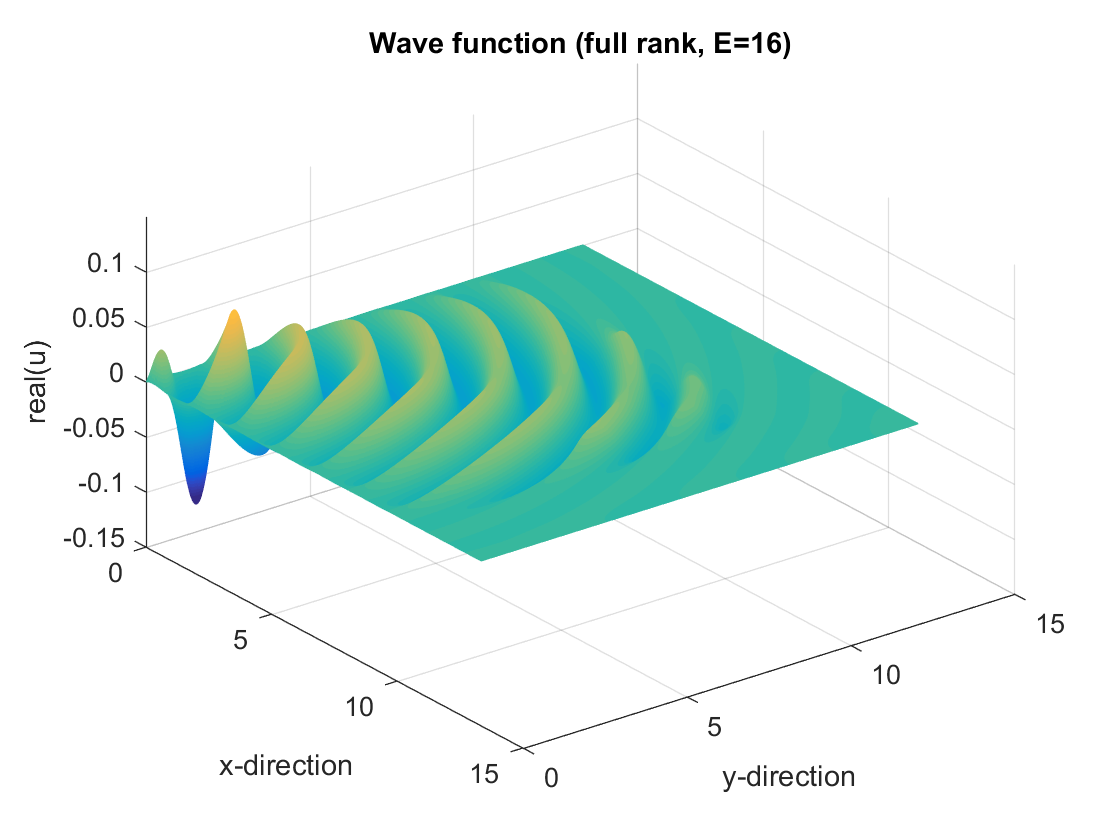}\\
		\input{wave-svd-M=1000-E=2.tikz} & \input{wave-svd-M=1000-E=16.tikz}
	\end{tabular}
	\caption{The top figures show the wave function for two different
		energies (i.e. $E=2$ and $E=16$). The singular values of the
		matrix-representation of the discretized functions are shown in
		the bottom figures.  we show the solution of a 2D Helmholtz
		equation with a space-dependent wave number $k^2(x,y)$ given
		by $k^2(x,y) = E - e^{-|x-y|}$. The right hand side $f(x,y)$ on
		the finite domain $[0, b[^2$ is given by $f(x,y) = -
		e^{-x^2-y^2}$.  Finite difference discretization is done on a
		uniform mesh with $M=1000$ interior mesh points per
		dimension. At the boundaries $x=b$ and $y=b$ the domain is
		extended with exterior complex scaling under an angle
		$\frac{\pi}{6}$ where $33\%$ additional discretization points
		are added, so $n = 1333$. We have $b=10$.}
	\label{fig:wavefunction}
\end{figure}

In Figure~\ref{fig:wavefunction} we solve the Helmholtz equation with a
space-dependent wave number, $k(x,y)$, in the first quadrant where $x
\ge 0$ and $y \ge 0$. The equation is
\begin{equation}\label{driven}
\left(  - \Delta_2 - k^2(x,y) \right) u_{sc}(x,y) = f(x,y),
\end{equation}
where $\Delta_2$ is the 2D Laplacian and the solution $u_{sc}$
satisfies homogeneous boundary conditions $u_{sc}(x,0)=0$ for all
$x\ge 0$ and $u_{sc}(0,y)=0$ for all $y \ge 0$. On the other
boundaries we have outgoing boundary conditions.

The right hand side $f(x,y)$ has a support that is limited to $[0,b]^2
\subset [0,L]^2 \subset \mathbb{R}_+^2$, i.e. $f(x,y)=0$, for all
$x\ge b$ or $y\ge b$.

The wave number $k(x,y)$ can be split in a constant part, $k_0^2$ and
variable part $\chi(x,y)$. The variable part is also only non-zero on
  $[0,b[^2$
\begin{equation}
  k^2(x,y) = \begin{cases}
    k_0^2 \left(1  + \chi(x,y) \right)  \quad &\text{if} \quad  x<b \quad \text{and} \quad y<b,\\
    k_0^2   \quad &\text{if} \quad  x\ge b \quad \text{or} \quad y\ge b.\\
  \end{cases}
\end{equation}
We extend the domain with exterior complex scaling (ECS) absorbing
boundary condition \cite{mccurdy2004solving}.

The wave function $u_{sc}$ is discretized on the two dimensional mesh and can be
represented as a matrix $\mv{A}\in \mathbb{C}^{n \times n}$. We can
compute the singular decomposition of this matrix, $\mv{A} = \mv{U}
\mv{\Sigma} \mv{V}^\H$ where $\mv{U},\mv{\Sigma}, \mv{V} \in
\mathbb{C}^{n \times n}$ where $\mv{U}^\H \mv{U} = \mv{I}$, $\mv{V}^\H
\mv{V} = \mv{I}$ and $\mv{\Sigma}$ is a diagonal matrix with the
singular value $\sigma_i$ on the diagonal.

The results are shown in Figure~\ref{fig:wavefunction} and show that
the singular values rapidly decrease. Thus the wave function can
efficiently be approximated by a truncated representation,
\begin{equation}\label{eq:truncated}
  \mv{A} \approx \sum_{i=1}^r \mv{u}_i \sigma_i \mv{v}^\H_i, 
\end{equation}
where $\mv{u}_i \in \mathbb{C}^n$ are columns of $\mv{U}$ and $\mv{v}_i \in \mathbb{C}^n$ are rows from $\mv{V}$ and $\sigma_i$ are largest $r$ singular values. Thus $\mv{A}$ is approximated by it low-rank representation with rank $r$.
This truncated decomposition drops all contributions with $\sigma_i<\tau$ below a threshold $\tau$, for example the expected discretisation error.

\begin{figure}
	\begin{center}
		\input{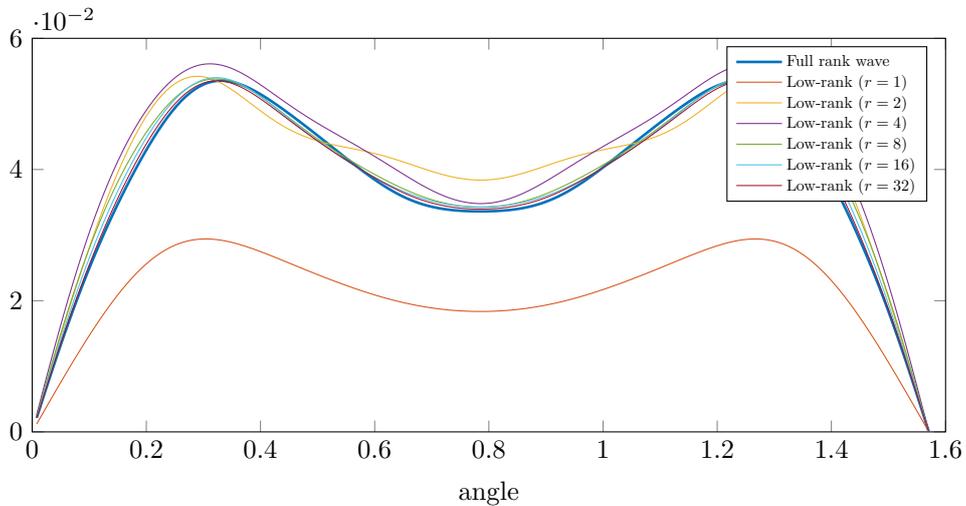}
	\end{center}
	\caption{The cross section based on different low-rank approximations of the wave function.}
	\label{fig:crosssection}
\end{figure}

Figure \ref{fig:crosssection} illustrates that a low-rank approximation to the wave function is sufficient
calculate an accurate approximation of the cross section.

\subsection{Determining the low-rank components directly}
\label{sec:lowrankcomponents}
In the example of the previous section, we have first calculated a
matrix representation of the solution, $\mv{A} \in \mathbb{C}^{n \times n}$,
and then approximated it by low-rank components.  The aim is now to
develop a method that calculates, directly, these components without
first calculating the full solution $\mv{A}$. This approach avoids expensive
calculations.

We start from the 2D Helmholtz equation as given in \eqref{driven}.
In matrix form this is given by
\begin{equation} \label{eq:matrixform}
 - \mv{D}_{xx} \mv{A} - \mv{A} \mv{D}^\T_{yy} - \mv{K} \circ \mv{A} = \mv{F},
\end{equation}
where $\mv{D}_{xx} \in \mathbb{C}^{n \times n}$ and $\mv{D}_{yy} \in
\mathbb{C}^{n \times n}$ are sparse matrices that represent the discretization
of the second derivatives, $\mv{K}$ is the matrix that represents the
space-dependent wave number, $k^2(x,y)$, on the grid and $\mv{A} \in
\mathbb{C}^{n\times n}$ is the matrix that describes the unknown
partial wave. The right hand side $\mv{F} \in \mathbb{C}^{n \times
  n}$ is given. The Hadamard product, $\circ$, multiplies the
matrices point wise, element by element.

We now make the approximation  $\mv{A} \approx \mv{U}\mv{V}^\H$, with low-rank matrices
$\mv{U} \in \mathbb{C}^{n \times r}$ and $\mv{V} \in \mathbb{C}^{n
  \times r}$ where $r \ll n$ and write
\begin{equation} \label{eq:matrixequation}
  - \mv{D}_{xx} \mv{U} \mv{V}^\H - \mv{U} \mv{V}^\H \mv{D}_{yy} - \mv{K} \circ \mv{U} \mv{V}^\H = \mv{F}.
\end{equation}
We start with a guess for $\mv{V} \in \mathbb{C}^{n \times r}$ with
orthogonal columns such that $\mv{V}^\H \mv{V}=\mv{I}_r$. We can now 
multiply, \eqref{eq:matrixequation}, from the right,  by $\mv{V}$ and obtain
\begin{equation}\label{eq:2dvarwavenr}
 - \mv{D}_{xx} \mv{U} \mv{V}^\H \mv{V} - \mv{U} \mv{V}^\H \mv{D}_{yy} \mv{V}  -  \left(\mv{K} \circ \mv{U}\mv{V}^\H\right) \mv{V} = \mv{F}\mv{V}. 
\end{equation}
where $\mv{U} \in \mathbb{C}^{n \times r}$ is now the remaining unknown.  We use
the vectorizing identities $\mvec{\mv{A} \circ \mv{B} } = \mvec{\mv{A}
}\circ \mvec{\mv{B}}$ for $\mv{A}, \mv{B}~\in~\mathbb{C}^{l \times p}$
and $\mvec{\mv{A}\mv{B}} = \left(\mv{I}_m \otimes \mv{A} \right) \mvec{\mv{B}} =
\left(\mv{B}^\T \otimes \mv{I}_k \right) \mvec{\mv{A}}$, for $\mv{A} \in
\mathbb{C}^{k \times l}$ and $\mv{B} \in \mathbb{C}^{l \times m}$ and
obtain
\begin{equation}
 - ( \mv{I}_r \otimes \mv{D}_{xx}) \mvec{\mv{U}}  - \left( (\mv{V}^\H \mv{D}_{yy} \mv{V})^\T \otimes \mv{I}_r  \right)\mvec{\mv{U}} - \mvec{\left(\mv{K} \circ \mv{U}\mv{V}^\H \right) \mv{V} }  = \mvec{\mv{F}\mv{V}}.
\end{equation}
We can simplify the last term of the left hand side and write
\begin{equation}\label{eq:vecThirdTerm}
  \begin{aligned}
    \mvec{  \left(\mv{K} \circ \mv{U}\mv{V}^\H \right) \mv{V} } &= \left(\mv{V}^\T \otimes \mv{I} \right) \mvec{\mv{K} \circ \mv{U} \mv{V}^\H},\\
    &= \left(\mv{V}^\T \otimes \mv{I} \right) \left( \mvec{\mv{K}} \circ \mvec{\mv{U} \mv{V}^\H} \right),\\
    &= \left(\mv{V}^\T \otimes \mv{I} \right) \diag{\mvec{\mv{K}}}  \mvec{\mv{U} \mv{V}^\H},\\
    &= \left(\mv{V}^\T \otimes \mv{I} \right) \diag{\mvec{\mv{K}}}  \left(\left(\mv{V}^\H\right)^\T \otimes \mv{I} \right) \mvec{\mv{U}}.\\
  \end{aligned}
\end{equation}
This results in
\begin{equation}\label{eq:first}
\left[ - ( \mv{I} \otimes \mv{D}_{xx})  - \left( \left(\mv{V}^\H \mv{D}_{yy} \mv{V}\right)^\T \otimes \mv{I}  \right) -  \left(\mv{V}^\T \otimes \mv{I} \right) \diag{\mvec{\mv{K}}}  \left(\left(\mv{V}^\H\right)^\T \otimes \mv{I} \right)\right] \mvec{\mv{U}} = \mvec{\mv{F}\mv{V}}.
\end{equation}
This is a linear system for the remaining unknown columns of the matrix $\mv{U} \in \mathbb{C}^{n \times r}$.

In \eqref{eq:truncated}, we have approximated $\mv{A}$ as
$\mv{U}\mv{\Sigma}\mv{V}^\H$ where $\mv{U},\mv{V} \in \mathbb{C}^{n
  \times r}$and $\mv{\Sigma} \in \mathbb{C}^{r \times r}$ are
truncated matrices. With an orthogonal guess for $\mv{V}$, we solve
for a $\mv{U}$ in \eqref{eq:first}.  Since we approximate $\mv{A}$
now by the product $\mv{U}\mv{V}^\H$, hence $\mv{U}$, solution of
\eqref{eq:first}, includes the diagonal matrix with singular values.

We can now do a QR decomposition of $\mv{U}$ to arrive at a guess for
the orthogonal matrix $\mv{U}$.

The next step is to improve the guess for $\mv{V}$ in a similar way. The
equation \eqref{eq:matrixequation} becomes, when we multiply 
from the left by $\mv{U}^\H$,
\begin{equation}
 - \mv{U}^\H \mv{D}_{xx} \mv{U}\mv{V}^\H - \mv{U}^\H\mv{U}\mv{V}^\H \mv{D}_{yy} - \mv{U}^\H \left(\mv{K} \circ \mv{U}\mv{V}^\H\right) = \mv{U}^\H \mv{F}.
\end{equation}
Using the vectorizing identities this results in
\begin{equation}\label{eq:second}
\left[ - \left(\mv{I} \otimes  \mv{U}^\H \mv{D}_{xx}\mv{U}\right) - \left(\mv{D}_{yy}^\T \otimes \mv{I} \right) -  \left( \mv{I} \otimes \mv{U}^\H\right) \diag{\mvec{\mv{K}}}  \left( \mv{I} \otimes \mv{U} \right) \right] \mvec{\mv{V}^\H}   = \mvec{\mv{U}^\H \mv{F}},
\end{equation}
where we use that 
\begin{equation}
  \begin{aligned}
    \mv{U}^\H \left(  \mv{K} \circ \mv{U}\mv{V}^\H\right)  & = \left( \mv{I} \otimes \mv{U}^\H \right) \mvec{\mv{K}} \circ \mvec{\mv{U}\mv{V}^\H} \\
      & = \left( \mv{I} \otimes \mv{U}^\H \right) \diag{\mvec{\mv{K}}}  \left( \mv{I} \otimes \mv{U} \right) \mvec{\mv{V}^\H}. \\
  \end{aligned}
\end{equation}

Combining the equations \eqref{eq:first} and \eqref{eq:second} we can
now propose an algorithm that updates $\mv{U}$ and $\mv{V}$ in an
alternating way. The steps are described in the
Algorithm~\ref{alg:2d}.

\begin{algorithm}[t]
  \SetAlgoLined
  Choose $\mv{V} \in \mathbb{C}^{n \times r}$ as initial guess\; 
  $[\mv{V},\mv{R}] = \qr{\mv{V},0}$\; 
  \While{not converged}{
    Solve $\left[ - \left( \mv{I} \otimes \mv{D}_{xx}\right)  - \left( \left(\mv{V}^\H \mv{D}_{yy} \mv{V}\right)^\T \otimes \mv{I}  \right) -  \left(\mv{V}^\T \otimes \mv{I} \right) \diag{\mvec{\mv{K}}}  \left(\left(\mv{V}^\H\right)^\T \otimes \mv{I} \right)\right] \mvec{\mv{U}} = \mvec{\mv{F}\mv{V}}$\; 
    $[\mv{U},\mv{R}] = \qr{\mv{U},0}$\;
    Solve $\left[ - \left(\mv{I} \otimes  \mv{U}^\H \mv{D}_{xx} \mv{U} \right) - \left(\mv{D}_{yy}^\T \otimes \mv{I}\right)  -  \left( \mv{I} \otimes \mv{U}^\H\right) \diag{\mvec{\mv{K}}}  \left( \mv{I} \otimes \mv{U} \right) \right] \mvec{\mv{V}^\H}   = \mvec{\mv{U}^\H \mv{F}}$\;
    $[\mv{V},\mv{R}] = \qr{\mv{V},0}$\;
  }
  $\mv{A} = \mv{U} \mv{R}^\H \mv{V}^\H$\;
  \caption{Solve for the low-rank matrix decomposition of the solution $\mv{A} \approx \mv{U}\mv{V}^\H$ of a 2D Helmholtz problem with space-dependent wave number.}
  \label{alg:2d}
\end{algorithm}

\subsection{Comparison between coupled channel and a low-rank decomposition}
We now compare the coupled channel approach from section
\ref{sec:singleionization} with the low-rank approach from the
previous section, Sec.~\ref{sec:lowrankcomponents}.  In the coupled channel
calculation we use the eigenfunctions of one-particle subsystems as a
basis for the expansion. After integration over one of the coordinates
this results in a coupled set of one-dimensional equations. Let us
illustrate that equation \eqref{eq:first} and \eqref{eq:second} reduce
to a the coupled channel equations when we choose the eigenfunctions,
from \eqref{eq:eigenstate_H1} and \eqref{eq:eigenstate_H2}, as columns
for $\mv{U}$.

Let use take a look at the equation \eqref{eq:second} and assume that
$\mv{K} = E\mv{I} -V_1(\mv{x})-V_2(\mv{y})-V_{12}(\mv{x},\mv{y})$ and
that the columns of $\mv{U}$ are the eigenstates of
$-\mv{D}_{xx}+V_1(\mv{x})$.  We can then write
\begin{equation*}
    - \left( \mv{I} \otimes \mv{U}^\H\right) \diag{\mvec{\mv{K}}}  \left( \mv{I} \otimes \mv{U} \right)
    =  - E\,  \mv{I} \otimes \mv{I} + V_2(\mv{y}) \otimes \mv{I}  + \mv{I} \otimes \mv{U}^\H V_1(\mv{x}) \mv{U} + \left( \mv{I} \otimes \mv{U}^\H\right) \diag{ \mvec{V_{12}(\mv{x},\mv{y})}}  \left( \mv{I} \otimes \mv{U} \right).
\end{equation*}
Then equation \eqref{eq:second} becomes
\begin{equation*}
\left[  \mv{I} \otimes  \mv{U}^\H \left(-\mv{D}_{xx} + V_1(\mv{x})\right) \mv{U}  +  \left(-\mv{D}_{yy}^\T + V_2(\mv{y}) -E \mv{I}\right) \otimes \mv{I}  -   \left( \mv{I} \otimes \mv{U}^\H\right) \diag{ \mvec{V_{12}(\mv{x},\mv{y})}} \left( \mv{I} \otimes \mv{U} \right)  \right] \mvec{\mv{V}^\H}   = \mvec{\mv{U}^\H \mv{F}}.
\end{equation*}
When we use that the columns of $\mv{U}$ are eigenfunctions of
$-\mv{D}_{xx}+V_1(\mv{x})$ with eigenvalues $\mu_i$, equation
\eqref{eq:second} becomes
\begin{equation}
\left[  \mv{I} \otimes \diag{\mv{\mu}} + \left(-\mv{D}_{yy}^\T + V_2(\mv{y}) -E \, \mv{I}\right) \otimes \mv{I} -  \left( \mv{I} \otimes \mv{U}^\H\right) \diag{\mvec{V_{12}(\mv{x}, \mv{y})}}  \left( \mv{I} \otimes \mv{U} \right)  \right] \mvec{\mv{V}^\H}   = \mvec{\mv{U}^\H \mv{F}}.
\end{equation}
The term $\left( \mv{I} \otimes \mv{U}^\H\right)
\diag{\mvec{V_{12}(\mv{x}, \mv{y})}} \left( \mv{I} \otimes \mv{U}
\right)$ couples the columns of $\mv{V}$. It should be interpreted as
a discretized version of $\int \phi^*_i(x)V(x,y)\phi_j(x) dx$, where
we integrate over one of the coordinates. The columns of $\mv{U}$ are
the $\phi_i$ represented on a integration grid.

However, in general, the columns are $\mv{U}$ are not eigenfunctions of the
operator. The term $\mv{I} \otimes \diag{\mv{\mu}}$ then becomes a matrix that
also couples the different components of $\mv{V}$.

In short, each iteration we are solving a generalized coupled channel
equation.
\subsection{Convergence with projection operators}
\label{sec:projection2d}
We will now write both linear systems, \eqref{eq:first} for
$\mvec{\mv{U}}$, and, \eqref{eq:second} for $\mvec{\mv{V}}$, as
projection operators applied to the residual of the matrix equation,
\eqref{eq:matrixform}.

We denote by $\mv{L}$ the discretized 2D Helmholtz operator on the full grid
\begin{equation} \label{eq:linearoperator}
	\mv{L} = \left( \mv{I} \otimes (-\mv{D}_{xx}) \right) + \left( (-\mv{D}_{yy}) \otimes \mv{I} \right) - \diag{\mvec{\mv{K}}} \left(\mv{I} \otimes \mv{I} \right).
\end{equation}
We can now explicitly write equation \eqref{eq:first} in terms of projections and this linear operator:
\begin{equation} \label{eq:solutionU}
	\left(\mv{V}^\T \otimes \mv{I} \right) \mv{L} \left(\conj{\mv{V}} \otimes \mv{I} \right) \mvec{\mv{U}} = \left(\mv{V}^\T \otimes \mv{I} \right)\mvec{\mv{F}}.
\end{equation}

The residual matrix, $\mv{R}$, is given by
\begin{equation}
	\begin{split}
	  \mv{R} &= \mv{F} - (-\mv{D}_{xx}) \mv{U} \mv{V}^\H - \mv{U} \mv{V}^\H (-\mv{D}_{yy}) + \mv{K} \circ \left(\mv{U}\mv{V}^\H\right).
        \end{split}
\end{equation}
In vector form this reads, using that $\mv{U}$ is a solution of equation \eqref{eq:solutionU},
\begin{equation}
  \begin{split}
		\mvec{\mv{R}} &= \mvec{\mv{F}} - \left( \conj{\mv{V}} \otimes (- \mv{D}_{xx}) + (-\mv{D}_{yy}^\T) \conj{\mv{V}} \otimes \mv{I} - \diag{\mvec{\mv{K}}} \left( \conj{\mv{V}} \otimes \mv{I} \right)  \right) \mvec{\mv{U}}, \\
		&= \left(\mv{I} - \left( \conj{\mv{V}} \otimes (-\mv{D}_{xx}) + (-\mv{D}_{yy}^\T) \conj{\mv{V}} \otimes \mv{I} - \diag{\mvec{\mv{K}}} \left( \conj{\mv{V}} \otimes \mv{I} \right)  \right) \left[ \left(\mv{V}^\T \otimes \mv{I} \right) \mv{L} \left(\conj{\mv{V}} \otimes \mv{I} \right) \right]^{-1} \left(\mv{V}^\T \otimes \mv{I} \right) \right) \mvec{\mv{F}} \\
		&= \left(\mv{I} - \mv{L} \left( \conj{\mv{V}} \otimes \mv{I} \right) \left[ \left(\mv{V}^\T \otimes \mv{I} \right) \mv{L} \left(\conj{\mv{V}} \otimes \mv{I} \right) \right]^{-1} \left(\mv{V}^\T \otimes \mv{I} \right) \right) \mvec{\mv{F}} \\
		&= P_{\mv{V}}\mvec{\mv{F}},
		\end{split}
\end{equation}
where $P_{\mv{V}}$ is given by
\begin{equation}\label{eq:projectorVarOnV}
	\begin{split}
		P_{\mv{V}} &:= \mv{I} - \mv{L} \left( \conj{\mv{V}} \otimes \mv{I} \right) \left[ \left(\mv{V}^\T \otimes \mv{I} \right) \mv{L} \left(\conj{\mv{V}} \otimes \mv{I} \right) \right]^{-1} \left(\mv{V}^\T \otimes \mv{I} \right) \\
		&:= \mv{I} - \mv{X}.
	\end{split}
\end{equation}
The operator $P_{\mv{V}}$ is a projection operator. Indeed, observe
that the terms between the two inverses cancel, in the next equation,
against one of the inverse factors:
\begin{equation}
	\begin{split}
		\mv{X}^2 &= \mv{L} \left( \conj{\mv{V}} \otimes \mv{I} \right) \left[ \left(\mv{V}^\T \otimes \mv{I} \right) \mv{L} \left(\conj{\mv{V}} \otimes \mv{I} \right) \right]^{-1} \left(\mv{V}^\T \otimes \mv{I} \right) \mv{L} \left( \conj{\mv{V}} \otimes \mv{I} \right) \left[ \left(\mv{V}^\T \otimes \mv{I} \right) \mv{L} \left(\conj{\mv{V}} \otimes \mv{I} \right) \right]^{-1} \left(\mv{V}^\T \otimes \mv{I} \right) \\
		&= \mv{L} \left( \conj{\mv{V}} \otimes \mv{I} \right) \left[ \left(\mv{V}^\T \otimes \mv{I} \right) \mv{L} \left(\conj{\mv{V}} \otimes \mv{I} \right) \right]^{-1} \left(\mv{V}^\T \otimes \mv{I} \right)\\
		&= \mv{X}.
	\end{split}
\end{equation}
We then have that $P_{\mv{V}}^2=(1-\mv{X})(1-\mv{X}) =
1-2\mv{X}+\mv{X}^2 = 1- \mv{X} = P_{\mv{V}}$.

This projection operator removes all components from the residual matrix
that can be corrected by the subspace spanned by $\mv{V}$. It is
similar as a deflation operator, often used in preconditioning
\cite{gaul2012preconditioned}.

A similar derivation results in a projection operator $P_{\mv{U}}$ for
the update of $\mv{V}$:
\begin{equation}\label{eq:projectorVarOnU}
	\begin{split}
		P_{\mv{U}} &= \mv{I} - \mv{L} \left( \mv{I} \otimes \mv{U} \right) \left[\left(\mv{I} \otimes \mv{U}^\H \right) \mv{L} \left(\mv{I} \otimes \mv{U} \right)\right]^{-1} \left(\mv{I} \otimes \mv{U}^\H \right) \\
		&= \mv{I} - \mv{Y}.
	\end{split}
\end{equation}
This is again a projection operator.

So, Algorithm~\ref{alg:2d} repeatedly projects the residual matrix, $\mv{R}$, on a
subspace. Alternating between a subspace that is orthogonal to the
subspace spanned by the columns of $\mv{V}^{(k)}$ at iteration $k$ and
a subspace that is orthogonal to the columns of $\mv{U}^{(k)}$ at
iteration $k$. The residual matrix $\mv{R}^{(k)}$ after $k$ iterations
is the result of a series of projections
\begin{equation}
\mv{R}^{(k)} = P_{\mv{U}^{(k)}}P_{\mv{V}^{(k)}} \, P_{\mv{U}^{(k-1)}}P_{\mv{V}^{(k)}} \ldots \,P_{\mv{U}^{(0)}}P_{\mv{V}^{(0)}} \,\mv{R}^{(0)}.
\end{equation}
It is similar to the \textit{method of Alternating Projections}
\cite{xu2002method} that goes back to Neumann
\cite{neumann1950functional}, where a solution is projected on two
alternating subspaces resulting in a solution that lies in the
intersection between the two spaces.  However, here the columns for
$\mv{U}^{(k)}$ and $\mv{V}^{(k)}$ are changing each iteration.

However, when the rank of $\mv{U}^{(k)}$ and $\mv{V}^{(k)}$ is 
sufficiently large, the only intersection between the changing 
subspaces  is the $\mv{0}$ matrix. 
So the residual converges to $0$.

\section{Solving for the low-rank tensor approximation of 3D Helmholtz equations}
\label{sec:3dhelmholtz}
Algorithm~\ref{alg:2d}, as introduced in section~\ref{sec:lowrankcomponents}
for two dimensional problems, can be extended to higher dimensions.
We illustrate now this extension for 3D Helmholtz problems. First, we will 
discuss the problem with a constant wave number and then extend the results 
to space-dependent wave numbers.

We solve the Helmholtz equation with a constant or space-dependent wave number,
$k(x,y,z)$, in the first quadrant, where $x\ge 0, y \ge 0$ and $z \ge 0$. The equation is
\begin{equation}\label{eq:driven3d}
	\left(  -\Delta_3 - k^2(x,y,z) \right) u_{sc}(x,y,z) = f(x,y,z),
\end{equation}
where $\Delta_3$ is the 3D Laplacian and the solution $u_{sc}$
satisfies homogeneous boundary conditions $u_{sc}(0,y,z)=0$ for all
$y,z\ge 0$, $u_{sc}(x,0,z)=0$ for all $x,z \ge 0$ and $u_{sc}(x,y,0)=0$ 
for all $x,y \ge 0$. On the other boundaries we have outgoing boundary conditions.

The right hand side $f(x,y,z)$ has a support that is limited to $[0,b]^3
\subset [0,L]^3 \subset \mathbb{R}_+^3$, i.e. $f(x,y,z)=0$, for all
$x \ge b$, $y \ge b$ or $z \ge b$.

The wave number $k(x,y,z)$ can be split in a constant part, $k_0^2$, and
variable part $\chi(x,y,z)$. The variable part is also only non-zero on
$[0,b[^3$
\begin{equation}
	k^2(x,y,z) = \begin{cases}
		k_0^2 \left(1  + \chi(x,y,z) \right)  \quad &\text{if} \quad  x<b \quad \text{and} \quad y<b \quad \text{and} \quad z<b,\\
		k_0^2   \quad &\text{if} \quad  x\ge b \quad \text{or} \quad y\ge b \quad \text{or} \quad z\ge b.\\
	\end{cases}
\end{equation}
We extend the domain with exterior complex scaling (ECS) absorbing
boundary condition \cite{mccurdy2004solving}.

The wave function is discretized on a three
dimensional mesh with $n_1 \times n_2 \times n_3$ unknowns
and can be represented by a Tucker tensor decomposition \cite{tucker1966some}
$\mt{M} \in \C^{n_1 \times n_2 \times n_3}$ with multi-linear rank 
$\mv{r} = (r_1, r_2, r_3)$:
\begin{equation}
	\mt{M} = \mt{G} \times_1 \mv{U}_1 \times_2 \mv{U}_2 \times_3 \mv{U}_3 \in \C^{n_1 \times n_2 \times n_3}.
\end{equation}
Here the tensor $\mt{G} \in \C^{r_1 \times r_2 \times r_3}$ is called
the core tensor and the factor matrices $\mv{U}_i \in \C^{n_i \times r_i}$
have orthonormal columns for $i = 1, 2, 3$.  Here
$r_i$ refers to the rank for each direction and $n_i$ to the number of
mesh points in each direction. So, to store this tensor only one core
tensor and $d$ factor matrices need to be stored, so the
storage costs\footnote{Here, we used the notation $r^d := \prod_{i=1}^{d} r_i$ and $r := \sqrt[d]{r^d}$,
to deal with possible unequal number of discretization points $n_i$ or ranks $r_i$ in different directions for a Tucker tensor.
}
scales $\Oh{r^d + dnr}$.

Let $\mathcal{L}$ be the discretization of
the three dimensional Helmholtz operator as given in
\eqref{eq:driven3d}. Observe that the operator $\mathcal{L}$ can be
written as a sum of Kronecker-products, where the matrix representation
$\mv{L}$ is of the following form 
\begin{equation}\label{eq:operatorLasSumKroneckerProducts}
	\mv{L} = -\mv{I} \otimes \mv{I} \otimes \mv{D}_{xx} - \mv{I} \otimes \mv{D}_{yy} \otimes \mv{I}  - \mv{D}_{zz} \otimes \mv{I} \otimes \mv{I} - \diag{\mvec{\mt{K}}} \mv{I} \otimes \mv{I} \otimes \mv{I}.
\end{equation}
Here $\mv{D}_{xx} \in \mathbb{C}^{n_1 \times n_1}, \mv{D}_{yy} \in \mathbb{C}^{n_2 \times n_2}$ and $\mv{D}_{zz} \in \mathbb{C}^{n_3 \times n_3}$ are sparse matrices that represent the discretization
of the second derivatives and $\mt{K}$ is the tensor that represents the
constant or space-dependent wave number, $k^2(x,y, z)$, discretized on the grid.

\subsection{Helmholtz equation with constant wave number}
First, consider the 3D Helmholtz problem with a constant wave number, so $k^2(x,y,z) \equiv k^2$.
The application of the Helmholtz operator $\mathcal{L}$ on
tensor $\mt{M}$ in Tucker tensor format is given by
\begin{equation}\label{eq:start3d}
	\begin{split}
		\mathcal{L} \mt{M} &= \mt{F} \\
		\mathcal{L} \mt{M} &= -\mt{G} \times_1 \mv{D}_{xx} \mv{U}_1 \times_2 \mv{U}_2 \times_3 \mv{U}_3\\
		&- \mt{G} \times_1 \mv{U}_1 \times_2 \mv{D}_{yy} \mv{U}_2 \times_3 \mv{U}_3\\
		&- \mt{G} \times_1 \mv{U}_1 \times_2 \mv{U}_2 \times_3 \mv{D}_{zz} \mv{U}_3\\
		&- k^2 \mt{G} \times_1 \mv{U}_1 \times_2 \mv{U}_2 \times_3 \mv{U}_3\\
		&= \mt{F}
	\end{split}
\end{equation}
where $\mv{U}_i^\H \mv{U}_i = \mv{I}$ for $i = 1,2,3$ and $\mt{F}$ is a tensor representation of the right hand side function $f$ discretized on the grid.

\subsubsection{Solving for product of basis functions and core terms (version 1)}\label{sec:version1}
Similar to the two dimensional case, we can derive equations to
iteratively solve for the factors $\mv{U}_1$, $\mv{U}_2$ and
$\mv{U}_3$. To derive the equations for $\mv{U}_1$ we start from
\eqref{eq:start3d} and multiply with $\mv{U}_2$ and $\mv{U}_3$ in the
second and third direction, respectively:
\begin{equation*}
	\mathcal{L}\mt{M} \times_2 \mv{U}_2^\H \times_3 \mv{U}_3^\H = \mt{F} \times_2 \mv{U}_2^\H \times_3 \mv{U}_3^\H.
\end{equation*}
For a review of tensors, tensor decompositions and tensor 
operations like this tensor-times-matrix product denoted by 
the symbol $\times_i$ we refer to \cite{kolda2009tensor}.

Using explicitly the Tucker tensor representation for $\mt{M}$ and that the columns
of $\mv{U}_i$ are orthonormal, the following expression is derived for
the Helmholtz operator applied on a tensor in Tucker format
\begin{equation*}
	\mathcal{L}\mt{M} \times_2 \mv{U}_2^\H \times_3 \mv{U}_3^\H
	= \mt{G} \times_1 (-\mv{D}_{xx} - k^2 \mv{I}) \mv{U}_1
	- \mt{G} \times_1 \mv{U}_1 \times_2 \mv{U}_2^\H \mv{D}_{yy} \mv{U}_2
	- \mt{G} \times_1 \mv{U}_1 \times_3 \mv{U}_3^\H \mv{D}_{zz} \mv{U}_3.
\end{equation*}
Writing this tensor equation in the first unfolding leads to a matrix
equation, recall the first unfolding is given by $\mv{M}_{(1)} =
\conj{\mv{U}_1} \mv{G}_{(1)} \left( \mv{U}_3 \otimes \mv{U}_2
\right)^\H$, see also \cite{kolda2009tensor}:
\begin{equation}\label{eq:firstUnfolding}
	\conj{\left(-\mv{D}_{xx} - k^2 \mv{I} \right) \mv{U}_1} \mv{G}_{(1)}
	- \conj{\mv{U}_1} \mv{G}_{(1)} \left(I \otimes \mv{U}_2^\H \mv{D}_{yy} \mv{U}_2\right)^\H
	- \conj{\mv{U}_1} \mv{G}_{(1)} \left(\mv{U}_3^\H \mv{D}_{zz} \mv{U}_3  \otimes I\right)^\H
	=  \mv{F}_{(1)} \left(\mv{U}_3 \otimes \mv{U}_2\right).
\end{equation}
To solve this equation for $\mv{U}_1$, it is written in vectorized form as
\begin{equation}\label{eq:eqnU1}
	\left\{
	\mv{I} \otimes \conj{\left(-\mv{D}_{xx} - k^2 \mv{I} \right)} +
	\left[ 
	-\conj{\left(\mv{I} \otimes \mv{U}_2^\H \mv{D}_{yy} \mv{U}_2\right)} - \conj{\left(\mv{U}_3^\H \mv{D}_{zz} \mv{U}_3 \otimes \mv{I} \right)}
	\right] \otimes \mv{I}
	\right\} \mvec{\underbrace{\conj{\mv{U}_1}\mv{G}_{(1)}}_{\mv{X}_1}} = \mvec{\mv{F}_{(1)} \left(\mv{U}_3 \otimes \mv{U}_2\right)}.
\end{equation}
Observe that this is a square system with $n_1 \times r_2 r_3$ unknowns,
where the solution in matrix form $\mv{X}_1$ could have, in general, a
rank $r > r_1$.
In a similar way, equations for $\mv{U}_2$ and $\mv{U}_3$ are
derived by multiplying \eqref{eq:start3d} with the other factor
matrices in the appropriate directions:
\begin{equation}\label{eq:eqnU2}
	\left\{
	\mv{I} \otimes \conj{\left(-\mv{D}_{yy} - k^2 \mv{I} \right)} +
	\left[ 
	\conj{-\left(\mv{I} \otimes \mv{U}_1^\H \mv{D}_{xx} \mv{U}_1\right)} - \conj{\left(\mv{U}_3^\H \mv{D}_{zz} \mv{U}_3 \otimes \mv{I} \right)}
	\right]  \otimes \mv{I}
	\right\} \mvec{\underbrace{\conj{\mv{U}_2}\mv{G}_{(2)}}_{\mv{X}_2}} 
	= \mvec{\mv{F}_{(2)} \left(\mv{U}_3 \otimes \mv{U}_1\right)},
\end{equation}
and
\begin{equation}\label{eq:eqnU3}
	\left\{
	\mv{I} \otimes \conj{\left(-\mv{D}_{zz} - k^2 \mv{I} \right)} +
	\left[- 
	\conj{\left(\mv{I} \otimes \mv{U}_1^\H \mv{D}_{xx} \mv{U}_1\right)} - \conj{\left(\mv{U}_2^\H \mv{D}_{yy} \mv{U}_2 \otimes \mv{I} \right)}
	\right] \otimes \mv{I}
	\right\} \mvec{\underbrace{\conj{\mv{U}_3}\mv{G}_{(3)}}_{\mv{X}_3}} 
	= \mvec{\mv{F}_{(3)} \left(\mv{U}_2 \otimes \mv{U}_1\right)}.
\end{equation}
Alternating between solving for $\mv{U}_1$, $\mv{U}_2$ and $\mv{U}_3$
using \eqref{eq:eqnU1}, \eqref{eq:eqnU2} or \eqref{eq:eqnU3} results
in algorithm that approximates the low-rank solutions for three
dimensional problems as given in \eqref{eq:driven3d}. This algorithm
is summarized in Algorithm~\ref{alg:const3dv1}.
Also in the three dimensional case the orthogonality of the columns 
of $\mv{U}_1$, $\mv{U}_2$ and $\mv{U}_3$ are maintained by 
additional QR factorizations.
\begin{algorithm}[t]
  \SetAlgoLined
	[$\mt{G}, \mv{U}_1, \mv{U}_2, \mv{U}_3$] = hosvd(initial guess)\;
	\While{not converged}{
		\For{i = 1, 2, 3}{
			Solve for $\mv{X}_i = \conj{\mv{U}_i}\mv{G}_{(i)} \in \C^{n \times r^{d-1}}$ using \eqref{eq:eqnU1}, \eqref{eq:eqnU2} or \eqref{eq:eqnU3}\;
			$\conj{\mv{U}_{i}} \mv{G}_{(i)} = \qr{\mv{X}_i(:,~ 1:r_i), 0}$;
		}
	}
	$\mt{G} = \texttt{reconstruct} \left[\mv{G}_{(i)}, ~i\right]$\;
	$\mt{M} = \mt{G} \times_1 \mv{U}_1  \times_2 \mv{U}_2  \times_3 \mv{U}_3$\;
    \caption{Solve for the low-rank tensor decomposition of the solution $\mt{M}$ of a 3D Helmholtz with constant wave number (version 1).}
    \label{alg:const3dv1}
\end{algorithm}
Observe that we solve for a large matrix $\mv{X}_i \in \C^{n_i \times
  r_1r_2r_3/r_i}$. So, in general the rank of this matrix could be
$\min\left(n_i,~ r_1r_2r_3/r_i\right)$. But it is also known that
$\mv{X}_i = \conj{\mv{U}_i}\mv{G}_{(i)}$ which leads to the fact that
the rank of $\mv{X}_i$ should be at most $r_i$. Selecting the
first $r_i$ columns of $\mv{X}_i$ and computing its QR decomposition
is sufficient to derive a new orthonormal basis as factor matrix
$\conj{\mv{U}}_i$.

\begin{figure}[t]
	\centering
	\begin{subfigure}{0.4\textwidth}
		\includegraphics[width=\textwidth]{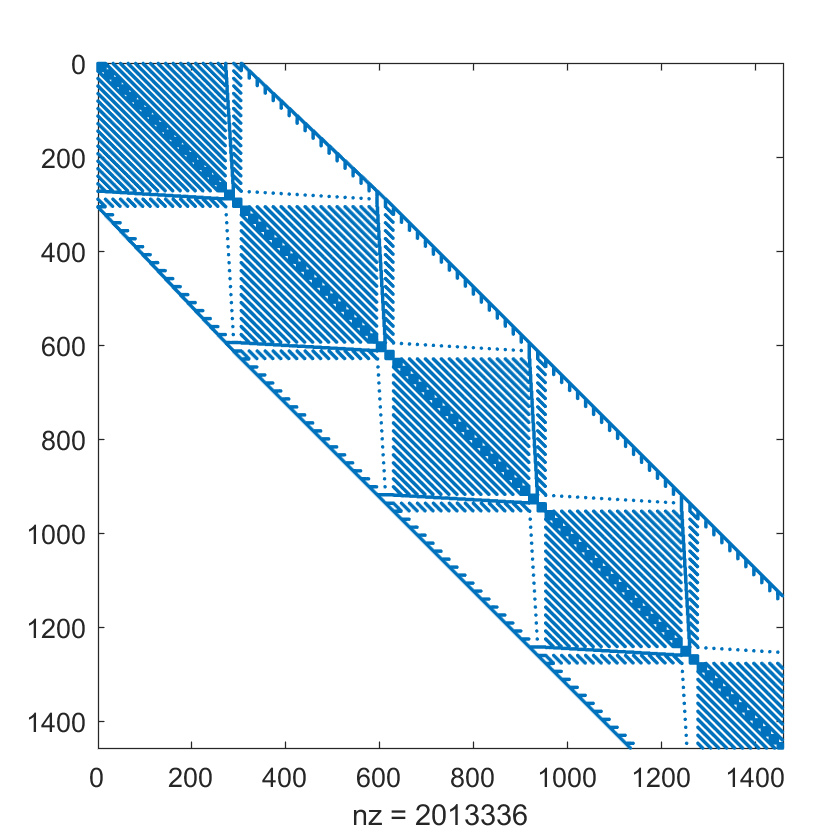}
		\caption{Sparsity pattern of the top of the symmetric reverse Cuthill-McKee permutation of the system matrix to solve for $\mv{X}_i$ using \eqref{eq:eqnU1}. Note: only the first 4.5 of the 168 blocks are shown.}
		\label{fig:spytop_symrcm_UGsystem_alg1}
  	\end{subfigure}
  \quad
	\begin{subfigure}{0.4\textwidth}
		\includegraphics[width=\textwidth]{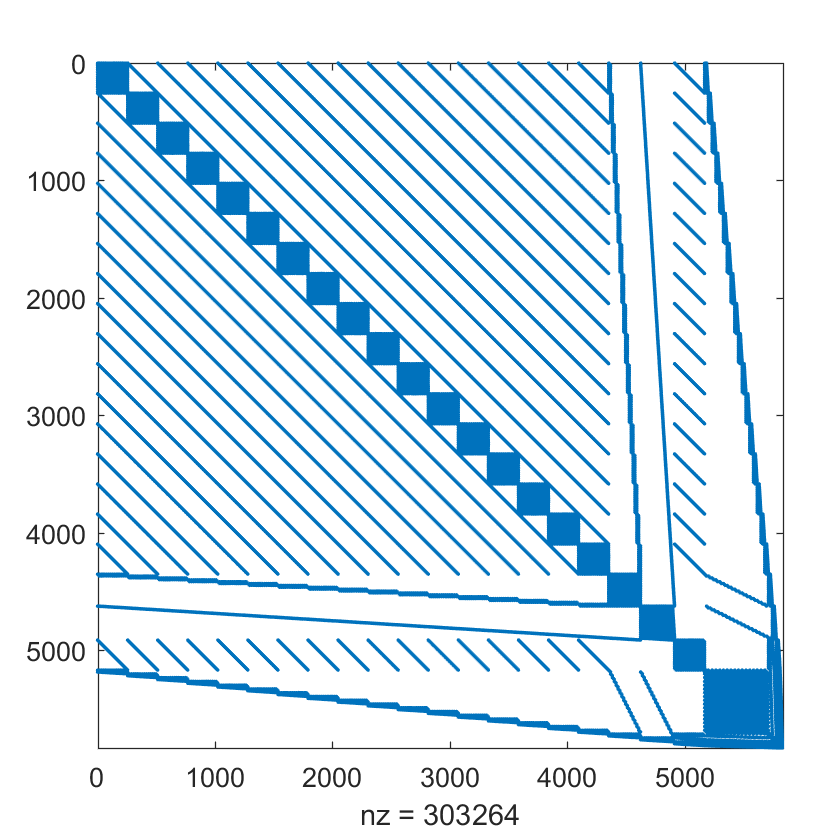}
		\caption{Sparsity pattern of the symmetric reverse Cuthill-McKee permutation of the system matrix to solve for $\mv{G}_{(1)}$ using \eqref{eq:eqnGv2}.}
		\label{fig:spytop_symrcm_Gsystem_alg2}
	\end{subfigure}
	\caption{Sparsity patterns of the symmetric reverse Cuthill-McKee permutation of certain system matrices ($d = 3, n = 168, r = 18$).}
\end{figure}

Finally, observe that solving for $\mv{X}_i$ using \eqref{eq:eqnU1},
\eqref{eq:eqnU2} or \eqref{eq:eqnU3} is computationally not
efficient. In all iterations, we solve for a total of $d nr^{d-1}$
unknowns, while there are only $r^d + dnr$ unknowns in the Tucker
tensor factorization. Furthermore, solving equations \eqref{eq:eqnU1},
\eqref{eq:eqnU2} and \eqref{eq:eqnU3} is also expensive. Indeed,
computing a symmetric reverse Cuthill-McKee permutation of the system
matrix one observes a matrix with a bandwidth $\Oh{r^{d-1}}$. For
example when $d = 3, n_i = n = 168, r_i = r = 18$ one obtains the
sparsity pattern on the diagonal of the matrix as shown in Figure
\ref{fig:spytop_symrcm_UGsystem_alg1}. So solving a system as given in
\eqref{eq:eqnU1}, \eqref{eq:eqnU2} or \eqref{eq:eqnU3} has a
computational cost of $\Oh{nr^{2(d-1)}}$.

\subsubsection{Solving for the basis functions and the core tensor separately (version 2)}\label{sec:version2}
To circumvent solving the large systems in \eqref{eq:eqnU1},
\eqref{eq:eqnU2} and \eqref{eq:eqnU3}, we can pre-compute the
QR factorization of the unfolding of the core tensor, $\mv{G}_{(i)}$,
and project the equations onto the obtained $\mv{Q}_i$.  Indeed, this
will further reduce the number of unknowns in these linear systems to exactly
the number of unknowns that are needed for the factor matrices
$\conj{\mv{U}}_i$, for $i = 1,2,3$.

Let us discuss the details. We start again from equation \eqref{eq:firstUnfolding} and use the QR
factorization of $\mv{G}_{(1)}^\H$, $\mv{Q}_1 \mv{R}_1^\H = \qr{\mv{G}^\H_{(1)}}$. This yields
\begin{equation*}
	\conj{\left(-\mv{D}_{xx} - k^2 \mv{I} \right) \mv{U}_1} \mv{R}_1\mv{Q}_1^\H
	- \conj{\mv{U}_1} \mv{R}_1\mv{Q}_1^\H \left(I \otimes \mv{U}_2^\H \mv{D}_{yy} \mv{U}_2\right)^\H
	- \conj{\mv{U}_1} \mv{R}_1\mv{Q}_1^\H \left(\mv{U}_3^\H \mv{D}_{zz} \mv{U}_3  \otimes I\right)^\H
	=  \mv{F}_{(1)} \left(\mv{U}_3 \otimes \mv{U}_2\right).
\end{equation*}
Post multiplication of this equation by $\mv{Q}_1$ yields
\begin{equation*}
	\conj{\left(-\mv{D}_{xx} - k^2 \mv{I} \right) \mv{U}_1} \mv{R}_1
	- \conj{\mv{U}_1} \mv{R}_1\mv{Q}_1^\H \left(I \otimes \mv{U}_2^\H \mv{D}_{yy} \mv{U}_2\right)^\H \mv{Q}_1
	- \conj{\mv{U}_1} \mv{R}_1\mv{Q}_1^\H \left(\mv{U}_3^\H \mv{D}_{zz} \mv{U}_3  \otimes I\right)^\H \mv{Q}_1
	=  \mv{F}_{(1)} \left(\mv{U}_3 \otimes \mv{U}_2\right)\mv{Q}_1.
\end{equation*}
To solve this equation for $\mv{U}_1$, it is written is vectorized form
as
\begin{equation}\label{eq:eqnU1v2}
\left\{
\mv{I} \otimes \conj{\left(-\mv{D}_{xx} - k^2 \mv{I} \right)} +
\mv{Q}_1^\T\left[ -\conj{\left(\mv{I} \otimes \mv{U}_2^\H \mv{D}_{yy} \mv{U}_2\right)} - \conj{\left(\mv{U}_3^\H \mv{D}_{zz} \mv{U}_3 \otimes \mv{I} \right)}
\right]\conj{\mv{Q}_1} \otimes \mv{I}
\right\} \mvec{\underbrace{\conj{\mv{U}_1}\mv{R}_1}_{\mv{X}_1}} = \mvec{\mv{F}_{(1)} \left(\mv{U}_3 \otimes \mv{U}_2\right) \mv{Q}_1}.
\end{equation}
In a similar way, the update equations for $\mv{U}_2$ and $\mv{U}_3$ are
derived by multiplying \eqref{eq:start3d} with the other factor
matrices in the appropriate dimensions and using the QR factorizations
of $\mv{G}_{(i)}^\H$:
\begin{equation}\label{eq:eqnU2v2}
\left\{
\mv{I} \otimes \conj{\left(-\mv{D}_{yy} - k^2 \mv{I} \right)} +
\mv{Q}_2^\T\left[- 
\conj{\left(\mv{I} \otimes \mv{U}_1^\H \mv{D}_{xx} \mv{U}_1\right)} - \conj{\left(\mv{U}_3^\H \mv{D}_{zz} \mv{U}_3 \otimes \mv{I} \right)}
\right]\conj{\mv{Q}_2}  \otimes \mv{I}
\right\} \mvec{\underbrace{{\conj{\mv{U}_2}\mv{R}_2}}_{\mv{X}_2}} 
= \mvec{\mv{F}_{(2)} \left(\mv{U}_3 \otimes \mv{U}_1\right) \mv{Q}_2}.
\end{equation}
and
\begin{equation}\label{eq:eqnU3v2}
\left\{
\mv{I} \otimes \conj{\left(-\mv{D}_{zz} - k^2 \mv{I} \right)} +
\mv{Q}_3^\T\left[ -
\conj{\left(\mv{I} \otimes \mv{U}_1^\H \mv{D}_{xx} \mv{U}_1\right)} - \conj{\left(\mv{U}_2^\H \mv{D}_{yy} \mv{U}_2 \otimes \mv{I} \right)}
\right]\conj{\mv{Q}_3} \otimes \mv{I}
\right\} \mvec{\underbrace{\conj{\mv{U}_3}\mv{R}_3}_{\mv{X}_3}} 
= \mvec{\mv{F}_{(3)} \left(\mv{U}_2 \otimes \mv{U}_1\right) \mv{Q}_3}.
\end{equation}
All these equations are cheap to solve. Indeed,
$\mvec{\conj{\mv{U}_i}\mv{R}_i}$ has length $n_i r_i$. Computing a
symmetric reverse Cuthill-McKee permutation of these system matrices one
observes a matrix with a bandwidth $\Oh{r}$, so solving these
equations has a computational cost $\Oh{nr^2}$.

Of course, this only updates the factor matrices as basis vectors in each
direction. As a single final step, we still have to compute the core
tensor $\mt{G}$. This will be the computationally most expensive
part.

Core tensor $\mt{G}$ can be obtained by multiplying \eqref{eq:start3d}
with all the $d$ factor matrices in the matching directions.
Unfolding this equation in a certain direction (eg. the first folding)
leads again to a matrix equation. In vectorized form, it is given by
\begin{equation}\label{eq:eqnGv2}
	\left\{
	\mv{I} \otimes \conj{\mv{U}_1^\H\left(-\mv{D}_{xx} - k^2 \mv{I} \right)\mv{U}_1} +
	\left[ 
	-\conj{\left(\mv{I} \otimes \mv{U}_2^\H \mv{D}_{yy} \mv{U}_2\right)} - \conj{\left(\mv{U}_3^\H \mv{D}_{zz} \mv{U}_3 \otimes \mv{I} \right)}
	\right] \otimes \mv{I}
	\right\} \mvec{\mv{G}_{(1)}} = \mvec{\mv{U}_1^\T \mv{F}_{(1)} \left(\mv{U}_3 \otimes \mv{U}_2\right)}.
\end{equation}
Indeed, considering again an example where $d = 3, n_i = n = 168, r_i
= r = 18$ one obtains a matrix with a sparsity pattern that is shown
in  Figure~\ref{fig:spytop_symrcm_Gsystem_alg2}. Hence, this matrix has
not a limited bandwidth anymore. It coupled all functions to all other
functions. Although this equation has to be solved only once in the 
algorithm, when the rank increases, it will rapidly dominate the computational 
cost of this algorithm.

\begin{algorithm}[t]
	[$\mt{G}, \mv{U}_1, \mv{U}_2, \mv{U}_3$] = hosvd(initial guess)\;
	\While{not converged}{
		\For{i = 1, 2, 3}{
			$\mv{Q}_i \widetilde{\mv{R}} = \qr{\mv{G}_{(i)}^\H, 0}$\;
			Solve for $\mv{X}_i = \conj{\mv{U}_i}\mv{R}_i \in \C^{n_i \times r_i}$ using \eqref{eq:eqnU1v2}, \eqref{eq:eqnU2v2} or \eqref{eq:eqnU3v2}\;
			$\conj{\mv{U}_{i}} \mv{R}_i = \qr{\mv{X}_i, 0}$\;
			$\mt{G} = \texttt{reconstruct} \left[\mv{R}_i\mv{Q}_i^\H, ~i\right]$;
		}
	}
	Solve for $\mv{G}_{(1)} \in \C^{r_1 \times r_2r_3}$ using \eqref{eq:eqnGv2}\;
	$\mt{G} = \texttt{reconstruct} \left[\mv{G}_{(1)}, ~1\right]$\;
	$\mt{M} = \mt{G} \times_1 \mv{U}_1  \times_2 \mv{U}_2  \times_3 \mv{U}_3$\;
    \caption{Solve for the low-rank tensor decomposition of the solution $\mt{M}$ of a 3D Helmholtz problem with constant wave number (version 2).}
    \label{alg:const3dv2}
\end{algorithm}

\subsubsection{Efficient combination of version 1 and version 2 into new algorithm (version 3)}\label{sec:version3}
In the first version of the algorithm, (see section~\ref{sec:version1}),
an update for $\mv{G}_{(i)}$ is computed for each direction in each
iteration. This leads to a too expensive algorithm.
Then we changed the algorithm such that the costs for the updates in each
direction is reduced, (see section~\ref{sec:version2}). But, in that version
almost all information for a full update of core tensor $\mt{G}$ is
lost. Therefore a final, but potential too expensive, equation needs to be solved.

Observe that the expensive computation for the full core tensor, in 
version 2, can now be replaced by a single solve per iteration as done in
version 1. This leads to a third version of the algorithm.
It avoids repeatedly solving the large systems (like version 1) and it
does not solve too expensive systems (like version 2). The computational
complexity of this algorithm is equal to the complexity of version 1,
so $\Oh{nr^{2(d-1)}}$.
Furthermore, the systems that need to be solved, each iteration, have 
exactly the same number of unknowns as the representation of the 
tensor in low-rank Tucker tensor format.
In summary, this final version of the algorithm is given by Algorithm~\ref{alg:const3dv3}.

\begin{algorithm}[t]
	[$\mt{G}, \mv{U}_1, \mv{U}_2, \mv{U}_3$] = hosvd(initial guess)\;
	\While{not converged}{
		\For{i = 1, 2}{
			$\mv{Q}_i \widetilde{\mv{R}} = \qr{\mv{G}_{(i)}^\H, 0}$\;
			Solve for $\mv{X}_i = \conj{\mv{U}_i}\mv{R}_i \in \C^{n_i \times r_i}$ using \eqref{eq:eqnU1v2} or \eqref{eq:eqnU2v2}\;
			$\conj{\mv{U}_{i}} \mv{R}_i = \qr{\mv{X}_i, 0}$\;
			$\mt{G} = \texttt{reconstruct} \left[\mv{R}_i\mv{Q}_i^\H, ~i\right]$;
		}
		Solve for $\mv{X}_3 = \conj{\mv{U}_3}\mv{G}_{(3)} \in \C^{n_3 \times r^{d-1}}$ using \eqref{eq:eqnU3}\;
		$\conj{\mv{U}_{3}} \mv{G}_{(3)} = \qr{\mv{X}_3, 0}$\;
		$\mt{G} = \texttt{reconstruct} \left[\mv{G}_{(3)},~3\right]$\;
	}
	$\mt{M} = \mt{G} \times_1 \mv{U}_1  \times_2 \mv{U}_2  \times_3 \mv{U}_3$\;
    \caption{Solve for the low-rank tensor decomposition of the solution $\mt{M}$ of a 3D Helmholtz problem with constant wave number (version 3).}
    \label{alg:const3dv3}
\end{algorithm}

\subsubsection{Numerical comparison of three versions for 3D Helmholtz equation}
Consider a three dimensional domain $\Omega = [-10, ~10]^3$ that is
discretized with $M=100$ equidistant mesh points per direction in the
interior of the domain. The domain is extended with exterior complex
scaling to implement the absorbing boundary conditions. Hence, in total
there are $n = n_x = n_y = n_z = 168$ unknowns per direction. As
constant wave number we use $\omega = 2$ and a right hand side $f(x,y,z)
= -e^{-x^2-y^2-z^2}$.
By symmetry, we expect a low rank factorization with a equal low-rank 
in each direction, so we fix $r = r_x = r_y = r_x$.

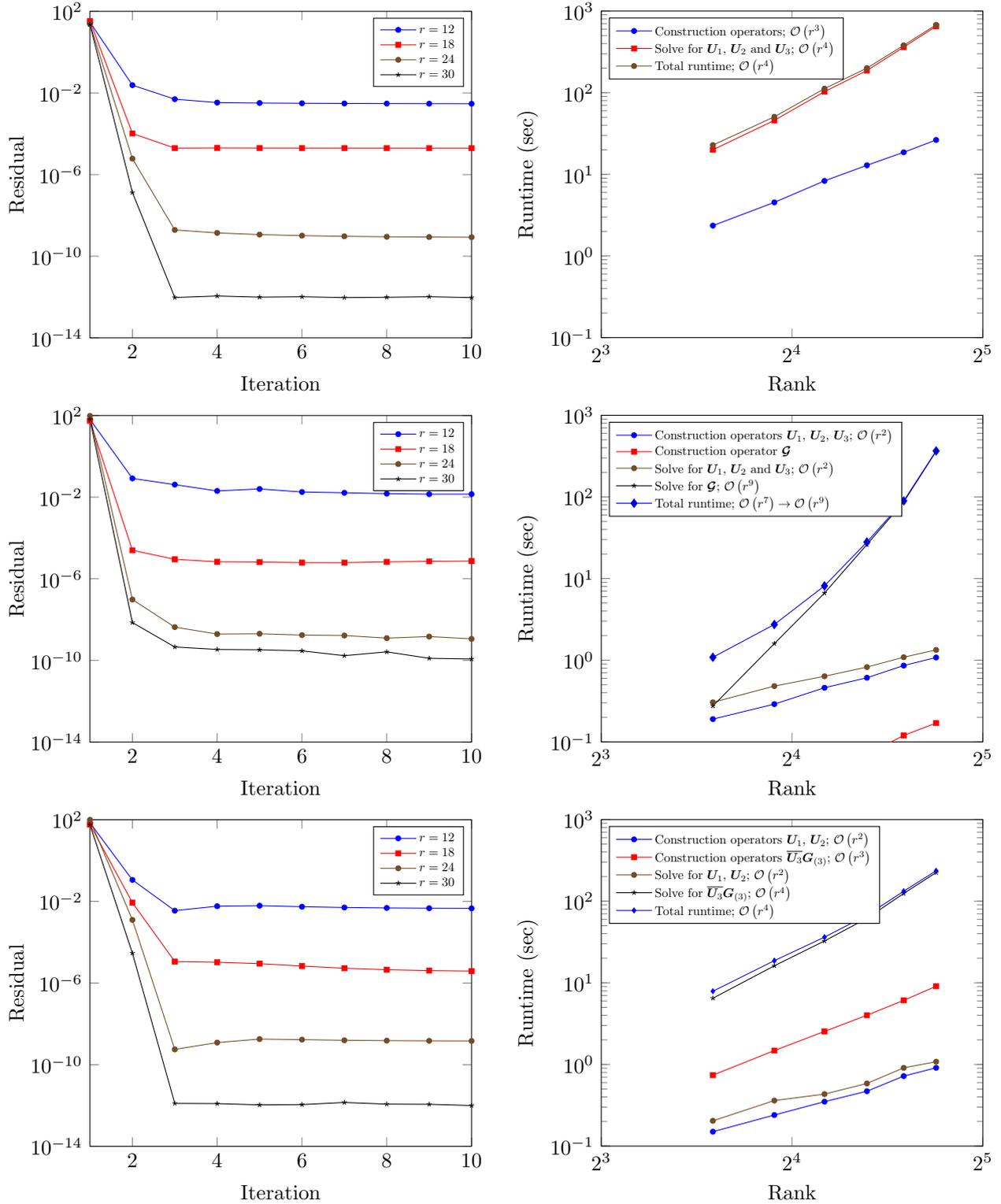
\begin{figure}
	\centering
	\input{residual-v1-M=100.tikz} ~ \input{runtime-v1-3d-M=100.tikz} \par
	\input{residual-v2-M=100.tikz} ~ \input{runtime-v2-3d-M=100.tikz} \par
	\input{residual-v3-M=100.tikz} ~ \input{runtime-v3-3d-M=100.tikz}
        \caption{Left: Plot of residual per iteration for constant wave number in 3D Helmholtz problem.
        	Right: Plot of runtime of most time consuming parts for constant wave number in 3D Helmholtz problem. Both problems have $M=100$.
        	Top: Algorithm~\ref{alg:const3dv1} (version 1),
        	middle: Algorithm~\ref{alg:const3dv2} (version 2),
        	bottom: Algorithm~\ref{alg:const3dv3} (version 3).}
	\label{fig:residual+runtime-const3d}
\end{figure}

The convergence of the residuals of the three versions are given in
the left column of Figure~\ref{fig:residual+runtime-const3d}.
It is clear that all three versions converge to a good low-rank 
approximation of the full solution. By increasing the maximal
attainable rank $r$, a better low-rank solution is obtained, as 
expected. Remarkably,  for $r=30$, in version 2, the final residual
is larger then the residuals obtained by both other algorithms while
the compute-time for version 2 is larger than the other algorithms.

The compute-time for the most time-consuming parts in the different
versions of the algorithm can be measured as a function of the maximal
attainable rank $r$. For the three versions of the algorithm the
runtimes are shown in the right column of
Figure~\ref{fig:residual+runtime-const3d}. For all parts the expected
and measured dependence on the rank $r$ are given. For all versions of
the algorithm 10 iterations are applied.

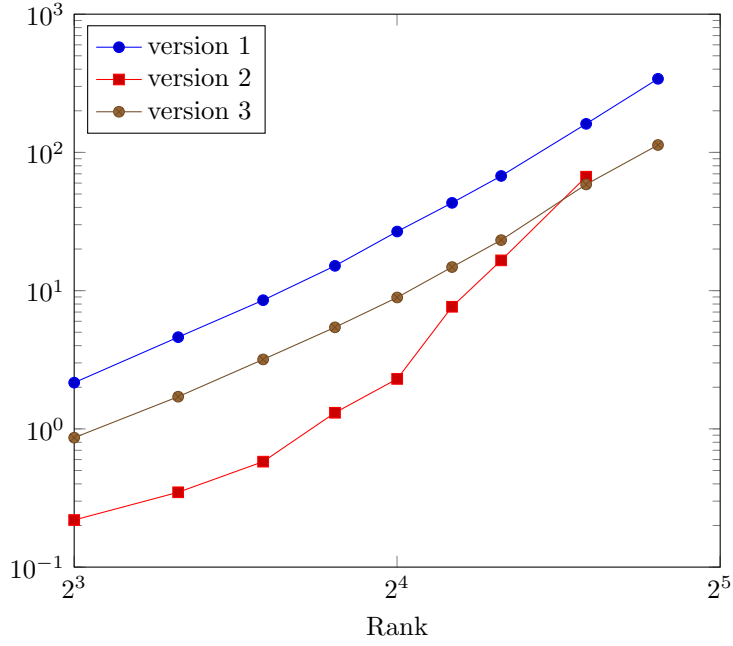
\begin{figure}
	\centering
 	\input{runtimes-v123-M=100-iter=4.tikz}
    \caption{Plot of runtime for 4 iteration with constant wavenumber in 3D using the three different algorithms ($M=100$)
      	\label{fig:runtime-const3d}}
\end{figure}

Comparing the total runtime for the three different versions one obtains
results as shown in Figure~\ref{fig:runtime-const3d}. Indeed, as
expected version 3 is approximately 3 times faster than version 1 and
the runtime scales similar in rank $r$. Further, for small rank $r$
version 2 is faster than both other versions. But when the rank
increases the expensive solve for the core tensor $\mt{G}$ starts to
dominate the runtime. The total runtime will increase dramatically.

\subsection{Projection operator for constant wave number}
Also in three dimensions we can write the linear systems
\eqref{eq:eqnU1} for $\mv{U}_1$, \eqref{eq:eqnU2} for 
$\mv{U}_2$ and \eqref{eq:eqnU3} for $\mv{U}_3$ as
projection operators applied to the residual of the tensor equation,
\eqref{eq:start3d}.

Consider a tensor $\mt{M}$ in Tucker format and factorized as $\mt{M}
= \mt{G} \times_1 \mv{U}_1 \times_2 \mv{U}_2 \times_3 \mv{U}_3$, with
unknowns $\mt{G}, \mv{U}_1, \mv{U}_2$ and $\mv{U}_3$. Discretization
of \eqref{eq:driven3d} leads to a linear operator $\mathcal{L}$
applied on tensors. Its matrix representation $\mv{L}$ has a sum
of Kronecker products structure, as given in
\eqref{eq:operatorLasSumKroneckerProducts}.

Solving for an unknown factors $\mv{U}_1$, $\mv{U}_2$ or $\mv{U}_3$
(and the core-tensor $\mt{G}$) using \eqref{eq:eqnU1},
\eqref{eq:eqnU2} or \eqref{eq:eqnU3} can be interpreted as a
projection operator applied on the residual.  For example,
\eqref{eq:eqnU1} can be interpreted as
\begin{equation}
	\conj{\left(\mv{U}_3^\H \otimes \mv{U}_2^\H \otimes \mv{I} \right) \mv{L} \left(\mv{U}_3 \otimes \mv{U}_2 \otimes \mv{I}\right)} \mvec{\conj{\mv{U}_1}\mv{G}_{(1)}} = \left( \mv{U}_3^\T \otimes \mv{U}_2^\T \otimes \mv{I} \right) \mvec{\mv{F}_{(1)}}.
\end{equation}
The residual, in tensor format, is given by
\begin{equation}
	\begin{split}
		\mt{R} &= \mt{F} - \mathcal{L}\mt{M}, \\
		&= \mt{F} - \mt{G} \times_1 \left(-\mv{D}_{xx} - k^2 \mv{I}\right) \mv{U}_1 \times_2 \mv{U}_2 \times_3 \mv{U}_3 + \mt{G} \times_1 \mv{U}_1 \times_2 \mv{D}_{yy}\mv{U}_2 \times_3 \mv{U}_3 + \mt{G} \times_1 \mv{U}_1 \times_2 \mv{U}_2 \times_3 \mv{D}_{zz}\mv{U}_3.
	\end{split}
\end{equation}
Writing this tensor equation in the first unfolding leads to the following matrix equation
\begin{equation}
		\mv{R}_{(1)} = \mv{F}_{(1)} - \conj{\left(-\mv{D}_{xx} - k^2 \mv{I}\right)\mv{U}_1}\mv{G}_{(1)} \left(\mv{U}_3 \otimes \mv{U}_2 \right)^\H + \conj{\mv{U}_1}\mv{G}_{(1)} \left(\mv{U}_3 \otimes \mv{D}_{yy}\mv{U}_2 \right)^\H + \conj{\mv{U}_1}\mv{G}_{(1)} \left(\mv{D}_{zz}\mv{U}_3 \otimes \mv{U}_2 \right)^\H,
\end{equation}
which can be vectorized as
\begin{equation}
	\begin{split}
		\mvec{\mv{R}_{(1)}} &= \mvec{\mv{F}_{(1)}} - \left(\conj{\left(\mv{U}_3 \otimes \mv{U}_2 \otimes \left(-\mv{D}_{xx}-k^2\mv{I}\right)\right)} - \conj{\left(\mv{U}_3 \otimes \mv{D}_{yy}\mv{U}_2 \otimes \mv{I}\right)} - \conj{\left(\mv{D}_{zz}\mv{U}_3 \otimes \mv{U}_2 \otimes \mv{I}\right)}\right)\mvec{\conj{\mv{U}_1}\mv{G}_{(1)}} \\
		&= \mvec{\mv{F}_{(1)}} - \conj{\mv{L} \left(\mv{U}_3 \otimes \mv{U}_2 \otimes \mv{I}\right)}\mvec{\conj{\mv{U}_1}\mv{G}_{(1)}} \\
		&= \mvec{\mv{F}_{(1)}} - \conj{\mv{L} \left(\mv{U}_3 \otimes \mv{U}_2 \otimes \mv{I}\right)} \left[\conj{\left(\mv{U}_3^\H \otimes \mv{U}_2^\H \otimes \mv{I} \right) \mv{L} \left(\mv{U}_3 \otimes \mv{U}_2 \otimes \mv{I}\right)}\right]^{-1}\left( \mv{U}_3^\T \otimes \mv{U}_2^\T \otimes \mv{I} \right) \mvec{\mv{F}_{(1)}} \\
		&= P_{23}\mvec{\mv{F}_{(1)}},
	\end{split}
\end{equation}
where operator $P_{23}$ is given by
\begin{equation}\label{eq:projector3dconst}
	\begin{split}
		P_{23} &= \mv{I} - \conj{\mv{L} \left(\mv{U}_3 \otimes \mv{U}_2 \otimes \mv{I}\right)} \left[\conj{\left(\mv{U}_3^\H \otimes \mv{U}_2^\H \otimes \mv{I} \right) \mv{L} \left(\mv{U}_3 \otimes \mv{U}_2 \otimes \mv{I}\right)}\right]^{-1}\left( \mv{U}_3^\T \otimes \mv{U}_2^\T \otimes \mv{I} \right) \\
		&= \mv{I} - \mv{X}.
	\end{split}
\end{equation}
This operator $P_{23}$ is indeed a projection operator. Observe that the terms between the two inverses cancel against one of the inverse factors:
\begin{equation*}
	\begin{split}
		\mv{X}^2 &= \conj{\mv{L} \left(\mv{U}_3 \otimes \mv{U}_2 \otimes \mv{I}\right)} \left[\conj{\left(\mv{U}_3^\H \otimes \mv{U}_2^\H \otimes \mv{I} \right) \mv{L} \left(\mv{U}_3 \otimes \mv{U}_2 \otimes \mv{I}\right)}\right]^{-1}\left( \mv{U}_3^\T \otimes \mv{U}_2^\T \otimes \mv{I} \right)\conj{\mv{L} \left(\mv{U}_3 \otimes \mv{U}_2 \otimes \mv{I}\right)} \left[\conj{\left(\mv{U}_3^\H \otimes \mv{U}_2^\H \otimes \mv{I} \right) \mv{L} \left(\mv{U}_3 \otimes \mv{U}_2 \otimes \mv{I}\right)}\right]^{-1}\left( \mv{U}_3^\T \otimes \mv{U}_2^\T \otimes \mv{I} \right) \\
		&= \conj{\mv{L} \left(\mv{U}_3 \otimes \mv{U}_2 \otimes \mv{I}\right)} \left[\conj{\left(\mv{U}_3^\H \otimes \mv{U}_2^\H \otimes \mv{I} \right) \mv{L} \left(\mv{U}_3 \otimes \mv{U}_2 \otimes \mv{I}\right)}\right]^{-1}\left( \mv{U}_3^\T \otimes \mv{U}_2^\T \otimes \mv{I} \right) \\
		&= \mv{X}.
	\end{split}
\end{equation*}
This operator is a natural extension to higher dimensions of the two
dimensional operators as derived in section \ref{sec:projection2d}.

A similar derivation results in projection operators $P_{13}$ and $P_{12}$ for the updates in $\mv{U}_2$ and $\mv{U}_3$, respectively.
\begin{equation}\label{eq:projectors3dconst}
	\begin{split}
		P_{23} &= \mv{I} - \conj{\mv{L} \left(\mv{U}_3 \otimes \mv{U}_2 \otimes \mv{I}\right)} \left[\conj{\left(\mv{U}_3^\H \otimes \mv{U}_2^\H \otimes \mv{I} \right) \mv{L} \left(\mv{U}_3 \otimes \mv{U}_2 \otimes \mv{I}\right)}\right]^{-1}\left( \mv{U}_3^\T \otimes \mv{U}_2^\T \otimes \mv{I} \right), \\
		P_{13} &= \mv{I} - \conj{\mv{L} \left(\mv{U}_3 \otimes \mv{I} \otimes \mv{U}_1\right)} \left[\conj{\left(\mv{U}_3^\H \otimes \mv{I} \otimes \mv{U}_1^\H \right) \mv{L} \left(\mv{U}_3 \otimes \mv{I} \otimes \mv{U}_1\right)}\right]^{-1}\left( \mv{U}_3^\T \otimes \mv{I} \otimes \mv{U}_1^\T \right), \\
		P_{12} &= \mv{I} - \conj{\mv{L} \left(\mv{I} \otimes \mv{U}_2 \otimes \mv{U}_1\right)} \left[\conj{\left(\mv{I} \otimes \mv{U}_2^\H \otimes \mv{U}_1^\H \right) \mv{L} \left(\mv{I} \otimes \mv{U}_2 \otimes \mv{U}_1\right)}\right]^{-1}\left( \mv{I} \otimes \mv{U}_2^\T \otimes \mv{U}_1^\T \right).
	\end{split}
\end{equation}
The successive application of these projection operators on the
residual results in an updated residual that lies in the intersection
of all subspaces.

\subsection{Helmholtz equation with space-dependent wave number}
The presented algorithms  with constant wave number can be
extended to space-dependent wave numbers. So, lets
consider a 3D Helmholtz problem  where $\mt{K} = k^2(x,y,z)$ represents the space-dependent
wave number on the discretized mesh.

Further, we assume that a Canonical Polyadic decomposition of the
space-dependent wave number tensor 
$\mt{K}$ is known, i.e.
\begin{equation}\label{eq:CP4K}
	\mt{K} = \sum_{i=1}^s \sigma_i \left(\mv{v}_i^{(1)} \circ \mv{v}_i^{(2)} \circ \cdots \circ	\mv{v}_i^{(d)}\right),
\end{equation}
where $s \in \N_+$ is the CP-rank of $\mt{K}$ and $\mv{v}_i^{(j)} \in \C^{n_j}$ for $i = 1,2,\ldots, s; j = 1,2,\ldots, d$ are vectors. Further, $\sigma_i$ is a tensor generalization of a singular value and $\circ$ denotes the vector outer product.

The application of the space-dependent Helmholtz operator $\mathcal{L}$ on tensor $\mt{M}$ is given by
\begin{equation}\label{eq:start3dvar}
	\begin{split}
		\mathcal{L} \mt{M} &= \mt{F} \\
		\mathcal{L} \mt{M} &= -\mt{G} \times_1 \mv{D}_{xx} \mv{U}_1 \times_2 \mv{U}_2 \times_3 \mv{U}_3\\
		&- \mt{G} \times_1 \mv{U}_1 \times_2 \mv{D}_{yy} \mv{U}_2 \times_3 \mv{U}_3\\
		&- \mt{G} \times_1 \mv{U}_1 \times_2 \mv{U}_2 \times_3 \mv{D}_{zz} \mv{U}_3\\
		&- \mt{K} \circ \left(\mt{G} \times_1 \mv{U}_1 \times_2 \mv{U}_2 \times_3 \mv{U}_3 \right)\\
		&= \mt{F},
	\end{split}
\end{equation}
where $\mv{U}_i^\H \mv{U}_i = \mv{I}$ for $i = 1,2,3$ and $\mt{F}$
is a tensor representation of the right hand side function $f$ 
discretized on the used grid.
Here $\circ$ denotes the Hadamard product for tensors.

In a similar way as in the three dimensional constant wave number
case, we can derive equations to iteratively solve for the factors
$\mv{U}_1$, $\mv{U}_2$ and $\mv{U}_3$. We start from \eqref{eq:start3dvar} and multiply with
$\mv{U}_2$ and $\mv{U}_3$ in the second and third direction,
respectively. Using that the columns of $\mv{U}_i$ are orthonormal,
the following expression is derived:
\begin{equation*}
	\mathcal{L} \mt{M}  \times_2 \mv{U}_2^\H \times_3 \mv{U}_3^\H = -\mt{G} \times_1 \mv{D}_{xx}
	- \mt{G} \times_1 \mv{U}_1 \times_2 \mv{U}_2^\H\mv{D}_{yy} \mv{U}_2
	- \mt{G} \times_1 \mv{U}_1 \times_3 \mv{U}_3^\H\mv{D}_{zz} \mv{U}_3
	- \left[\mt{K} \circ \left(\mt{G} \times_1 \mv{U}_1 \times_2 \mv{U}_2 \times_3 \mv{U}_3 \right) \right] \times_2 \mv{U}_2^\H \times_3 \mv{U}_3^\H
\end{equation*}
Written in the first unfolding, the multiplication with $\mv{U}_2^\H$ and $\mv{U}_3^\H$ in, respectively, the second and third direction is equivalent to post-multiplication with the matrix
$\left(\mv{I} \otimes \mv{U}_2^\H \right)^\H \left(\mv{U}_3^\H \otimes \mv{I} \right)^\H = \left(\mv{U}_3^\H \otimes \mv{U}_2^\H \right)^\H = \left(\mv{U}_3 \otimes \mv{U}_2 \right)$.

Most of the terms are equal to the case where we had a constant wave
number, see also \eqref{eq:start3d}. Let us focus on the last term 
that contains the Hadamard product with the space-dependent wave number, i.e.:
\begin{equation}\label{eq:hadamardProduct}
	\begin{split}
		\mt{K} &\circ \left(\mt{G} \times_1 \mv{U}_1 \times_2 \mv{U}_2 \times_3 \mv{U}_3 \right).
	\end{split}
\end{equation}
For Hadamard products of tensors, $\mt{Z} = \mt{X} \circ \mt{Y}$, the following property for the $k$-th unfolding holds $\mv{Z}_{(k)} = \mv{X}_{(k)} \circ \mv{Y}_{(k)}$.
Thus, written in the first unfolding \eqref{eq:hadamardProduct} is given by
\begin{equation}\label{eq:hadamardProduct1stUnfolding}
	\begin{split}
		\mv{K}_{(1)} &\circ \mv{M}_{(1)}\\
		\mv{K}_{(1)} &\circ \left( \conj{\mv{U}_1} \mv{G}_{(1)} \left( \mv{U}_3 \otimes \mv{U}_2 \right)^\H \right).
	\end{split}
\end{equation}

As the Hadamard product-term \eqref{eq:hadamardProduct1stUnfolding} is written in the first unfolding and multiplication with $\mv{U}_2^\H$ and $\mv{U}_3^\H$ in respectively the second and third dimension results in
\begin{equation}\label{eq:lastTerm}
	\left[\underbrace{\mv{K}_{(1)}}_{\mv{K}} \circ \underbrace{\conj{\mv{U}_1}\mv{G}_{(1)}}_{\mv{U}} \underbrace{\left( \mv{U}_3 \otimes \mv{U}_2 \right)^\H}_{\mv{V}^\H} \right] \underbrace{\left(\mv{U}_3 \otimes \mv{U}_2 \right)}_{\mv{V}}.
\end{equation}

The derivation of the other terms of \eqref{eq:start3dvar} are equal to the constant wave number case.
The equation in the first unfolding leads to a matrix equation:
\begin{equation}\label{eq:firstUnfoldingVar}
	-\conj{\mv{D}_{xx}\mv{U}_1} \mv{G}_{(1)}
	- \conj{\mv{U}_1} \mv{G}_{(1)} \left(I \otimes \mv{U}_2^\H \mv{D}_{yy} \mv{U}_2\right)^\H
	- \conj{\mv{U}_1} \mv{G}_{(1)} \left(\mv{U}_3^\H \mv{D}_{zz} \mv{U}_3  \otimes I\right)^\H
	- \left[\mv{K}_{(1)} \circ \conj{\mv{U}_1}\mv{G}_{(1)} \left( \mv{U}_3 \otimes \mv{U}_2 \right)^\H \right] \left(\mv{U}_3 \otimes \mv{U}_2 \right)
	=  \mv{F}_{(1)} \left(\mv{U}_3 \otimes \mv{U}_2\right).
\end{equation}

Vectorization of the last term, i.e. \eqref{eq:lastTerm}, results again in an expression for the space-dependent wave number of the form $\left(\mv{K} \circ \mv{U}\mv{V}^\H\right)\mv{V}$, similar to the two dimensional case which was given in \eqref{eq:2dvarwavenr}.
Using again \eqref{eq:vecThirdTerm}, the vectorization of this expression is given by
\begin{equation}\label{eq:tmp1}
	\left(\mv{U}_3^\T \otimes \mv{U}_2^\T \otimes \mv{I} \right) \diag{\mvec{\mv{K}_{(1)}}} \conj{\left(\mv{U}_3 \otimes 	\mv{U}_2 \otimes \mv{I} \right)} \mvec{\conj{\mv{U}_1}\mv{G}_{(1)}}.
\end{equation}

Because $\mt{K}$ is known in a Canonical Polyadic tensor (CP tensor) decomposition\footnote{Otherwise a Canonical Polyadic tensor decomposition can be computed using for example an CP-ALS algorithm \cite{kolda2009tensor}.}, as given in \eqref{eq:CP4K}, we have
\begin{equation}
	\diag{\mvec{\mv{K}_{(1)}}} = \sum_{i=1}^{s} \sigma_i \diag{\mv{v}_{i}^{(3)}} \otimes \diag{\mv{v}_{i}^{(2)}} \otimes \diag{\mv{v}_{i}^{(1)}}.
\end{equation}

So, using the CP tensor representation of the space-dependent wave number the vectorization in \eqref{eq:tmp1} simplifies even further:
\begin{equation*}
	\sum_{i=1}^{s} \sigma_i \left(\mv{U}_3^\T \otimes \mv{U}_2^\T \otimes \mv{I} \right) \left[ \diag{\mv{v}_{i}^{(3)}} \otimes \diag{\mv{v}_{i}^{(2)}} \otimes \diag{\mv{v}_{i}^{(1)}} \right] \conj{\left(\mv{U}_3 \otimes \mv{U}_2 \otimes \mv{I}\right)} \mvec{\conj{\mv{U}_1}\mv{G}_{(1)}}
\end{equation*}
which reduced to
\begin{equation*}
	\underbrace{\sum_{i=1}^{s} \sigma_i
		\left(\mv{U}_3^\T \diag{\mv{v}_{i}^{(3)}} \conj{\mv{U}_3} \right) \otimes
		\left(\mv{U}_2^\T \diag{\mv{v}_{i}^{(2)}} \conj{\mv{U}_2} \right) \otimes
		\left(\diag{\mv{v}_{i_1}^{(1)}}\right)}_{K_1} \mvec{\conj{\mv{U}_1}\mv{G}_{(1)}}.
\end{equation*}
In this way the $K_1$ operator is defined and can be applied to $\mvec{\conj{\mv{U}_1}\mv{G}_{(1)}}$. Observe that this expansion is only advantageous if the space-dependent wave number has low rank, which is typical the case for our applications.

In a similar way, the $K_2$ and $K_3$ operators can be derived:
\begin{equation}\label{eq:Koperators}
	\begin{split}
		K_1 &= \sum_{i=1}^{s}\sigma_i \left(\mv{U}_3^\T \diag{\mv{v}_{i}^{(3)}} \conj{\mv{U}_3} \right) \otimes \left(\mv{U}_2^\T \diag{\mv{v}_{i}^{(2)}} \conj{\mv{U}_2} \right) \otimes \left(\diag{\mv{v}_{i}^{(1)}}\right) \\ 
		K_2 &= \sum_{i=1}^{s}\sigma_i \left(\mv{U}_3^\T \diag{\mv{v}_{i}^{(3)}} \conj{\mv{U}_3} \right) \otimes \left(\mv{U}_1^\T \diag{\mv{v}_{i}^{(1)}} \conj{\mv{U}_1} \right) \otimes \left(\diag{\mv{v}_{i}^{(2)}}\right) \\
		K_3 &= \sum_{i=1}^{s}\sigma_i \left(\mv{U}_2^\T \diag{\mv{v}_{i}^{(2)}} \conj{\mv{U}_2} \right) \otimes \left(\mv{U}_1^\T \diag{\mv{v}_{i}^{(1)}} \conj{\mv{U}_1} \right) \otimes \left(\diag{\mv{v}_{i}^{(3)}}\right).
	\end{split}
\end{equation}

So, we find the following linear system to solve for $\mvec{\conj{\mv{U}_1}\mv{G}_{(1)}}$:
\begin{equation}\label{eq:eqnU1var}
	\left\{
	-\mv{I} \otimes \conj{\mv{D}_{xx}} -
	\left[ 
	\conj{\left(\mv{I} \otimes \mv{U}_2^\H \mv{D}_{yy} \mv{U}_2\right)} - \conj{\left(\mv{U}_3^\H \mv{D}_{zz} \mv{U}_3 \otimes \mv{I}\right)}
	\right] \otimes \mv{I}
	- K_1
	\right\}
	\mvec{\underbrace{\conj{\mv{U}_1}\mv{G}_{(1)}}_{\mv{X}_1}} = \mvec{\mv{F}_{(1)} \left(\mv{U}_3 \otimes \mv{U}_2\right)}.
\end{equation}
Observe this is a square system with $n_1 \times r_2 r_3$ unknowns
(where the solution in matrix form $\mv{X}_1$ is typical for rank $r >
r_1$).  In a similar way, update equations for $\mv{U}_2$ and
$\mv{U}_3$ are derived by multiplying \eqref{eq:start3dvar} with the
other factor matrices in the appropriate directions:
\begin{equation}\label{eq:eqnU2var}
	\left\{
	-\mv{I} \otimes \conj{\mv{D}_{yy}} +
	\left[ -
	\conj{\left(I \otimes \mv{U}_1^\H \mv{D}_{xx} \mv{U}_1\right)} - \conj{\left(\mv{U}_3^\H \mv{D}_{zz} \mv{U}_3 \otimes \mv{I}\right)}
	\right]  \otimes \mv{I}
	- K_2
	\right\}
	\mvec{\underbrace{\conj{\mv{U}_2}\mv{G}_{(2)}}_{\mv{X}_2}} = \mvec{\mv{F}_{(2)} \left(\mv{U}_3 \otimes \mv{U}_1\right)}.
\end{equation}
and
\begin{equation}\label{eq:eqnU3var}
	\left\{-\mv{I} \otimes \conj{\mv{D}_{zz}} +
	\left[-
	\conj{\left(\mv{I} \otimes \mv{U}_1^\H \mv{D}_{xx} \mv{U}_1\right)} - \conj{\left(\mv{U}_2^\H \mv{D}_{yy} \mv{U}_2 \otimes I\right)}
	\right] \otimes I
	- K_3
	\right\} \mvec{\underbrace{\conj{\mv{U}_3}\mv{G}_{(3)}}_{\mv{X}_3}} = \mvec{\mv{F}_{(3)} \left(\mv{U}_2 \otimes \mv{U}_1\right)}.
\end{equation}
Alternating between solving for $\mv{U}_1$, $\mv{U}_2$ and $\mv{U}_3$ using \eqref{eq:eqnU1var}, \eqref{eq:eqnU2var} or \eqref{eq:eqnU3var} results in an algorithm to approximate low-rank tensor solutions for three dimensional problems as given in \eqref{eq:driven3d}. Also in this case the orthogonality of the columns of $\mv{U}_1$, $\mv{U}_2$ and $\mv{U}_3$ are maintained by additional QR factorizations.
So, we derive the algorithm as formulated in Algorithm \ref{alg:var3dv1}. The generalization for dimensions $d > 3$ is straight forward.

\begin{algorithm}[t]
	\SetAlgoLined
	[$\mt{G}, \mv{U}_1, \mv{U}_2, \mv{U}_3$] = hosvd(initial guess)\;
	[$\mv{\Sigma}, \mv{V}_1, \mv{V}_2, \mv{V}_3$] = cp\_als($\mt{K}$)\;
	\While{not converged}{
		\For{i = 1, 2, 3}{
			Compute $K_i$ using \eqref{eq:Koperators}\;
			Solve for $\mv{X}_i = \conj{\mv{U}_i}\mv{G}_{(i)} \in \C^{n_i \times r^{d-1}}$ using \eqref{eq:eqnU1var}, \eqref{eq:eqnU2var} or \eqref{eq:eqnU3var}\;
			$\conj{\mv{U}_{i}} \mv{G}_{(i)} = \qr{\mv{X}_i(:,~ 1:r_i), 0}$;
		}
	}
	$\mt{G} = \texttt{reconstruct} \left[\mv{G}_{(i)}, ~i\right]$\;
	$\mt{M} = \mt{G} \times_1 \mv{U}_1  \times_2 \mv{U}_2  \times_3 \mv{U}_3$\;
	\caption{Solve for the low-rank tensor decomposition of the solution $\mt{M}$ of a 3D Helmholtz problem with space-dependent wave number (version 1).}
	\label{alg:var3dv1}
\end{algorithm}

Similar to the discussion for the constant wave number algorithms,
observe that we solve again for a large matrix $\mv{X}_i \in \C^{n_i \times	r_1r_2r_3/r_i}$.
So, in general the rank of this matrix could be
$\min\left(n_i,~ r_1r_2r_3/r_i\right)$. But it is also known that
$\mv{X}_i = \conj{\mv{U}_i}\mv{G}_{(i)}$ leads to the fact that
the rank of $\mv{X}_i$ should be at most $r_i$.  So selecting the
first $r_i$ columns of $\mv{X}_i$ and computing its QR decomposition
is sufficient to derive a new orthonormal basis as factor matrix
$\conj{\mv{U}}_i$.

Algorithm \ref{alg:var3dv1} is exactly the space-dependent wave number
equivalent of Algorithm \ref{alg:const3dv1}. The same ideas can be applied
to derive space-dependent wave number alternatives of the algorithms
corresponding to version 2 and 3.
Again, to circumvent solving large systems, we can pre-compute the 
QR factorization of $\mv{G}_{(i)}$ and project these equations onto 
the obtained $\mv{Q}_i$. Indeed, this will reduce the number of 
unknowns in these linear systems to exactly the number of unknowns 
as needed for the factor matrices $\conj{\mv{U}}_1$ and
$\conj{\mv{U}}_2$.

Let us discuss the details. We start again from equation \eqref{eq:firstUnfoldingVar} and use the QR
factorization of $\mv{G}_{(1)}^\H$, $\mv{Q}_1 \mv{R}_1^\H = \qr{\mv{G}^\H_{(1)}}$. This yields
\begin{equation*}
	-\conj{\mv{D}_{xx}\mv{U}_1} \mv{R}_1\mv{Q}_1^\H
	- \conj{\mv{U}_1} \mv{R}_1\mv{Q}_1^\H \left(I \otimes \mv{U}_2^\H \mv{D}_{yy} \mv{U}_2\right)^\H
	- \conj{\mv{U}_1} \mv{R}_1\mv{Q}_1^\H \left(\mv{U}_3^\H \mv{D}_{zz} \mv{U}_3  \otimes I\right)^\H
	- \left[\mv{K}_{(1)} \circ \conj{\mv{U}_1}\mv{R}_1\mv{Q}_1^\H \left( \mv{U}_3 \otimes \mv{U}_2 \right)^\H \right] \left(\mv{U}_3 \otimes \mv{U}_2 \right)
	=  \mv{F}_{(1)} \left(\mv{U}_3 \otimes \mv{U}_2\right).
\end{equation*}
Post multiplication of the left hand side of this equation by $\mv{Q}_1$ yields
\begin{equation*}
	-\conj{\mv{D}_{xx}\mv{U}_1} \mv{R}_1
	- \conj{\mv{U}_1} \mv{R}_1\mv{Q}_1^\H \left(I \otimes \mv{U}_2^\H \mv{D}_{yy} \mv{U}_2\right)^\H\mv{Q}_1
	- \conj{\mv{U}_1} \mv{R}_1\mv{Q}_1^\H \left(\mv{U}_3^\H \mv{D}_{zz} \mv{U}_3  \otimes I\right)^\H\mv{Q}_1
	- \left[\mv{K}_{(1)} \circ \conj{\mv{U}_1}\mv{R}_1\mv{Q}_1^\H \left( \mv{U}_3 \otimes \mv{U}_2 \right)^\H \right] \left(\mv{U}_3 \otimes \mv{U}_2 \right)\mv{Q}_1.
\end{equation*}
To solve this equation for $\mv{U}_1$, it is written in vectorized form as
\begin{equation}\label{eq:eqnU1v2forvar}
	\left\{
	-\mv{I} \otimes \conj{\mv{D}_{xx}} +
	\mv{Q}_1^\T\left[- 
	\conj{\left(\mv{I} \otimes \mv{U}_2^\H \mv{D}_{yy} \mv{U}_2\right)} - \conj{\left(\mv{U}_3^\H \mv{D}_{zz} \mv{U}_3 \otimes \mv{I} \right)}
	\right]\conj{\mv{Q}_1} \otimes \mv{I}
	-\mv{Q}_1^\T K_1 \conj{\mv{Q}_1}
	\right\} \mvec{\underbrace{\conj{\mv{U}_1}\mv{R}_1}_{\mv{X}_1}} = \mvec{\mv{F}_{(1)} \left(\mv{U}_3 \otimes \mv{U}_2\right) \mv{Q}_1}.
\end{equation}
In a similar way, the update equations for $\mv{U}_2$ and $\mv{U}_3$ are
derived by multiplying \eqref{eq:start3dvar} with the other factor
matrices in the appropriate dimensions and using the QR factorizations
of $\mv{G}_{(i)}^\H$:
\begin{equation}\label{eq:eqnU2v2forvar}
	\left\{
	-\mv{I} \otimes \conj{\mv{D}_{yy}} +
	\mv{Q}_2^\T\left[ 
	-\conj{\left(\mv{I} \otimes \mv{U}_1^\H \mv{D}_{xx} \mv{U}_1\right)} - \conj{\left(\mv{U}_3^\H \mv{D}_{zz} \mv{U}_3 \otimes \mv{I} \right)}
	\right]\conj{\mv{Q}_2}  \otimes \mv{I}
	-\mv{Q}_2^\T K_2 \conj{\mv{Q}_2}
	\right\} \mvec{\underbrace{{\conj{\mv{U}_2}\mv{R}_2}}_{\mv{X}_2}} 
	= \mvec{\mv{F}_{(2)} \left(\mv{U}_3 \otimes \mv{U}_1\right) \mv{Q}_2}.
\end{equation}
and
\begin{equation}\label{eq:eqnU3v2forvar}
	\left\{
	-\mv{I} \otimes \conj{\mv{D}_{zz}} +
	\mv{Q}_3^\T\left[ 
	-\conj{\left(\mv{I} \otimes \mv{U}_1^\H \mv{D}_{xx} \mv{U}_1\right)} - \conj{\left(\mv{U}_2^\H \mv{D}_{yy} \mv{U}_2 \otimes \mv{I} \right)}
	\right]\conj{\mv{Q}_3} \otimes \mv{I}
	\right\} \mvec{\underbrace{\conj{\mv{U}_3}\mv{R}_3}_{\mv{X}_3}} 
	-\mv{Q}_3^\T K_3 \conj{\mv{Q}_3}
	= \mvec{\mv{F}_{(3)} \left(\mv{U}_2 \otimes \mv{U}_1\right) \mv{Q}_3}.
\end{equation}

All these equations are cheap to solve. Indeed,
$\mvec{\conj{\mv{U}_i}\mv{R}_i}$ has length $n_i r_i$. Computing a
symmetric reverse Cuthill-McKee permutation of the system matrix one
observes a matrix with a bandwidth $\Oh{r}$, so solving these
equations has a computational cost $\Oh{nr^2}$.

Alternating between solving for $\mv{U}_1$, $\mv{U}_2$ and $\mv{U}_3$ using \eqref{eq:eqnU1v2forvar}, \eqref{eq:eqnU2v2forvar} or \eqref{eq:eqnU3var} results again in an algorithm to approximate low-rank solutions for three dimensional space-dependent Helmholtz problems. Also in this case the orthogonality of the columns of $\mv{U}_1$, $\mv{U}_2$ and $\mv{U}_3$ are maintained by additional QR factorizations.
So, we derive the algorithm as formulated in Algorithm \ref{alg:var3dv3}.
Algorithm \ref{alg:var3dv3} is exactly the space-dependent wave number equivalent of Algorithm \ref{alg:const3dv3}.

\begin{algorithm}[t]
	\SetAlgoLined
	[$\mt{G}, \mv{U}_1, \mv{U}_2, \mv{U}_3$] = hosvd(initial guess)\;
	[$\mv{\Sigma}, \mv{V}_1, \mv{V}_2, \mv{V}_3$] = cp\_als($\mt{K}$)\;
	\While{not converged}{
		\For{i = 1, 2}{
			Compute $K_i$ using \eqref{eq:Koperators}\;
			$\mv{Q}_i \widetilde{\mv{R}} = \qr{\mv{G}_{(i)}^\H, 0}$\;
			Solve for $\mv{X}_i = \conj{\mv{U}_i}\mv{R}_i \in \C^{n_i \times r_i}$ using \eqref{eq:eqnU1v2forvar} or \eqref{eq:eqnU2v2forvar}\;
			$\conj{\mv{U}_{i}} \mv{R}_i = \qr{\mv{X}_i, 0}$\;
			$\mt{G} = \texttt{reconstruct} \left[\mv{R}_i\mv{Q}_i^\H, ~i\right]$;
		}
		Compute $K_3$ using \eqref{eq:Koperators}\;
		Solve for $\mv{X}_3 = \conj{\mv{U}_3}\mv{G}_{(3)} \in \C^{n_3 \times r^{d-1}}$ using \eqref{eq:eqnU3var}\;
		$\conj{\mv{U}_{3}} \mv{G}_{(3)} = \qr{\mv{X}_3, 0}$\;
		$\mt{G} = \texttt{reconstruct} \left[\mv{G}_{(3)},~3\right]$\;
	}
	$\mt{M} = \mt{G} \times_1 \mv{U}_1  \times_2 \mv{U}_2  \times_3 \mv{U}_3$\;
	\caption{Solve for the low-rank tensor decomposition of the solution $\mt{M}$ of a 3D Helmholtz problem with space-dependent wave number (version 3).}
	\label{alg:var3dv3}
\end{algorithm}

\subsection{Projection operator for space-dependent wave number}
Consider a tensor $\mt{M}$ in Tucker tensor format and factorized as
$\mt{M} = \mt{G} \times_1 \mv{U}_1 \times_2 \mv{U}_2 \times_3 \mv{U}_3$,
with unknowns $\mt{G}, \mv{U}_1, \mv{U}_2$ and $\mv{U}_3$.
Discretization of \eqref{eq:driven3d} with a space-dependent wave number 
leads to a linear operator $\mathcal{L}$ applied on tensors. Its matrix
representation $\mv{L}$ has again a structure as given in
\eqref{eq:operatorLasSumKroneckerProducts}.

Solving for the unknown factors $\mv{U}_1$, $\mv{U}_2$ or $\mv{U}_3$
(and the core-tensor $\mt{G}$) using \eqref{eq:eqnU1var},
\eqref{eq:eqnU2var} or \eqref{eq:eqnU3var} can, again, be interpreted as a
projection operator applied on the residual. For example,
\eqref{eq:eqnU1var} can be interpreted as
\begin{equation}
	\conj{\left(\mv{U}_3^\H \otimes \mv{U}_2^\H \otimes \mv{I} \right) \mv{L} \left(\mv{U}_3 \otimes \mv{U}_2 \otimes \mv{I}\right)} \mvec{\conj{\mv{U}_1}\mv{G}_{(1)}} = \left( \mv{U}_3^\T \otimes \mv{U}_2^\T \otimes \mv{I} \right) \mvec{\mv{F}_{(1)}},
\end{equation}
The residual in tensor format is given by
\begin{equation}
	\begin{split}
		\mt{R} = \mt{F} &- \mathcal{L}\mt{M}, \\
		= \mt{F} 
		&+ \mt{G} \times_1 \mv{D}_{xx} \mv{U}_1 \times_2 \mv{U}_2 \times_3 \mv{U}_3 \\
		&+ \mt{G} \times_1 \mv{U}_1 \times_2 \mv{D}_{yy}\mv{U}_2 \times_3 \mv{U}_3 \\
		&+ \mt{G} \times_1 \mv{U}_1 \times_2 \mv{U}_2 \times_3 \mv{D}_{zz}\mv{U}_3 \\
		&+ \mt{K} \circ \left(\mt{G} \times_1 \mv{U}_1 \times_2 \mv{U}_2 \times_3 \mv{U}_3 \right).
	\end{split}
\end{equation}

Writing this tensor equation in the first unfolding leads to the following matrix equation
\begin{equation}
	\begin{split}
	\mv{R}_{(1)} = \mv{F}_{(1)} 
	&+ \conj{\mv{D}_{xx}\mv{U}_1}\mv{G}_{(1)} \left(\mv{U}_3 \otimes \mv{U}_2 \right)^\H \\
	&+ \conj{\mv{U}_1}\mv{G}_{(1)} \left(\mv{U}_3 \otimes \mv{D}_{yy}\mv{U}_2 \right)^\H \\
	&+ \conj{\mv{U}_1}\mv{G}_{(1)} \left(\mv{D}_{zz}\mv{U}_3 \otimes \mv{U}_2 \right)^\H \\
	&+ \mv{K}_{(1)} \circ \left( \conj{\mv{U}_1} \mv{G}_{(1)} \left( \mv{U}_3 \otimes \mv{U}_2 \right)^\H \right)
	\end{split}
\end{equation}
which can be matricized as
\begin{equation*}
	\begin{split}
		\mvec{\mv{R}_{(1)}} &= \mvec{\mv{F}_{(1)}} \\
		&- \left(-\conj{\left(\mv{U}_3 \otimes \mv{U}_2 \otimes \mv{D}_{xx}\right)} 
		- \conj{\left(\mv{U}_3 \otimes \mv{D}_{yy}\mv{U}_2 \otimes \mv{I}\right)} 
		- \conj{\left(\mv{D}_{zz}\mv{U}_3 \otimes \mv{U}_2 \otimes \mv{I}\right)} 
		- \diag{\mvec{\mv{K}_{(1)}}} \conj{\left(\mv{U}_3 \otimes 	\mv{U}_2 \otimes \mv{I} \right)}\right)\mvec{\conj{\mv{U}_1}\mv{G}_{(1)}}.
	\end{split}
\end{equation*}
Rewriting this results in exactly the same structure and projection operator as in the constant wave number case:
\begin{equation}
	\begin{split}
		\mvec{\mv{R}_{(1)}} &= \ldots \\
		&= \mvec{\mv{F}_{(1)}} - \conj{\mv{L} \left(\mv{U}_3 \otimes \mv{U}_2 \otimes \mv{I}\right)}\mvec{\conj{\mv{U}_1}\mv{G}_{(1)}} \\
		&= \mvec{\mv{F}_{(1)}} - \conj{\mv{L} \left(\mv{U}_3 \otimes \mv{U}_2 \otimes \mv{I}\right)} \left[\conj{\left(\mv{U}_3^\H \otimes \mv{U}_2^\H \otimes \mv{I} \right) \mv{L} \left(\mv{U}_3 \otimes \mv{U}_2 \otimes \mv{I}\right)}\right]^{-1}\left( \mv{U}_3^\T \otimes \mv{U}_2^\T \otimes \mv{I} \right) \mvec{\mv{F}_{(1)}} \\
		&= P_{23}\mvec{\mv{F}_{(1)}}
	\end{split}
\end{equation}
where projection operator $P_{23}$ is similar to the projector in the constant wave number case, see \eqref{eq:projector3dconst}, and now given by
\begin{equation}
	\begin{split}
		P_{23} &= \mv{I} - \conj{\mv{L} \left(\mv{U}_3 \otimes \mv{U}_2 \otimes \mv{I}\right)} \left[\conj{\left(\mv{U}_3^\H \otimes \mv{U}_2^\H \otimes \mv{I} \right) \mv{L} \left(\mv{U}_3 \otimes \mv{U}_2 \otimes \mv{I}\right)}\right]^{-1}\left( \mv{U}_3^\T \otimes \mv{U}_2^\T \otimes \mv{I} \right) \\
		&= \mv{I} - \mv{X}.
	\end{split}
\end{equation}

A similar derivation results in projection operators $P_{13}$ and $P_{12}$ for the updates in $\mv{U}_2$ and $\mv{U}_3$, respectively. Both are also the same as in the constant wave number case, as given in \eqref{eq:projectors3dconst}.

\section{Numerical results} \label{sec:NumericalResults}
In this section, we demonstrate the promising results of the derived
algorithms with some numerical experiments in two and three
dimensions. Furthermore, we consider discretizations of the Helmholtz
equation with constant and space-dependent wave numbers.

\subsection{2D Helmholtz problem with space-dependent wave number}
First, we consider a 2D Helmholtz problem with a space-dependent wave 
number given by $k^2(x,y) = 2 + e^{-x^2 -y^2}$.

For this example the two dimensional domain $\Omega = [-10, ~10]^2$ is discretized with $M=1000$ equidistant mesh points per direction in the interior of the domain. Further it is extended with exterior complex scaling to implement the absorbing boundary conditions. In total, the number of discretization points per directions equals $n = n_1 = n_2 = 1668$. As external force $f(x, y) = -e^{-x^2-y^2}$ is applied.

In this space-dependent wave number example it is known that the matrix representation of the semi-exact solution of the Helmholtz equation on the full grid has a low rank. Indeed, approximating the semi-exact solution with a low-rank matrix with rank $r=17$ is in this case sufficient to obtain an error below the threshold $\tau = 10^{-6}$.

\begin{figure}
  \centering
  \begin{tabular}{cc}
	\input{error-variable-2d-r=17.tikz} & 
	\input{singularvalues-variable-2d.tikz}
  \end{tabular}
        \caption{Plot of error and residual (left) and singular values (right) per iteration for space-dependent wave number in 2D Helmholtz problem ($M=1000, r=17$).}
        \label{fig:error+singularvalues-var2d-r=17}
\end{figure}
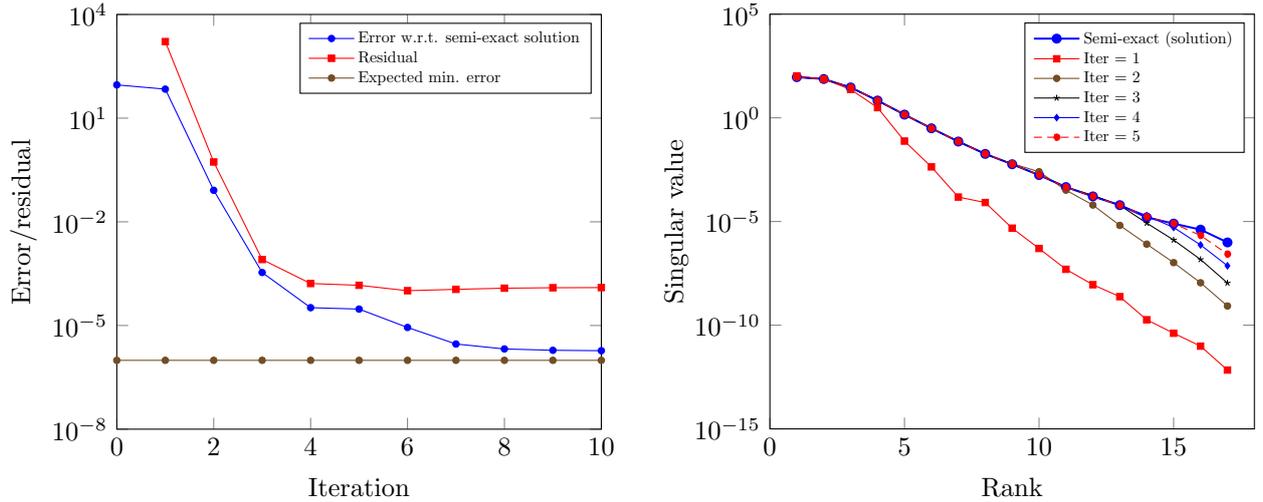

Starting with an random (orthonormalized) initial guess for $\mv{V}^{(0)} \in \C^{n \times r}$ only a small number of iterations of Algorithm \ref{alg:2d} is needed to obtain an error similar to the specified threshold $\tau$. As shown in Figure~\ref{fig:error+singularvalues-var2d-r=17} both the residual and the error with respect to the semi-exact solution decay in only a few iterations (i.e. in this example 4-8 iterations) to a level almost similar to the expected tolerance.

The singular values of the approximation $\mv{A}^{(k)} =
\mv{U}^{(k)}{\mv{R}^{(k)}}^\H {\mv{V}^{(k)}}^\H$ in iteration $i$ can
be computed and are shown for increasing iterations in
Figure~\ref{fig:error+singularvalues-var2d-r=17}. As expected the low-
rank approximations recover the singular values of the full grid
semi-exact solution. In fact $\mv{R}^{(k)}$ converges towards
$\text{diag}(\sigma_i)$.

The numerical rank of the matrix representation of the solution of a
Helmholtz problem with a space-dependent wave number is unknown in
advance. But, the presented algorithm is stable with respect to over-
and underestimation of the numerical rank of the solution. Figure
\ref{fig:error+residuals-var2d-M=1000} shows both the error and
residual per iteration and illustrates this statement by approximating
the same semi-exact solution with increasing ranks $r \in \{ 12, 18, 24, 36 \}$.

\begin{figure}
  \centering
 \begin{tabular}{cc}
	\input{error-variable-2d-M=1000.tikz} &
	\input{residuals-variable-2d-M=1000.tikz}
 \end{tabular}
  \caption{Plot of errors (left) and residuals (right) per iteration for space-dependent wave number in 2D Helmholtz problem with increasing ranks ($M=1000$)}
  \label{fig:error+residuals-var2d-M=1000}
\end{figure}
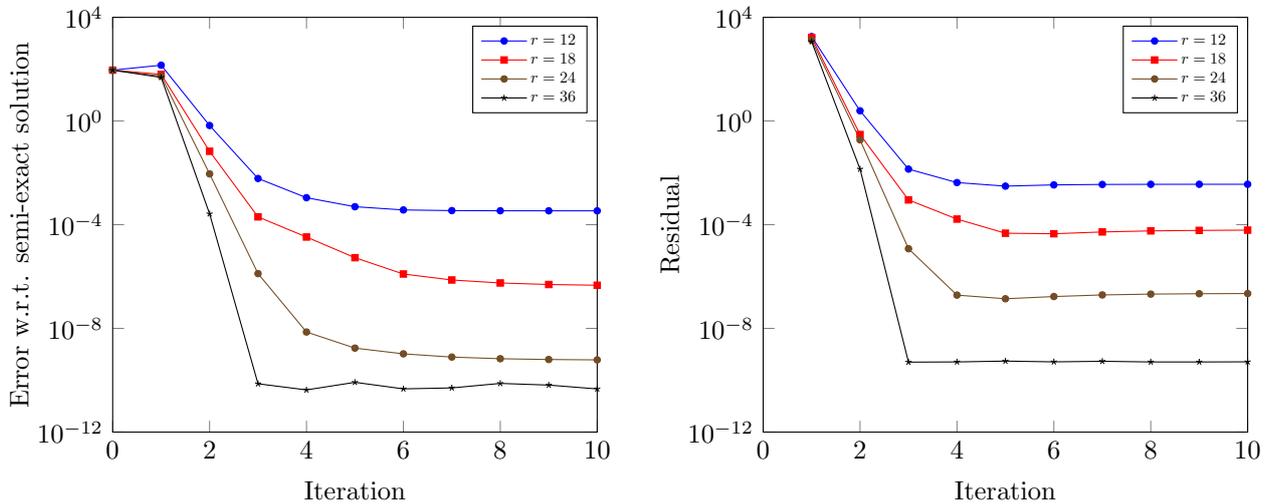

\begin{figure}
	\centering
 \begin{tabular}{cc}
	\input{runtime-variable-2d-M=1000.tikz} &
	\input{runtime-variable-2d-M=1000+loglog.tikz}
 \end{tabular}
        \caption{Plot of runtime for 10 iterations for space-dependent wave number in 2D Helmholtz problem with increasing ranks ($M=1000$). Left: linlin-scale, right: loglog-scale.}
        \label{fig:runtime-variable-M=1000}
\end{figure}
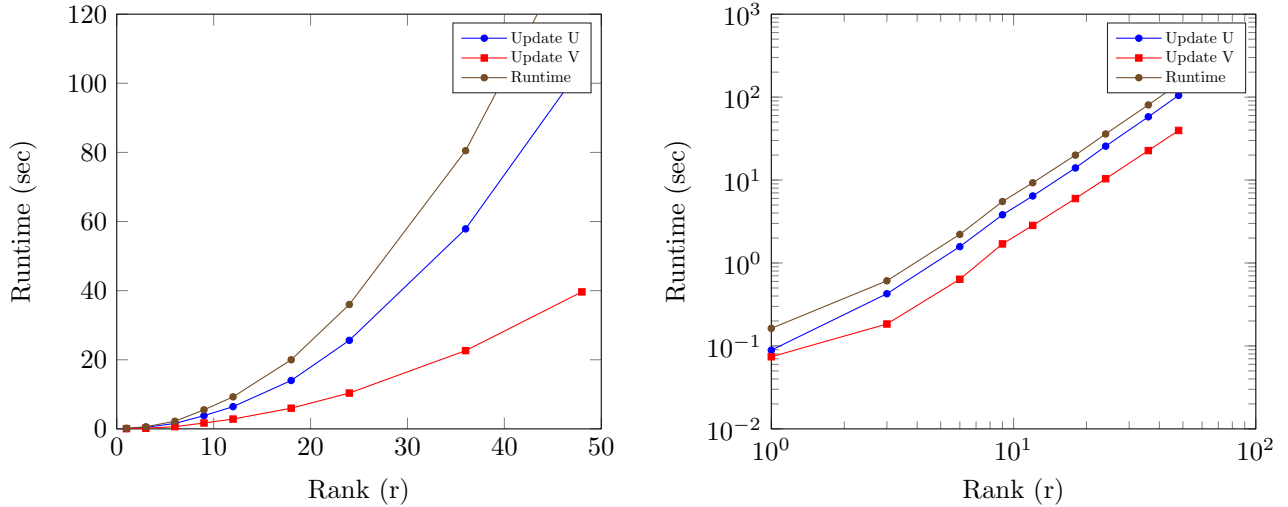

In contrast to the constant coefficient wave number case the
convergence with space-dependent wave number depends also
on the maximal attainable rank. For increasing maximal attainable
ranks the number of needed iterations decreases. This is especially
observed when the error is considered, but it can also be seen in the
figure where the residuals are shown, Fig.~\ref{fig:error+residuals-var2d-M=1000}.

\subsection{3D Helmholtz problem with space-dependent wave number}
In this example we solve a 3D Helmholtz problem with a space-dependent
wave number discretized on a DVR-grid \cite{rescigno2000numerical}. All
three versions of the 3D algorithm for space-dependent wave numbers
can successfully be applied.

\begin{figure}
	\centering
	\begin{subfigure}{.45\textwidth}
		\centering
		\input{cprank-orderDvr=7-wavenumber.tikz} \\
		\caption{CP-rank of space-dependent wave number for 3D Helmholtz problem.}
		\label{fig:cprank_wavenumber}
	\end{subfigure}
	\quad
	\begin{subfigure}{.45\textwidth}
		\centering
		\input{residual-all-orderDvr=7-solrank=16-waverank=32.tikz} \\
		\caption{Residual per iteration ($r = 16, s = 32$).}
		\label{fig:residual-all-orderDvr=7-solrank=16-waverank=32}
	\end{subfigure}
	\caption{Low rank approximation to space-dependent wave number and residuals for version 1, version 2 and version 3 of 3D Helmholtz problem with space-dependent wave number (orderDvr = 7).}
	\label{fig:residual-all-orderDvr=7}
\end{figure}
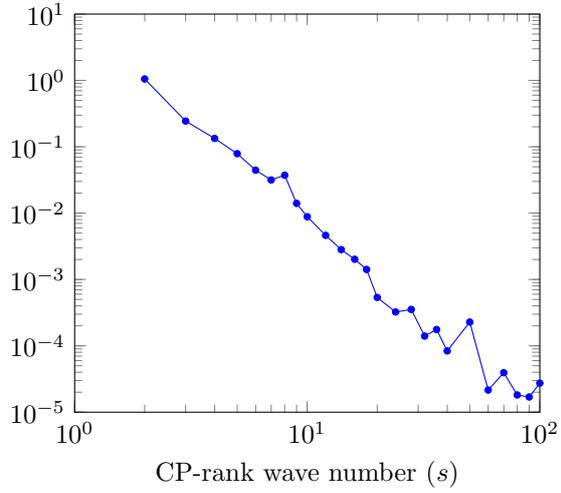
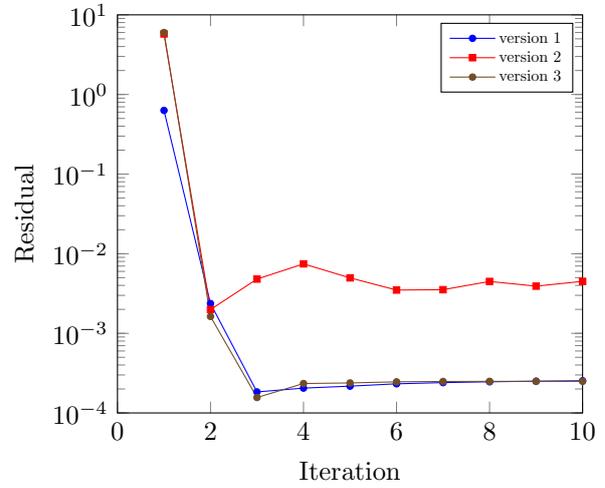

First, to reduce computational cost of construction of the operators $K_1$, $K_2$ and $K_3$, see \eqref{eq:Koperators}, a CP-decomposition of the space-dependent wave number is constructed.
As shown in Figure \ref{fig:cprank_wavenumber} the space-dependent wave number can be well-approximated by a small number of rank-1 tensors. For the examples discussed in this section we used a CP-rank $s = 32$ to approximate this space-dependent wave number. Hence, the error in approximating the wave number is approximately $\Oh{10^{-4}}$.

For all three versions of the algorithm we use 10 iterations of the algorithm to converge to the low-rank solution. For example if we compute the low-rank solution (with $r = r_x = r_y = r_z = 16$) the residual after each iteration for all algorithms is shown in Figure \ref{fig:residual-all-orderDvr=7-solrank=16-waverank=32}.

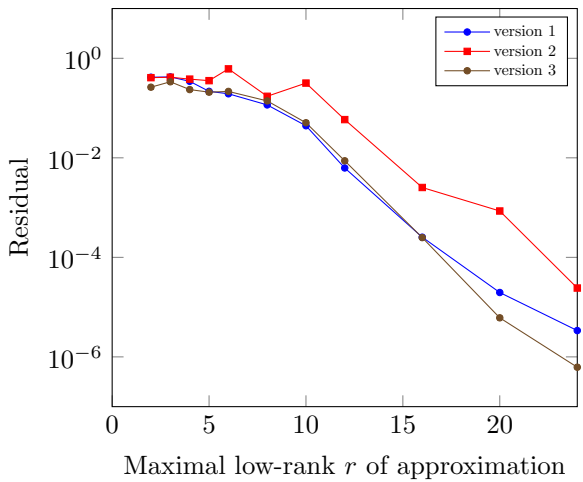
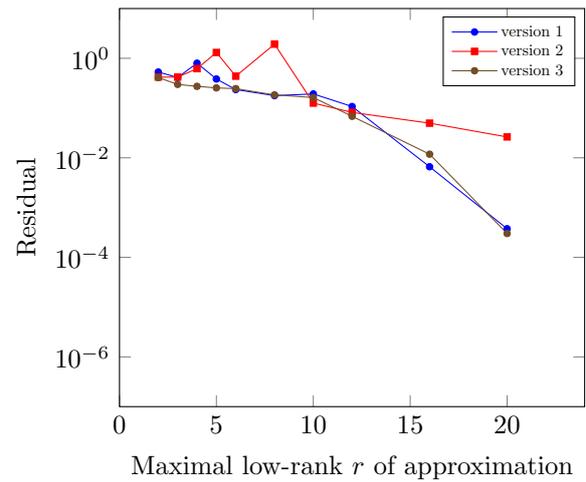
\begin{figure}
	\centering
	\begin{subfigure}{.45\textwidth}
		\centering
		\input{residual-all-orderDvr=7-waverank=32.tikz} \\
		\caption{orderDvr = 7}
		\label{fig:residual-all-orderDvr=7-waverank=32}
	\end{subfigure}
	\quad
	\begin{subfigure}{.45\textwidth}
		\centering
		\input{residual-all-orderDvr=14-waverank=32.tikz} \\
		\caption{orderDvr = 14}
		\label{fig:residual-all-orderDvr=14-waverank=32}
	\end{subfigure}
	\caption{Residual after iteration 10 iterations for all three versions of algorithm with $s = 32$.}
	\label{fig:residual-all-waverank=32}
\end{figure}

If we increase the maximal attainable rank $r$ of the low-rank approximation, indeed the residual decreases as shown in Figure \ref{fig:residual-all-orderDvr=7-waverank=32}. The residual for version 1 and version 3 are good, while version 2 cannot compete with both other versions by reducing the residual as far as the other versions. Therefore version 1 or version 3, as given in Algorithm~\ref{alg:var3dv1} or Algorithm~\ref{alg:var3dv3} are preferred.

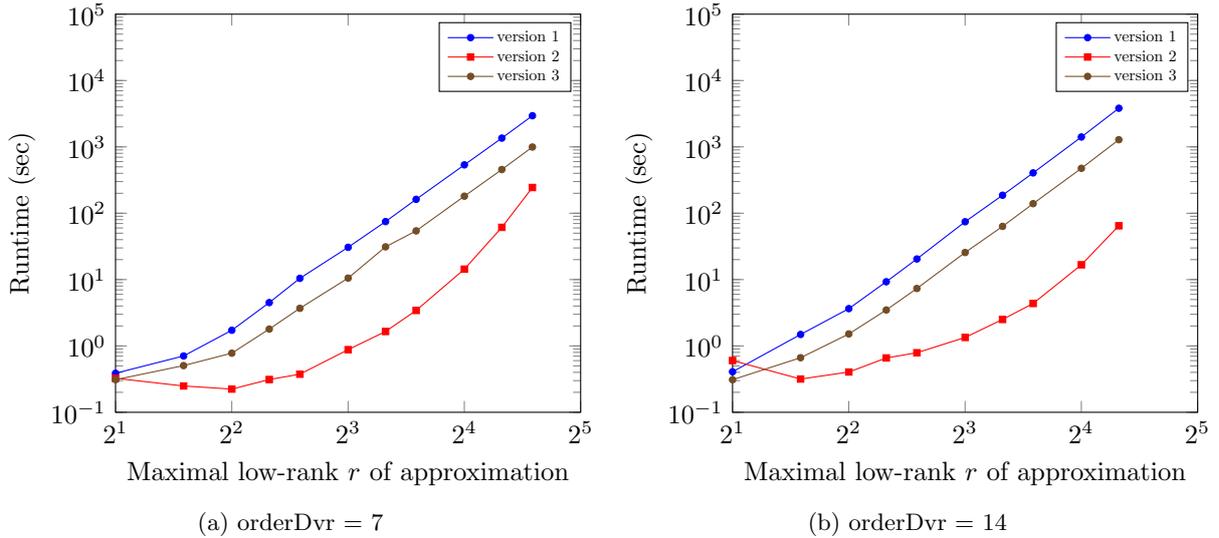
\begin{figure}
	\centering
	\begin{subfigure}{.45\textwidth}
		\centering
		\input{runtimes-all-orderDvr=7-waverank=32+loglog.tikz} \\
		\caption{orderDvr = 7}
		\label{fig:runtimes-all-orderDvr=7-waverank=32+loglog}
	\end{subfigure}
	\quad
	\begin{subfigure}{.45\textwidth}
		\centering
		\input{runtimes-all-orderDvr=14-waverank=32+loglog.tikz} \\
		\caption{orderDvr = 14}
		\label{fig:runtimes-all-orderDvr=14-waverank=32+loglog}
	\end{subfigure}
	\caption{Runtime of 10 iterations for all three versions of algorithm for 3D Helmholtz with space-dependent wave number of low-rank $s = 32$.}
	\label{fig:runtimes-all-waverank=32}
\end{figure}

Considering the runtimes of the three versions, similar results as before are observed. In this experiment with \mbox{orderDvr=7} the number of gridpoints equals to $n = 41$. For version 1 and 3, again a runtime of $\Oh{nr^4}$ is observed. The runtime for version 2 splits into two parts: $\Oh{nr^2 + r^9}$. Due to the small rank $r$ and the large number of iterations in these examples algorithm 2 is the fastest version. The runtimes for version 1 and version 3 differ indeed approximately a factor $d$, which makes version 3 better then version 1. The runtimes with orderDvr = 7 (i.e. $n=41$) are shown in Figure \ref{fig:runtimes-all-orderDvr=7-waverank=32+loglog} and with orderDvr = 14 (i.e. $n=90$) are shown in Figure \ref{fig:runtimes-all-orderDvr=14-waverank=32+loglog}.

Comparing the runtimes for orderDvr=7 and orderDvr=14 we see for version 2 (when the rank gets larger) indeed approximately the same runtime independent of orderDvr. Also versions 1 and 3 consume approximately twice as much time which is as expected by the linear dependence on $n$ for both algorithms.

\begin{figure}
	\centering
	\begin{subfigure}{\textwidth}
		\centering
		\includegraphics[width=0.85\textwidth]{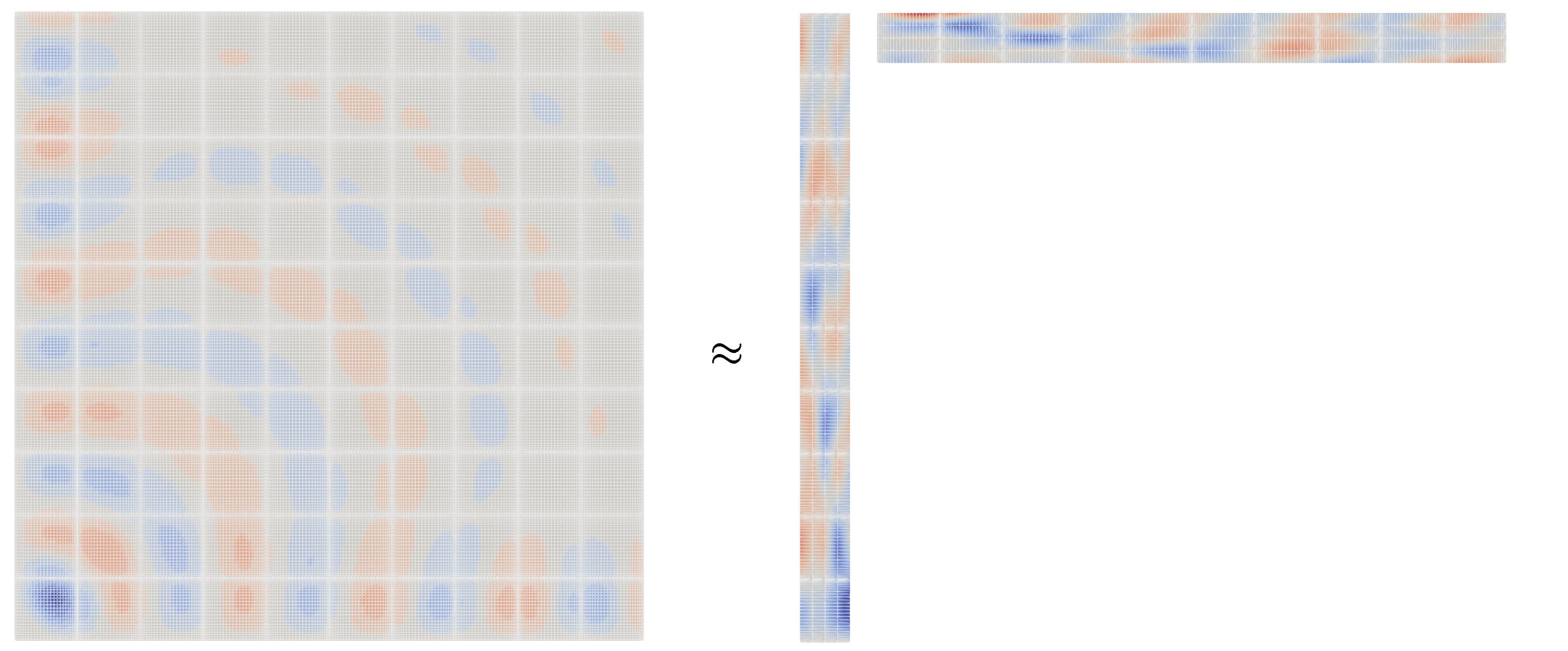}
		\caption{Impression of a low-rank matrix approximation in 2D.}
		\label{fig:impression2d}
	\end{subfigure}
	\par
	\begin{subfigure}{\textwidth}
		\centering
		\includegraphics[width=0.9\textwidth]{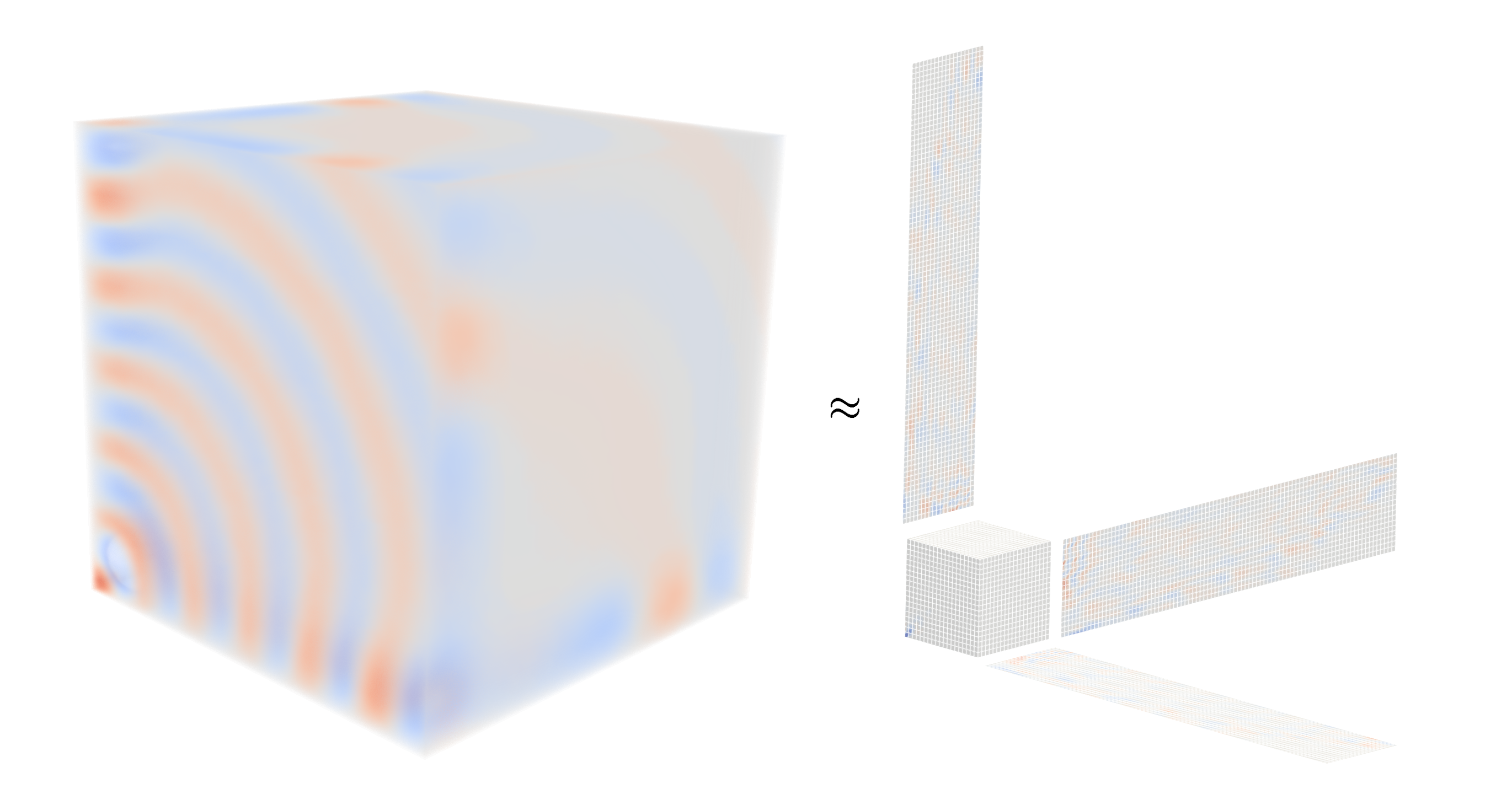}
		\caption{Impression of a low-rank tensor approximation in 3D.}
		\label{fig:impression3d}
	\end{subfigure}
	\caption{Impressions of a low-rank approximation of a matrix and a Tucker tensor representing the wave function as solution to a 2D and 3D Helmholtz problem with a space-dependent wave number.}
\end{figure}

\begin{figure}
	\centering
	\includegraphics[width=0.9\textwidth]{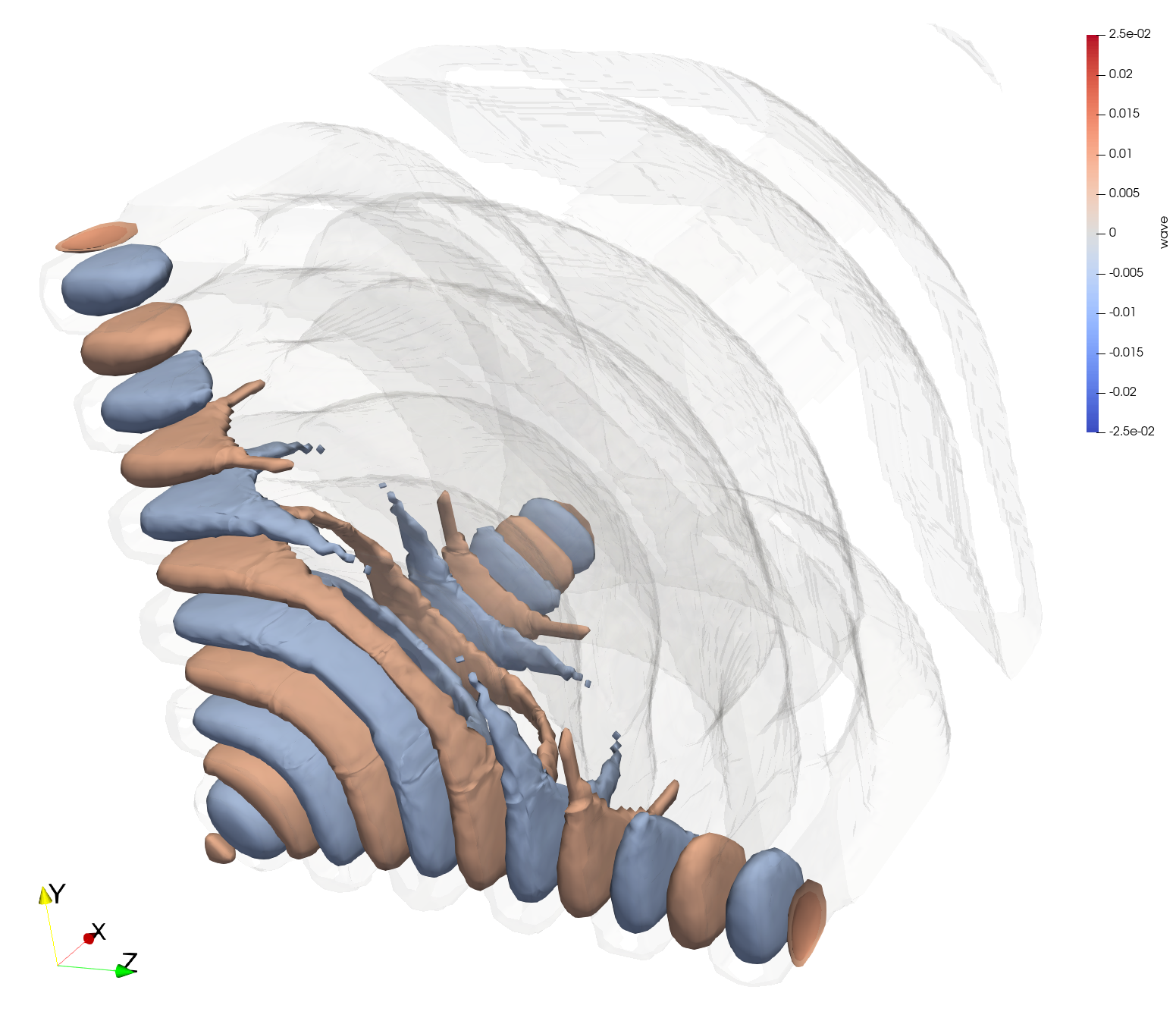}
	\caption{Visualization of a 3D wave as low-rank approximation
          to a 3D Helmholtz problem with space-dependent wave number
          with single, double and triple ionization.}
	\label{fig:impression3dv2}
\end{figure}

An impression of the low-rank approximation to the wave function is shown in Figures \ref{fig:impression3d} and \ref{fig:impression3dv2}. In this impression the single, double and triple ionization are visible and can be represented by a low-rank wave function.

\section{Discussion and conclusions}\label{sec:discussion}
In this paper we have analyzed the scattering solutions of a driven
Schr\"odinger equations. These describe a  break-up reaction where a
quantum system is fragmented into multiple fragments.  These problems
are equivalent to solving a Helmholtz equations with space-dependent wave
numbers.

We have shown, first in 2D and then in 3D, that the wave function of
multiple ionization can be well approximated by a low-rank solution.
In 2D, the waves can be represented as a product of two low-rank 
matrices and 3D as a low-rank Tucker tensor decomposition.

We propose a method that determines these low-rank components
of the solution directly.  We write the solution as a product of low-rank
components and assume that a guess for all but one component is given.
We then write a linear system for the remaining unknown component.
This is the repeated until each of the components is updated.

This procedure can be interpreted as a series of projections of the
residual on a subspaces and a correction within that subspace.

In theory, the generalization for dimensions $d > 3$ is straightforward. But
for dimensions $d>3$ it starts to be beneficial to change to a Tensor
Train factorization \cite{oseledets2011tensor}. It is expected that similar 
strategy can also be applied to tensors in Tensor Train format.

As demonstrated by the numerical experiments, the presented algorithms
are able to exploit the low-rank structure of the solutions. This
gives the advantage to reduce the number of unknowns and shorten the
computational time to solve the Helmholtz equation.

In two dimensions, the low-rank representation of the solution can be
represented by only $2nr$ unknowns instead of the full grid of
$n^2$ unknowns. Also the linear systems to solve per iteration have
only $nr$ unknowns.

In high-dimensional Helmholtz equations, the low-rank Tucker tensor
decomposition represents the solution with $\Oh{r^d + dnr}$
unknowns. So, the total number of unknowns is reduced, but it is still
exponential in the dimension $d$. For increasing dimensions this leads,
again, to systems with a number of unknowns exponential in $d$. Maybe
other Tucker-like tensor decompositions with a number of unknowns only
polynomial in $d$ can resolve this problem and make the presented
algorithm also applicable for higher dimensions.

\bibliographystyle{plain}
\bibliography{refs}

\end{document}

%% file: wave-svd-M=1000-E=2.tikz
%
%
%
\begin{tikzpicture}

\begin{axis}[%
width=0.35\textwidth,
scale only axis,
xmin=0,
xmax=200,
xlabel={Index $i$},
ymode=log,
ymin=1e-9,
ymax=1e2,
ylabel={Singular value $\sigma_i$}
]
\addplot [color=blue]
  table[row sep=crcr]{%
1	25.2391069768743\\
2	14.7486597425292\\
3	7.08038216723123\\
4	3.19201010434919\\
5	1.34076278785827\\
6	0.559522742971098\\
7	0.290571962278091\\
8	0.184311375509482\\
9	0.106086794013317\\
10	0.0763490835483271\\
11	0.0511304392483243\\
12	0.0345097877457911\\
13	0.0245876245470376\\
14	0.0181070650760933\\
15	0.0137632148741801\\
16	0.0106539974863224\\
17	0.00841763886975532\\
18	0.00670969606517275\\
19	0.00541477734480655\\
20	0.00439822028018232\\
21	0.0036126558714963\\
22	0.00299347654417121\\
23	0.00250951018012138\\
24	0.00212028307507985\\
25	0.00180660998864443\\
26	0.00154477197171625\\
27	0.0013274195423419\\
28	0.00114447424206621\\
29	0.00099211368071666\\
30	0.000864768727319691\\
31	0.000757556325529342\\
32	0.000666997492701201\\
33	0.000588657786976922\\
34	0.00052158321004709\\
35	0.000463065775731398\\
36	0.000413082642842699\\
37	0.000369499069688016\\
38	0.000332066477970409\\
39	0.000298997576821309\\
40	0.00027009601299721\\
41	0.000244284547736034\\
42	0.000221616053479412\\
43	0.000201413527822392\\
44	0.000183657204394097\\
45	0.000167754751466196\\
46	0.000153591850622663\\
47	0.000140761273944052\\
48	0.000129231565857336\\
49	0.000118806254668186\\
50	0.000109454409739157\\
51	0.000101018406875196\\
52	9.33956751257642e-05\\
53	8.64536344708777e-05\\
54	8.01103552994468e-05\\
55	7.43233185831988e-05\\
56	6.90442729094999e-05\\
57	6.42530447122629e-05\\
58	5.98738733922753e-05\\
59	5.58743523847707e-05\\
60	5.21769465307176e-05\\
61	4.87801947117721e-05\\
62	4.56388839837367e-05\\
63	4.27672491667743e-05\\
64	4.01177067036462e-05\\
65	3.7689087699374e-05\\
66	3.54271488018332e-05\\
67	3.33352765983949e-05\\
68	3.13798882323954e-05\\
69	2.95768142416458e-05\\
70	2.78988910932559e-05\\
71	2.63533988652314e-05\\
72	2.49069427557121e-05\\
73	2.35619122193325e-05\\
74	2.22948091513225e-05\\
75	2.11157547993773e-05\\
76	2.00104496279949e-05\\
77	1.89856471188971e-05\\
78	1.80239137364968e-05\\
79	1.71251528920825e-05\\
80	1.62747068832878e-05\\
81	1.54763778475813e-05\\
82	1.47234266301421e-05\\
83	1.40200035368952e-05\\
84	1.33586463166044e-05\\
85	1.27380466222525e-05\\
86	1.21500054381333e-05\\
87	1.1594023879612e-05\\
88	1.10673279049933e-05\\
89	1.05714588505001e-05\\
90	1.01043869645543e-05\\
91	9.66449184411818e-06\\
92	9.24772650918006e-06\\
93	8.85169698763425e-06\\
94	8.47514552374247e-06\\
95	8.11810607232529e-06\\
96	7.78080981574905e-06\\
97	7.46196870669996e-06\\
98	7.16004391566297e-06\\
99	6.87221595527182e-06\\
100	6.59788421539382e-06\\
101	6.33608956319761e-06\\
102	6.08790825536539e-06\\
103	5.8523282447401e-06\\
104	5.62926282627797e-06\\
105	5.41623040529459e-06\\
106	5.21287859412802e-06\\
107	5.01781058550922e-06\\
108	4.83212246030624e-06\\
109	4.65505860494167e-06\\
110	4.48723785127448e-06\\
111	4.326793198208e-06\\
112	4.17358456200211e-06\\
113	4.02605605789845e-06\\
114	3.8850760582175e-06\\
115	3.74999091688866e-06\\
116	3.62169288762612e-06\\
117	3.49891059168344e-06\\
118	3.3816918473537e-06\\
119	3.26857398051612e-06\\
120	3.16009926380518e-06\\
121	3.0556677497592e-06\\
122	2.95615779455785e-06\\
123	2.86078184870544e-06\\
124	2.76975371829925e-06\\
125	2.68185580614053e-06\\
126	2.59735582282932e-06\\
127	2.51566935694608e-06\\
128	2.4375165781791e-06\\
129	2.36244514540344e-06\\
130	2.29076966275902e-06\\
131	2.22159419155137e-06\\
132	2.15499947049941e-06\\
133	2.09042159078683e-06\\
134	2.0283642841051e-06\\
135	1.96857710945411e-06\\
136	1.9114108190805e-06\\
137	1.85628895404844e-06\\
138	1.80320684017063e-06\\
139	1.75165066358392e-06\\
140	1.70190064187645e-06\\
141	1.65381499605819e-06\\
142	1.60771367986698e-06\\
143	1.56328742570603e-06\\
144	1.52051828954703e-06\\
145	1.47897424319907e-06\\
146	1.43875385672586e-06\\
147	1.39975018281459e-06\\
148	1.36222126346646e-06\\
149	1.32604120946548e-06\\
150	1.2912207585987e-06\\
151	1.25743876738265e-06\\
152	1.22466653806716e-06\\
153	1.19280119433088e-06\\
154	1.16201327158844e-06\\
155	1.13228969434382e-06\\
156	1.10367118487856e-06\\
157	1.07595607552175e-06\\
158	1.04905311837742e-06\\
159	1.02284740343321e-06\\
160	9.97425451244367e-07\\
161	9.72821389946336e-07\\
162	9.49097626344689e-07\\
163	9.26158893134089e-07\\
164	9.0390378840215e-07\\
165	8.82216558522832e-07\\
166	8.61104663232479e-07\\
167	8.40613139429708e-07\\
168	8.20803362441094e-07\\
169	8.0166252688738e-07\\
170	7.83110974666743e-07\\
171	7.65047608558841e-07\\
172	7.47422742248498e-07\\
173	7.30264811707891e-07\\
174	7.1362178911558e-07\\
175	6.97530162217879e-07\\
176	6.81945269862077e-07\\
177	6.66797906511795e-07\\
178	6.52004728436495e-07\\
179	6.37571005925169e-07\\
180	6.2351769865535e-07\\
181	6.09905137906056e-07\\
182	5.96717918116961e-07\\
183	5.8392664000279e-07\\
184	5.71440882581762e-07\\
185	5.59243629722493e-07\\
186	5.47326613737421e-07\\
187	5.35750754596186e-07\\
188	5.24519860667794e-07\\
189	5.13641055141158e-07\\
190	5.03037397318156e-07\\
191	4.92681596593974e-07\\
192	4.82536890025848e-07\\
193	4.72652380540751e-07\\
194	4.6303671806195e-07\\
195	4.53723923587854e-07\\
196	4.44660718423728e-07\\
197	4.35822808376506e-07\\
198	4.27154184962186e-07\\
199	4.18684045102014e-07\\
200	4.1041641449523e-07\\
201	4.02398285731575e-07\\
202	3.94601093421034e-07\\
203	3.87014560188363e-07\\
204	3.79575746551937e-07\\
205	3.72294453802435e-07\\
206	3.65162330128565e-07\\
207	3.58228006262537e-07\\
208	3.51480757061189e-07\\
209	3.44928480934587e-07\\
210	3.38514285302007e-07\\
211	3.32233841729272e-07\\
212	3.26064419721558e-07\\
213	3.20046682584462e-07\\
214	3.14179407290419e-07\\
215	3.08485875760986e-07\\
216	3.02924145030275e-07\\
217	2.97485384561276e-07\\
218	2.92134409147787e-07\\
219	2.86898821371434e-07\\
220	2.81778316139337e-07\\
221	2.76805168859225e-07\\
222	2.71954592878871e-07\\
223	2.67222230575358e-07\\
224	2.62567052529412e-07\\
225	2.58001968089069e-07\\
226	2.53521689626224e-07\\
227	2.49159917653044e-07\\
228	2.44906053438058e-07\\
229	2.40765448809663e-07\\
230	2.36699653759329e-07\\
231	2.32709900132801e-07\\
232	2.28782477743928e-07\\
233	2.24946370951131e-07\\
234	2.21199224787361e-07\\
235	2.17556139888551e-07\\
236	2.13988021357152e-07\\
237	2.10490577171474e-07\\
238	2.07041885359054e-07\\
239	2.03661972622021e-07\\
240	2.00350731755394e-07\\
241	1.971293844266e-07\\
242	1.93980819242096e-07\\
243	1.90902227514941e-07\\
244	1.87867200481445e-07\\
245	1.84885291530272e-07\\
246	1.81953796429302e-07\\
247	1.79095114598399e-07\\
248	1.76302407271897e-07\\
249	1.73578830098552e-07\\
250	1.70899296290559e-07\\
251	1.68264571460478e-07\\
252	1.65666417405971e-07\\
253	1.63124001046498e-07\\
254	1.60636620297989e-07\\
255	1.58214075865835e-07\\
256	1.55837763986152e-07\\
257	1.53504155884938e-07\\
258	1.51199244802324e-07\\
259	1.48935754245227e-07\\
260	1.46714695090257e-07\\
261	1.4455003355843e-07\\
262	1.42431689103774e-07\\
263	1.40357120578489e-07\\
264	1.3830910979991e-07\\
265	1.36292818136186e-07\\
266	1.34307207707183e-07\\
267	1.32366938864984e-07\\
268	1.30469098264545e-07\\
269	1.28615626608625e-07\\
270	1.26790546249841e-07\\
271	1.24992775795321e-07\\
272	1.23217097153613e-07\\
273	1.21475483993161e-07\\
274	1.19769091439e-07\\
275	1.18104604125593e-07\\
276	1.16471064865442e-07\\
277	1.1486480650495e-07\\
278	1.13276297460521e-07\\
279	1.11712617776787e-07\\
280	1.10175498792713e-07\\
281	1.08674435380837e-07\\
282	1.07204661680725e-07\\
283	1.05764027028772e-07\\
284	1.04341072963964e-07\\
285	1.02937251337728e-07\\
286	1.01552219141235e-07\\
287	1.00195418196365e-07\\
288	9.88669183747541e-08\\
289	9.75683235949581e-08\\
290	9.62898041230832e-08\\
291	9.50287512079896e-08\\
292	9.37812185102661e-08\\
293	9.25542065962452e-08\\
294	9.13499176259496e-08\\
295	9.01734540689237e-08\\
296	8.90192258537762e-08\\
297	8.78837718516834e-08\\
298	8.67600868168691e-08\\
299	8.56510853639038e-08\\
300	8.4558562130784e-08\\
301	8.34891176644918e-08\\
302	8.24416784893806e-08\\
303	8.14149455885039e-08\\
304	8.04012387496627e-08\\
305	7.93994373857085e-08\\
306	7.84089256713981e-08\\
307	7.74356226115313e-08\\
308	7.64812847463992e-08\\
309	7.55478609786795e-08\\
310	7.46299024142499e-08\\
311	7.37241364108725e-08\\
312	7.28269628043574e-08\\
313	7.19417259055058e-08\\
314	7.10708437599271e-08\\
315	7.02184808408445e-08\\
316	6.9382822490826e-08\\
317	6.85612992750825e-08\\
318	6.77484936545059e-08\\
319	6.69444027337517e-08\\
320	6.61500843754805e-08\\
321	6.53701162773379e-08\\
322	6.46056226536852e-08\\
323	6.38565934163931e-08\\
324	6.31183033903348e-08\\
325	6.23881801185546e-08\\
326	6.16648270242932e-08\\
327	6.09513919954192e-08\\
328	6.02502350632875e-08\\
329	5.95637232168275e-08\\
330	5.88897654844155e-08\\
331	5.82254959760115e-08\\
332	5.75673557170893e-08\\
333	5.69159391533668e-08\\
334	5.6272982935867e-08\\
335	5.56418548548317e-08\\
336	5.50231771132227e-08\\
337	5.44159056693096e-08\\
338	5.38162254659242e-08\\
339	5.32222429035373e-08\\
340	5.2633805509374e-08\\
341	5.20536022426673e-08\\
342	5.14837158962157e-08\\
343	5.09253407390914e-08\\
344	5.03764563566624e-08\\
345	4.98343504387523e-08\\
346	4.92967672375399e-08\\
347	4.87644957924902e-08\\
348	4.82394348475041e-08\\
349	4.77240197685866e-08\\
350	4.72185903896772e-08\\
351	4.67216826630263e-08\\
352	4.62303191187043e-08\\
353	4.57430746353012e-08\\
354	4.52603963669944e-08\\
355	4.47844700352483e-08\\
356	4.4317098660164e-08\\
357	4.38588103059167e-08\\
358	4.340785083707e-08\\
359	4.29617963932193e-08\\
360	4.25191930188141e-08\\
361	4.20807776147141e-08\\
362	4.16483759221832e-08\\
363	4.12237800032727e-08\\
364	4.08072494245536e-08\\
365	4.03972630262291e-08\\
366	3.999150220907e-08\\
367	3.95887868903891e-08\\
368	3.918976924022e-08\\
369	3.87961829740175e-08\\
370	3.84096239799255e-08\\
371	3.80303215782283e-08\\
372	3.76568739263406e-08\\
373	3.72871389751881e-08\\
374	3.69200852152391e-08\\
375	3.65562629281195e-08\\
376	3.61973289302384e-08\\
377	3.58446698411892e-08\\
378	3.54986143163439e-08\\
379	3.51577879608679e-08\\
380	3.48203521110261e-08\\
381	3.44851971767535e-08\\
382	3.41529448903997e-08\\
383	3.38249464115308e-08\\
384	3.35026653733577e-08\\
385	3.31862778047129e-08\\
386	3.28747679240638e-08\\
387	3.25662384235878e-08\\
388	3.22598219196918e-08\\
389	3.19558067909497e-08\\
390	3.16556373116064e-08\\
391	3.1360457828327e-08\\
392	3.10707599089895e-08\\
393	3.07854390418242e-08\\
394	3.05030144023058e-08\\
395	3.02223338537934e-08\\
396	2.99438409400625e-08\\
397	2.96685505511227e-08\\
398	2.93978358679126e-08\\
399	2.91319670551749e-08\\
400	2.88703377960225e-08\\
401	2.8611289036948e-08\\
402	2.83539637148039e-08\\
403	2.80983367493858e-08\\
404	2.78455994223836e-08\\
405	2.75967379950726e-08\\
406	2.73524467995437e-08\\
407	2.71119785377661e-08\\
408	2.6874171643655e-08\\
409	2.66377847337264e-08\\
410	2.64029825682504e-08\\
411	2.61704485308091e-08\\
412	2.59414583178393e-08\\
413	2.571643898025e-08\\
414	2.54952172077201e-08\\
415	2.5276426274001e-08\\
416	2.50591665273559e-08\\
417	2.48430606940777e-08\\
418	2.46289874744515e-08\\
419	2.44177824649038e-08\\
420	2.42103235072961e-08\\
421	2.40062768198565e-08\\
422	2.38048466634994e-08\\
423	2.36047427083666e-08\\
424	2.34057900248465e-08\\
425	2.32083020335809e-08\\
426	2.30133965767614e-08\\
427	2.28216366905886e-08\\
428	2.26332810675304e-08\\
429	2.24473651833686e-08\\
430	2.22630122964838e-08\\
431	2.20794900039449e-08\\
432	2.18972837359278e-08\\
433	2.17170169878783e-08\\
434	2.15396638886277e-08\\
435	2.1365297777742e-08\\
436	2.11935842013473e-08\\
437	2.10233375808801e-08\\
438	2.0854059838358e-08\\
439	2.06856312293404e-08\\
440	2.05189015916741e-08\\
441	2.03544681681507e-08\\
442	2.01929474432863e-08\\
443	2.00338834036683e-08\\
444	1.98766141179525e-08\\
445	1.97201522183867e-08\\
446	1.95645197787612e-08\\
447	1.94100118097159e-08\\
448	1.92575470615127e-08\\
449	1.91074966201784e-08\\
450	1.8960034159049e-08\\
451	1.88143277223765e-08\\
452	1.86697278766177e-08\\
453	1.85256607340623e-08\\
454	1.83825655936556e-08\\
455	1.82409051434348e-08\\
456	1.81014787492812e-08\\
457	1.79643280056846e-08\\
458	1.78292419346733e-08\\
459	1.76952608440318e-08\\
460	1.7561986192433e-08\\
461	1.74292524580301e-08\\
462	1.72977231767057e-08\\
463	1.7167855564434e-08\\
464	1.70402192769424e-08\\
465	1.69145277872314e-08\\
466	1.67903347813108e-08\\
467	1.66667824089566e-08\\
468	1.65437968727878e-08\\
469	1.64214962321806e-08\\
470	1.63006118828815e-08\\
471	1.61814854329943e-08\\
472	1.60644255884566e-08\\
473	1.59488886473503e-08\\
474	1.58343693164837e-08\\
475	1.57202363452521e-08\\
476	1.56066885379191e-08\\
477	1.54939978341433e-08\\
478	1.5382855529421e-08\\
479	1.52734429777449e-08\\
480	1.51658254602009e-08\\
481	1.50593175505725e-08\\
482	1.49534884010753e-08\\
483	1.48479518331272e-08\\
484	1.47430886341536e-08\\
485	1.46392212367696e-08\\
486	1.45369491012225e-08\\
487	1.44362837913252e-08\\
488	1.43371058217998e-08\\
489	1.42386959884336e-08\\
490	1.41407600199938e-08\\
491	1.40431160978703e-08\\
492	1.39462376360881e-08\\
493	1.38504399794672e-08\\
494	1.37562079078377e-08\\
495	1.3663407852631e-08\\
496	1.3571803978002e-08\\
497	1.34807197399844e-08\\
498	1.33900042791394e-08\\
499	1.32996169249674e-08\\
500	1.32100643523629e-08\\
501	1.31216221198392e-08\\
502	1.30346663723869e-08\\
503	1.29489546595107e-08\\
504	1.28641966441972e-08\\
505	1.27797938429817e-08\\
506	1.26957109993705e-08\\
507	1.26119949448759e-08\\
508	1.2529145789417e-08\\
509	1.24473929060803e-08\\
510	1.23670200876521e-08\\
511	1.22877188909994e-08\\
512	1.22091874260789e-08\\
513	1.21309097995293e-08\\
514	1.20529304887447e-08\\
515	1.19753403389083e-08\\
516	1.18986158685279e-08\\
517	1.18229425993733e-08\\
518	1.17485382880363e-08\\
519	1.16750646708977e-08\\
520	1.16022316597788e-08\\
521	1.15295926896398e-08\\
522	1.14572352806601e-08\\
523	1.13852675199941e-08\\
524	1.13141332130927e-08\\
525	1.1243986373865e-08\\
526	1.11750042060784e-08\\
527	1.11068490731639e-08\\
528	1.10392524975406e-08\\
529	1.09718163388739e-08\\
530	1.09046422024382e-08\\
531	1.0837835081888e-08\\
532	1.07718062011525e-08\\
533	1.07066899516332e-08\\
534	1.06426481481059e-08\\
535	1.05793658803167e-08\\
536	1.05165986935092e-08\\
537	1.04539760276131e-08\\
538	1.03915885547449e-08\\
539	1.03295256621437e-08\\
540	1.02681653991282e-08\\
541	1.02076352456626e-08\\
542	1.01481012316263e-08\\
543	1.00892903487939e-08\\
544	1.00309847331568e-08\\
545	9.97282409046694e-09\\
546	9.91486889998579e-09\\
547	9.85718024692593e-09\\
548	9.80010321123638e-09\\
549	9.74376478242883e-09\\
550	9.68835082515835e-09\\
551	9.63364608773646e-09\\
552	9.57946647661961e-09\\
553	9.52545256004183e-09\\
554	9.47161743631578e-09\\
555	9.41798189631778e-09\\
556	9.3648530812843e-09\\
557	9.31235954919066e-09\\
558	9.26071356734211e-09\\
559	9.20977234161314e-09\\
560	9.15940092907013e-09\\
561	9.109243659516e-09\\
562	9.05925571583446e-09\\
563	9.00940127824535e-09\\
564	8.95994066344525e-09\\
565	8.91099704650041e-09\\
566	8.86280784671514e-09\\
567	8.81530968850689e-09\\
568	8.76843864296574e-09\\
569	8.72186354273396e-09\\
570	8.6754813497613e-09\\
571	8.6291838217999e-09\\
572	8.58316759180485e-09\\
573	8.53753698364221e-09\\
574	8.49254135341219e-09\\
575	8.44819159972017e-09\\
576	8.40451393019109e-09\\
577	8.36124097984477e-09\\
578	8.3182357314664e-09\\
579	8.27530853230189e-09\\
580	8.23257002104886e-09\\
581	8.19008214487055e-09\\
582	8.14808464446844e-09\\
583	8.1066401480146e-09\\
584	8.0658649693908e-09\\
585	8.02559763620769e-09\\
586	7.98572401924686e-09\\
587	7.9459938377672e-09\\
588	7.90641202131365e-09\\
589	7.86696609798738e-09\\
590	7.82785494054178e-09\\
591	7.78915791472128e-09\\
592	7.75105177048303e-09\\
593	7.71349346860358e-09\\
594	7.67646129095766e-09\\
595	7.63971828135876e-09\\
596	7.60318102033047e-09\\
597	7.56673129935543e-09\\
598	7.5304861283541e-09\\
599	7.49449707243477e-09\\
600	7.45895348346394e-09\\
601	7.42389238471633e-09\\
602	7.38940564188048e-09\\
603	7.35536542991168e-09\\
604	7.32169494012864e-09\\
605	7.28819412127749e-09\\
606	7.25486126743518e-09\\
607	7.22166581889154e-09\\
608	7.188754422987e-09\\
609	7.15617989520293e-09\\
610	7.12409420960757e-09\\
611	7.09249249865883e-09\\
612	7.06141194619446e-09\\
613	7.03070339668234e-09\\
614	7.00030247373536e-09\\
615	6.97006813773967e-09\\
616	6.94004280136571e-09\\
617	6.91022547780834e-09\\
618	6.88075338163262e-09\\
619	6.85166703490167e-09\\
620	6.82309126351362e-09\\
621	6.79501428590896e-09\\
622	6.76747601844371e-09\\
623	6.74037084302021e-09\\
624	6.71367573280236e-09\\
625	6.68727178005468e-09\\
626	6.66114583604342e-09\\
627	6.63526475925443e-09\\
628	6.61002759759282e-09\\
629	6.58650551234422e-09\\
630	6.56510029427929e-09\\
631	6.54235282896447e-09\\
632	6.51697006605059e-09\\
633	6.49024795602836e-09\\
634	6.47069197227885e-09\\
635	6.46215981967158e-09\\
636	6.4345601775403e-09\\
637	6.40629689593619e-09\\
638	6.37770035020072e-09\\
639	6.34890063300045e-09\\
640	6.32004458918199e-09\\
641	6.3092840289009e-09\\
642	6.29127637117728e-09\\
643	6.26270568224303e-09\\
644	6.23436516401349e-09\\
645	6.2061842668613e-09\\
646	6.17802438822729e-09\\
647	6.14977135526281e-09\\
648	6.12140496437313e-09\\
649	6.10761400938165e-09\\
650	6.09299694671349e-09\\
651	6.06466936039118e-09\\
652	6.0365504509868e-09\\
653	6.00873051903068e-09\\
654	5.9812137129986e-09\\
655	5.95389589194374e-09\\
656	5.92661772321473e-09\\
657	5.89926915505024e-09\\
658	5.87550272078862e-09\\
659	5.87184853926473e-09\\
660	5.84444387807019e-09\\
661	5.81718339163991e-09\\
662	5.7901892456408e-09\\
663	5.76352944366637e-09\\
664	5.73717034858262e-09\\
665	5.71097451543793e-09\\
666	5.68478288235674e-09\\
667	5.65851855790078e-09\\
668	5.6322168316161e-09\\
669	5.61768215526031e-09\\
670	5.60598647588107e-09\\
671	5.57995838760933e-09\\
672	5.55423884671894e-09\\
673	5.52885868500113e-09\\
674	5.50373581275576e-09\\
675	5.47870837822342e-09\\
676	5.45364193789385e-09\\
677	5.42850998868092e-09\\
678	5.40338646121878e-09\\
679	5.37839491196688e-09\\
680	5.35365944163433e-09\\
681	5.33743877368453e-09\\
682	5.32925696017184e-09\\
683	5.30516689059984e-09\\
684	5.28125944634886e-09\\
685	5.25737199064818e-09\\
686	5.23341868218987e-09\\
687	5.20942750222433e-09\\
688	5.18550229974247e-09\\
689	5.16177050028415e-09\\
690	5.13833710091206e-09\\
691	5.11523427203585e-09\\
692	5.09238252497541e-09\\
693	5.06962229377873e-09\\
694	5.04682186137099e-09\\
695	5.03746986596653e-09\\
696	5.02395728339116e-09\\
697	5.00110240323797e-09\\
698	4.97837724577192e-09\\
699	4.95590060171337e-09\\
700	4.93374363077374e-09\\
701	4.91188041590591e-09\\
702	4.89017953987173e-09\\
703	4.86848439867324e-09\\
704	4.84672008090836e-09\\
705	4.82492184106668e-09\\
706	4.80319349350305e-09\\
707	4.78165702125897e-09\\
708	4.76040882518452e-09\\
709	4.73947050538267e-09\\
710	4.72021618187424e-09\\
711	4.71875389508494e-09\\
712	4.69810162739278e-09\\
713	4.67739604394314e-09\\
714	4.65662810692007e-09\\
715	4.63587739939282e-09\\
716	4.61526087807627e-09\\
717	4.59488878897753e-09\\
718	4.57481969725113e-09\\
719	4.55501359621346e-09\\
720	4.53533335995552e-09\\
721	4.51563392624961e-09\\
722	4.49586038721893e-09\\
723	4.4760607898075e-09\\
724	4.45634014245494e-09\\
725	4.43681369617063e-09\\
726	4.41756563241112e-09\\
727	4.39860084903267e-09\\
728	4.38802777091277e-09\\
729	4.37981779834312e-09\\
730	4.3610630026491e-09\\
731	4.34223987641734e-09\\
732	4.32335898820335e-09\\
733	4.30450694696618e-09\\
734	4.28579713233298e-09\\
735	4.26732934760067e-09\\
736	4.24914540173672e-09\\
737	4.23118642023389e-09\\
738	4.21330942952607e-09\\
739	4.19538624690286e-09\\
740	4.1773867069295e-09\\
741	4.15937291239401e-09\\
742	4.14145072233632e-09\\
743	4.12372758916989e-09\\
744	4.10627256529035e-09\\
745	4.08906914605401e-09\\
746	4.07200034931264e-09\\
747	4.05492091287097e-09\\
748	4.04326244219714e-09\\
749	4.03776077940252e-09\\
750	4.02055208192624e-09\\
751	4.00338812163814e-09\\
752	3.98637791884239e-09\\
753	3.96960874992127e-09\\
754	3.95310095163438e-09\\
755	3.93677259491307e-09\\
756	3.92047739141815e-09\\
757	3.90411095884572e-09\\
758	3.88767156805137e-09\\
759	3.87123539545085e-09\\
760	3.85490835765319e-09\\
761	3.83878844015491e-09\\
762	3.82292499356161e-09\\
763	3.80727374570628e-09\\
764	3.79170192359705e-09\\
765	3.77608045523405e-09\\
766	3.76037247753959e-09\\
767	3.74463280196933e-09\\
768	3.72896051234301e-09\\
769	3.71344635614983e-09\\
770	3.69819747144998e-09\\
771	3.68835084281143e-09\\
772	3.6831673873489e-09\\
773	3.66825916816924e-09\\
774	3.65333097514932e-09\\
775	3.63831308438478e-09\\
776	3.62323594589399e-09\\
777	3.60818891120485e-09\\
778	3.5932762564982e-09\\
779	3.57858224844532e-09\\
780	3.5641267944492e-09\\
781	3.54982861404192e-09\\
782	3.53554403045236e-09\\
783	3.52117521622379e-09\\
784	3.50672676333679e-09\\
785	3.4922760974435e-09\\
786	3.47792608333398e-09\\
787	3.46376984939333e-09\\
788	3.44985005834078e-09\\
789	3.43611354796835e-09\\
790	3.42242367191762e-09\\
791	3.4086613646643e-09\\
792	3.39480600791828e-09\\
793	3.38092175176985e-09\\
794	3.36710836030949e-09\\
795	3.35346407908698e-09\\
796	3.34004822155234e-09\\
797	3.32683257317623e-09\\
798	3.32584859395888e-09\\
799	3.313692421779e-09\\
800	3.300494736482e-09\\
801	3.2871962662341e-09\\
802	3.2738478559046e-09\\
803	3.26054539750847e-09\\
804	3.2473899229878e-09\\
805	3.23445235422542e-09\\
806	3.22172443856211e-09\\
807	3.20909484048612e-09\\
808	3.19642200462375e-09\\
809	3.18364434269613e-09\\
810	3.17080088209673e-09\\
811	3.15798378985162e-09\\
812	3.14529557793016e-09\\
813	3.13281489062661e-09\\
814	3.12054805121502e-09\\
815	3.10839501084309e-09\\
816	3.09620977374195e-09\\
817	3.08391744107611e-09\\
818	3.07154817708278e-09\\
819	3.05919145553901e-09\\
820	3.04695071688303e-09\\
821	3.0349094264701e-09\\
822	3.02308296673053e-09\\
823	3.01137835078502e-09\\
824	2.99964682508587e-09\\
825	2.98780574375295e-09\\
826	2.97588075223254e-09\\
827	2.96396045903282e-09\\
828	2.9584799054192e-09\\
829	2.95214914231927e-09\\
830	2.94053277304282e-09\\
831	2.92913030755459e-09\\
832	2.91785006017614e-09\\
833	2.90654092024073e-09\\
834	2.89511854113801e-09\\
835	2.88360941565899e-09\\
836	2.87210353634263e-09\\
837	2.86070618928052e-09\\
838	2.84950363585976e-09\\
839	2.83851259144603e-09\\
840	2.82763586734347e-09\\
841	2.81672000388965e-09\\
842	2.80568583107712e-09\\
843	2.79456657366204e-09\\
844	2.78345592222895e-09\\
845	2.77246018434864e-09\\
846	2.76166333707929e-09\\
847	2.75107301130923e-09\\
848	2.74057937844243e-09\\
849	2.73002826972769e-09\\
850	2.71935386695406e-09\\
851	2.70860182337348e-09\\
852	2.69787119711597e-09\\
853	2.68726882586816e-09\\
854	2.6768727005073e-09\\
855	2.66667283874119e-09\\
856	2.65654037893182e-09\\
857	2.64632564208674e-09\\
858	2.63598558296042e-09\\
859	2.62558264240597e-09\\
860	2.61522183745702e-09\\
861	2.60500909163809e-09\\
862	2.59501039239998e-09\\
863	2.5891811877753e-09\\
864	2.58518743921777e-09\\
865	2.57538958738218e-09\\
866	2.56548263703789e-09\\
867	2.55545562279081e-09\\
868	2.54538952509289e-09\\
869	2.53539448973571e-09\\
870	2.5255722776073e-09\\
871	2.51596705731009e-09\\
872	2.50650041923798e-09\\
873	2.49700513560258e-09\\
874	2.48737925748082e-09\\
875	2.47765018435011e-09\\
876	2.46791592957799e-09\\
877	2.45828931612674e-09\\
878	2.44886141117096e-09\\
879	2.43963915193549e-09\\
880	2.43049580709237e-09\\
881	2.42126656216437e-09\\
882	2.41190016966704e-09\\
883	2.40246225714427e-09\\
884	2.39306332633265e-09\\
885	2.38381440263693e-09\\
886	2.37478291182323e-09\\
887	2.36591829194159e-09\\
888	2.35705166218138e-09\\
889	2.34805365941478e-09\\
890	2.33893487900295e-09\\
891	2.32979158076033e-09\\
892	2.3207406492072e-09\\
893	2.31188226415872e-09\\
894	2.30323656378079e-09\\
895	2.29467980796236e-09\\
896	2.28603512147908e-09\\
897	2.27724390133277e-09\\
898	2.26837410506977e-09\\
899	2.25954068095443e-09\\
900	2.25085903153671e-09\\
901	2.24239817641015e-09\\
902	2.23409777335193e-09\\
903	2.22577497072104e-09\\
904	2.2211478200746e-09\\
905	2.21730523125779e-09\\
906	2.20871510893507e-09\\
907	2.20011130706455e-09\\
908	2.19161575162458e-09\\
909	2.18332652454868e-09\\
910	2.175240868745e-09\\
911	2.16720062230343e-09\\
912	2.15903332933814e-09\\
913	2.15071728597188e-09\\
914	2.14234548879526e-09\\
915	2.13404241813908e-09\\
916	2.12592278373995e-09\\
917	2.11802650620242e-09\\
918	2.11023142117397e-09\\
919	2.10233842005892e-09\\
920	2.09427989498618e-09\\
921	2.08613241634943e-09\\
922	2.07802164391348e-09\\
923	2.07007248168788e-09\\
924	2.06235313606169e-09\\
925	2.05477331572658e-09\\
926	2.04712165008676e-09\\
927	2.03929461096315e-09\\
928	2.03135471586705e-09\\
929	2.02342947801494e-09\\
930	2.01565145178813e-09\\
931	2.00810597673466e-09\\
932	2.00072106643599e-09\\
933	1.99327665023147e-09\\
934	1.98564783473282e-09\\
935	1.97789196159375e-09\\
936	1.97014125681717e-09\\
937	1.96253540604934e-09\\
938	1.95516758878279e-09\\
939	1.9479642859634e-09\\
940	1.94069194093616e-09\\
941	1.93322331639421e-09\\
942	1.92562517512605e-09\\
943	1.91803882352085e-09\\
944	1.91061112215058e-09\\
945	1.90343190274848e-09\\
946	1.89639883712963e-09\\
947	1.8892599136325e-09\\
948	1.88191254222658e-09\\
949	1.87444987650853e-09\\
950	1.86702626905871e-09\\
951	1.85979416138524e-09\\
952	1.85788953238037e-09\\
953	1.85281911368063e-09\\
954	1.84593504347151e-09\\
955	1.83888413380159e-09\\
956	1.8316253712008e-09\\
957	1.82429003199272e-09\\
958	1.81704666401687e-09\\
959	1.81004438683535e-09\\
960	1.80328083297048e-09\\
961	1.796499152533e-09\\
962	1.78949001543809e-09\\
963	1.78230696097123e-09\\
964	1.77511781576959e-09\\
965	1.76810173242318e-09\\
966	1.76137075844316e-09\\
967	1.75478253225929e-09\\
968	1.74802874558797e-09\\
969	1.74103711535334e-09\\
970	1.73395373525247e-09\\
971	1.72697035696881e-09\\
972	1.72025679548705e-09\\
973	1.71379470465441e-09\\
974	1.70727166958423e-09\\
975	1.7004870344309e-09\\
976	1.69354168918764e-09\\
977	1.68663395831023e-09\\
978	1.67996195914299e-09\\
979	1.67359658508081e-09\\
980	1.66727165539809e-09\\
981	1.66069098069273e-09\\
982	1.65389923692619e-09\\
983	1.6470978837905e-09\\
984	1.64050175374839e-09\\
985	1.63423689104356e-09\\
986	1.62808697756043e-09\\
987	1.62169633729499e-09\\
988	1.61506116801121e-09\\
989	1.60838583530824e-09\\
990	1.60189768016624e-09\\
991	1.59575447343158e-09\\
992	1.5897692148753e-09\\
993	1.58354986863302e-09\\
994	1.57706614821412e-09\\
995	1.57052905359277e-09\\
996	1.56417685910201e-09\\
997	1.55818315724643e-09\\
998	1.55235993117814e-09\\
999	1.5462911364948e-09\\
1000	1.53994973589665e-09\\
1001	1.53355941476723e-09\\
1002	1.52736848859349e-09\\
1003	1.52155242257564e-09\\
1004	1.51588748830127e-09\\
1005	1.50994647778073e-09\\
1006	1.50373779350667e-09\\
1007	1.50329985628628e-09\\
1008	1.4975023982649e-09\\
1009	1.4914980155417e-09\\
1010	1.48588425437613e-09\\
1011	1.48036500414275e-09\\
1012	1.47452737997922e-09\\
1013	1.46844407832306e-09\\
1014	1.46237414976789e-09\\
1015	1.45658357812907e-09\\
1016	1.45118713002697e-09\\
1017	1.44578494426034e-09\\
1018	1.44002644039884e-09\\
1019	1.43406563814164e-09\\
1020	1.42817660794883e-09\\
1021	1.42263039937119e-09\\
1022	1.41744621258643e-09\\
1023	1.41211428659674e-09\\
1024	1.4064145289868e-09\\
1025	1.40057919832716e-09\\
1026	1.39489503565988e-09\\
1027	1.38962381906759e-09\\
1028	1.38460834253795e-09\\
1029	1.37928634014683e-09\\
1030	1.37363255467719e-09\\
1031	1.36793413310508e-09\\
1032	1.36249433723613e-09\\
1033	1.35751167174814e-09\\
1034	1.35256309359157e-09\\
1035	1.34719254218697e-09\\
1036	1.34158314647236e-09\\
1037	1.33605102617921e-09\\
1038	1.33091910169569e-09\\
1039	1.32616462970688e-09\\
1040	1.32112539047576e-09\\
1041	1.31567189497532e-09\\
1042	1.31012702858213e-09\\
1043	1.30484003755492e-09\\
1044	1.30007117714138e-09\\
1045	1.29531182575019e-09\\
1046	1.29004146078271e-09\\
1047	1.28452087828195e-09\\
1048	1.27913000249352e-09\\
1049	1.27425898344901e-09\\
1050	1.26966357637874e-09\\
1051	1.26454873731177e-09\\
1052	1.25904720857484e-09\\
1053	1.25358972780365e-09\\
1054	1.24864251774481e-09\\
1055	1.24412712784007e-09\\
1056	1.23908181061184e-09\\
1057	1.23355989365836e-09\\
1058	1.22806142248735e-09\\
1059	1.22313402802822e-09\\
1060	1.21866994371666e-09\\
1061	1.21356257112706e-09\\
1062	1.20797020765072e-09\\
1063	1.20248374773211e-09\\
1064	1.19773963271411e-09\\
1065	1.19325566906808e-09\\
1066	1.1879370620401e-09\\
1067	1.18226573287207e-09\\
1068	1.17694859130816e-09\\
1069	1.17255450602043e-09\\
1070	1.16778689991164e-09\\
1071	1.16218874604627e-09\\
1072	1.16174881443072e-09\\
1073	1.15655512381232e-09\\
1074	1.15177270043348e-09\\
1075	1.14749461256953e-09\\
1076	1.1421528710179e-09\\
1077	1.13641189110257e-09\\
1078	1.13123621807125e-09\\
1079	1.12716011336585e-09\\
1080	1.12217276625911e-09\\
1081	1.11642948970424e-09\\
1082	1.11101269686026e-09\\
1083	1.10696108796671e-09\\
1084	1.10231298437894e-09\\
1085	1.0966037693086e-09\\
1086	1.09107196808829e-09\\
1087	1.08705258455058e-09\\
1088	1.08265131224704e-09\\
1089	1.07696297315857e-09\\
1090	1.07140211765041e-09\\
1091	1.06754312521286e-09\\
1092	1.06324964383928e-09\\
1093	1.05754102806082e-09\\
1094	1.05202837015809e-09\\
1095	1.04852040293406e-09\\
1096	1.04413122209235e-09\\
1097	1.03836679653745e-09\\
1098	1.03302435453132e-09\\
1099	1.03003790966289e-09\\
1100	1.02528308506034e-09\\
1101	1.01946585303923e-09\\
1102	1.01458233638253e-09\\
1103	1.01200136430838e-09\\
1104	1.00668643467198e-09\\
1105	1.0008684119387e-09\\
1106	9.97054728383272e-10\\
1107	9.94076617429562e-10\\
1108	9.88346275961533e-10\\
1109	9.826492018678e-10\\
1110	9.80463616153138e-10\\
1111	9.76131651233967e-10\\
1112	9.70285781030968e-10\\
1113	9.65384595009464e-10\\
1114	9.63957206037687e-10\\
1115	9.58339090598565e-10\\
1116	9.52529567099832e-10\\
1117	9.4979077874185e-10\\
1118	9.46574602168864e-10\\
1119	9.40800852609121e-10\\
1120	9.35202124703908e-10\\
1121	9.34651945054489e-10\\
1122	9.29247985124345e-10\\
1123	9.23532201286755e-10\\
1124	9.20143354671349e-10\\
1125	9.17815080118749e-10\\
1126	9.12145082791419e-10\\
1127	9.06495323433442e-10\\
1128	9.05664288183068e-10\\
1129	9.00839740266165e-10\\
1130	8.95205380545038e-10\\
1131	8.91431356864627e-10\\
1132	8.89570785209724e-10\\
1133	8.83933976493174e-10\\
1134	8.78283312770828e-10\\
1135	8.77405027572313e-10\\
1136	8.72607648882599e-10\\
1137	8.66898492165003e-10\\
1138	8.63594869596767e-10\\
1139	8.61142802544089e-10\\
1140	8.55328251855124e-10\\
1141	8.49995241015572e-10\\
1142	8.49441402759182e-10\\
1143	8.4346917753292e-10\\
1144	8.38199086535864e-10\\
1145	8.37398205751283e-10\\
1146	8.36603962899668e-10\\
1147	8.31215748878511e-10\\
1148	8.2490964374078e-10\\
1149	8.23418656756095e-10\\
1150	8.18468349112733e-10\\
1151	8.11880940954764e-10\\
1152	8.10436653020513e-10\\
1153	8.05136883729382e-10\\
1154	7.9822574185795e-10\\
1155	7.97655442484624e-10\\
1156	7.9113670014422e-10\\
1157	7.85072533003105e-10\\
1158	7.8385801146846e-10\\
1159	7.76376193241335e-10\\
1160	7.72685464845622e-10\\
1161	7.68675011328363e-10\\
1162	7.60734042328511e-10\\
1163	7.60491806601911e-10\\
1164	7.52526631979796e-10\\
1165	7.48489156324017e-10\\
1166	7.44016673612011e-10\\
1167	7.36675135367943e-10\\
1168	7.35153407498129e-10\\
1169	7.25862159919772e-10\\
1170	7.25047397353586e-10\\
1171	7.16026175259605e-10\\
1172	7.13603614334542e-10\\
1173	7.05445931804533e-10\\
1174	7.02341486777893e-10\\
1175	6.93726503432576e-10\\
1176	6.91258737468441e-10\\
1177	6.80353110394984e-10\\
1178	6.79817927243851e-10\\
1179	6.69622373924827e-10\\
1180	6.59064313584751e-10\\
1181	6.48676741047392e-10\\
1182	6.38457481214361e-10\\
1183	6.28404385312945e-10\\
1184	6.18515318145077e-10\\
1185	6.08788166465432e-10\\
1186	5.99220834523904e-10\\
1187	5.89811242637107e-10\\
1188	5.80557330125127e-10\\
1189	5.7145705342133e-10\\
1190	5.62508384007625e-10\\
1191	5.53709312157891e-10\\
1192	5.45057842835668e-10\\
1193	5.38281143509769e-10\\
1194	5.36551994480266e-10\\
1195	5.28189809612025e-10\\
1196	5.19969334864045e-10\\
1197	5.11888640497838e-10\\
1198	5.03945809125416e-10\\
1199	4.96138934432693e-10\\
1200	4.88466128940841e-10\\
1201	4.80925519164799e-10\\
1202	4.73515243083184e-10\\
1203	4.66233452514765e-10\\
1204	4.59078316086038e-10\\
1205	4.52048011611698e-10\\
1206	4.45140730779213e-10\\
1207	4.38354677567106e-10\\
1208	4.31688068093934e-10\\
1209	4.25139128260132e-10\\
1210	4.18706096075008e-10\\
1211	4.1238722101238e-10\\
1212	4.06180758966716e-10\\
1213	4.00084975444302e-10\\
1214	3.94098144740194e-10\\
1215	3.88218546526709e-10\\
1216	3.82444465875808e-10\\
1217	3.76774192745495e-10\\
1218	3.71206018364586e-10\\
1219	3.65738232393501e-10\\
1220	3.60369127632201e-10\\
1221	3.55096987413752e-10\\
1222	3.49920089738603e-10\\
1223	3.44836702060248e-10\\
1224	3.39845075719071e-10\\
1225	3.34943442621293e-10\\
1226	3.30130009805547e-10\\
1227	3.25402952871418e-10\\
1228	3.20760407748866e-10\\
1229	3.16200462455784e-10\\
1230	3.11721147337622e-10\\
1231	3.07320420608713e-10\\
1232	3.0299615976384e-10\\
1233	2.98746138604608e-10\\
1234	2.94568012048224e-10\\
1235	2.90459296572948e-10\\
1236	2.86417340067633e-10\\
1237	2.82439298408608e-10\\
1238	2.78522103297923e-10\\
1239	2.74662429067474e-10\\
1240	2.70856658587016e-10\\
1241	2.67550828199972e-10\\
1242	2.67100847944261e-10\\
1243	2.63390699531753e-10\\
1244	2.59721533995637e-10\\
1245	2.56088282711398e-10\\
1246	2.52485496496198e-10\\
1247	2.48907379837181e-10\\
1248	2.45347856803538e-10\\
1249	2.41800671967423e-10\\
1250	2.38259525144715e-10\\
1251	2.34718231243673e-10\\
1252	2.31170889881291e-10\\
1253	2.27612059473093e-10\\
1254	2.24036903941333e-10\\
1255	2.20441304542984e-10\\
1256	2.16821934782222e-10\\
1257	2.13176282365063e-10\\
1258	2.09502640484542e-10\\
1259	2.05800061145763e-10\\
1260	2.02068288700731e-10\\
1261	1.98307683707899e-10\\
1262	1.94519137955904e-10\\
1263	1.90703996207169e-10\\
1264	1.86863979005009e-10\\
1265	1.83001114482125e-10\\
1266	1.79117680013796e-10\\
1267	1.75216147550009e-10\\
1268	1.71299142299916e-10\\
1269	1.67369403597949e-10\\
1270	1.6342975648451e-10\\
1271	1.59483082823329e-10\\
1272	1.55532303361949e-10\\
1273	1.51580358103083e-10\\
1274	1.47630193928471e-10\\
1275	1.43684751083934e-10\\
1276	1.39746954648959e-10\\
1277	1.35819706878282e-10\\
1278	1.31905880530719e-10\\
1279	1.28008313890401e-10\\
1280	1.24129806105431e-10\\
1281	1.20273113718708e-10\\
1282	1.16440949271217e-10\\
1283	1.12635976754396e-10\\
1284	1.08860811619115e-10\\
1285	1.05118018278211e-10\\
1286	1.01410109751175e-10\\
1287	9.77395450433236e-11\\
1288	9.41087309002212e-11\\
1289	9.05200184149927e-11\\
1290	8.69757049063545e-11\\
1291	8.34780325347715e-11\\
1292	8.00291871027489e-11\\
1293	7.66312997332375e-11\\
1294	7.32864453793031e-11\\
1295	6.99966436988059e-11\\
1296	6.67638573450134e-11\\
1297	6.35899941407313e-11\\
1298	6.04769057147247e-11\\
1299	5.74263878617594e-11\\
1300	5.44401807658222e-11\\
1301	5.15199692793794e-11\\
1302	4.86673825055542e-11\\
1303	4.58839945816745e-11\\
1304	4.31713244394651e-11\\
1305	4.0530836540321e-11\\
1306	3.79639402287406e-11\\
1307	3.54719910367227e-11\\
1308	3.30562901895355e-11\\
1309	3.07180852470516e-11\\
1310	2.84585702396032e-11\\
1311	2.62788861841573e-11\\
1312	2.41801214502755e-11\\
1313	2.21633119003043e-11\\
1314	2.02294417513138e-11\\
1315	1.83794435411414e-11\\
1316	1.66141989166999e-11\\
1317	1.49345390041229e-11\\
1318	1.33412450292086e-11\\
1319	1.18350488038707e-11\\
1320	1.04166333775009e-11\\
1321	9.08663367272502e-12\\
1322	7.84563713643502e-12\\
1323	6.69418453618184e-12\\
1324	5.63277068300557e-12\\
1325	4.66184530618332e-12\\
1326	3.78181399725162e-12\\
1327	2.99303918155136e-12\\
1328	2.29584131982787e-12\\
1329	1.69050030105804e-12\\
1330	1.17725724541823e-12\\
1331	7.56317175620545e-13\\
1332	4.27853883918298e-13\\
1333	1.92021831245545e-13\\
1334	4.90028275571477e-14\\
};
\addlegendentry{$E=2$}
\end{axis}
\end{tikzpicture}%

%% file: wave-svd-M=1000-E=16.tikz
%
%
%
\begin{tikzpicture}

\begin{axis}[%
width=0.35\textwidth,
scale only axis,
xmin=0,
xmax=200,
xlabel={Index $i$},
ymode=log,
ymin=1e-9,
ymax=1e2,
ylabel={Singular value $\sigma_i$}
]
\addplot [color=blue]
  table[row sep=crcr]{%
1	15.6924110537721\\
2	15.0603357416878\\
3	13.8245786953236\\
4	11.4527909887611\\
5	6.56659653725088\\
6	2.67421033906282\\
7	0.886968086039128\\
8	0.265232582763955\\
9	0.0824065838430515\\
10	0.0607646058861363\\
11	0.054282449311977\\
12	0.0266129213109444\\
13	0.0199848139700822\\
14	0.0163252200171428\\
15	0.0107407136576178\\
16	0.00893285384836976\\
17	0.00682934788408951\\
18	0.00556964674745719\\
19	0.00464134220955529\\
20	0.00370761313246779\\
21	0.00324648020393229\\
22	0.00250451047893929\\
23	0.00236491201041606\\
24	0.00190217568689019\\
25	0.00149929181461339\\
26	0.0013616554088354\\
27	0.00118139105373546\\
28	0.00098506863969418\\
29	0.000835288351797969\\
30	0.000727392327415357\\
31	0.000651820879772302\\
32	0.000582326747098598\\
33	0.000513242338886553\\
34	0.000452288178264501\\
35	0.000400126908834151\\
36	0.000355622947752988\\
37	0.00031782098561263\\
38	0.000285311989656859\\
39	0.000257432895860568\\
40	0.000232852943823615\\
41	0.000211142697888375\\
42	0.000191571772726683\\
43	0.000174098805956585\\
44	0.000158407643603325\\
45	0.000144446856093445\\
46	0.000131979530709394\\
47	0.000120902450173163\\
48	0.000111006944496932\\
49	0.000102171293146062\\
50	9.42096932881191e-05\\
51	8.70075752735131e-05\\
52	8.04489129273604e-05\\
53	7.44596708253171e-05\\
54	6.89983316699219e-05\\
55	6.40055478456916e-05\\
56	5.94681688321243e-05\\
57	5.53188800992416e-05\\
58	5.15551413740835e-05\\
59	4.81031944694403e-05\\
60	4.49601097089536e-05\\
61	4.20553565245448e-05\\
62	3.93902446588262e-05\\
63	3.69120554916681e-05\\
64	3.46283298275603e-05\\
65	3.25042231325751e-05\\
66	3.05464468772214e-05\\
67	2.87299082262576e-05\\
68	2.70553563759611e-05\\
69	2.55028208319478e-05\\
70	2.4066606379753e-05\\
71	2.2730933425788e-05\\
72	2.14863604118794e-05\\
73	2.03240456470493e-05\\
74	1.92350926640735e-05\\
75	1.8217855590934e-05\\
76	1.72637813849421e-05\\
77	1.63749191291035e-05\\
78	1.55413065528975e-05\\
79	1.47659972819786e-05\\
80	1.40370321702092e-05\\
81	1.33574794183827e-05\\
82	1.27147812492741e-05\\
83	1.21129608569088e-05\\
84	1.15414276817829e-05\\
85	1.10055167400623e-05\\
86	1.04968146624873e-05\\
87	1.00204946914088e-05\\
88	9.56913131068561e-06\\
89	9.14653215029623e-06\\
90	8.7456491877537e-06\\
91	8.36874314425623e-06\\
92	8.00980418920932e-06\\
93	7.6704897991047e-06\\
94	7.34683427807099e-06\\
95	7.04030647940867e-06\\
96	6.74857273405043e-06\\
97	6.47240655625682e-06\\
98	6.21023219875893e-06\\
99	5.96181390710457e-06\\
100	5.72591518973704e-06\\
101	5.50149044172138e-06\\
102	5.28781675001714e-06\\
103	5.08366701902631e-06\\
104	4.8890884879866e-06\\
105	4.70295940583457e-06\\
106	4.52583508664684e-06\\
107	4.35652724471307e-06\\
108	4.19571666682232e-06\\
109	4.04199337786176e-06\\
110	3.89596919886682e-06\\
111	3.75604743920286e-06\\
112	3.62281025847132e-06\\
113	3.49472738346071e-06\\
114	3.37250664653402e-06\\
115	3.25487112807028e-06\\
116	3.14263238108266e-06\\
117	3.03469917334382e-06\\
118	2.93184402955819e-06\\
119	2.83303028001933e-06\\
120	2.73889490290357e-06\\
121	2.64839386185873e-06\\
122	2.56202087868131e-06\\
123	2.47878333116685e-06\\
124	2.39912255449103e-06\\
125	2.32221645039872e-06\\
126	2.2485213134341e-06\\
127	2.17739735453441e-06\\
128	2.10927881148658e-06\\
129	2.04362819156135e-06\\
130	1.9807988851965e-06\\
131	1.92028447908037e-06\\
132	1.86232973353905e-06\\
133	1.80644158314793e-06\\
134	1.7527838681272e-06\\
135	1.70092868321193e-06\\
136	1.65102634072925e-06\\
137	1.60276236871942e-06\\
138	1.55629444100375e-06\\
139	1.51139757944181e-06\\
140	1.46820600911396e-06\\
141	1.42653244537341e-06\\
142	1.38645669887646e-06\\
143	1.34779437290137e-06\\
144	1.31056348128185e-06\\
145	1.27458614440177e-06\\
146	1.2398505900028e-06\\
147	1.2062239773105e-06\\
148	1.17370526110445e-06\\
149	1.14222222124241e-06\\
150	1.11178400430268e-06\\
151	1.08235424414448e-06\\
152	1.05392812482579e-06\\
153	1.02647503340836e-06\\
154	9.9995754768059e-07\\
155	9.74335659475905e-07\\
156	9.49541292795381e-07\\
157	9.25537398585953e-07\\
158	9.02253370756275e-07\\
159	8.79680719679134e-07\\
160	8.57766881835859e-07\\
161	8.36532828759357e-07\\
162	8.15935129971437e-07\\
163	7.96005094690613e-07\\
164	7.76690413714555e-07\\
165	7.58015513994089e-07\\
166	7.39908418910541e-07\\
167	7.22381983826441e-07\\
168	7.05351865825573e-07\\
169	6.88833778528441e-07\\
170	6.72752661596851e-07\\
171	6.57141955865867e-07\\
172	6.4194343812817e-07\\
173	6.27203280537175e-07\\
174	6.12869251476003e-07\\
175	5.98986808713354e-07\\
176	5.85496861054579e-07\\
177	5.72434491618105e-07\\
178	5.59728631639513e-07\\
179	5.47405177097004e-07\\
180	5.35390918049041e-07\\
181	5.23715732243491e-07\\
182	5.12318895832947e-07\\
183	5.01241152853258e-07\\
184	4.90434318019295e-07\\
185	4.79943357169069e-07\\
186	4.69722792893116e-07\\
187	4.59812629262669e-07\\
188	4.50161452947748e-07\\
189	4.40799258936077e-07\\
190	4.31667833754942e-07\\
191	4.22791577223524e-07\\
192	4.14115629301136e-07\\
193	4.0566903454943e-07\\
194	3.9740856109975e-07\\
195	3.89369523956866e-07\\
196	3.8151676089484e-07\\
197	3.73885585120798e-07\\
198	3.66441034555813e-07\\
199	3.59212261669537e-07\\
200	3.52159151624565e-07\\
201	3.45302794058379e-07\\
202	3.38600040605934e-07\\
203	3.32069793620708e-07\\
204	3.25674670096522e-07\\
205	3.1943812850756e-07\\
206	3.13331884256535e-07\\
207	3.07382591547354e-07\\
208	3.01566394005398e-07\\
209	2.9590794716236e-07\\
210	2.90381992999217e-07\\
211	2.85007163072197e-07\\
212	2.79754053809572e-07\\
213	2.74635658082891e-07\\
214	2.69622342503019e-07\\
215	2.64727426110051e-07\\
216	2.59927510477631e-07\\
217	2.55239733391659e-07\\
218	2.50646988333751e-07\\
219	2.46167575573067e-07\\
220	2.41786222471465e-07\\
221	2.37518370567314e-07\\
222	2.33346508822145e-07\\
223	2.29280833135668e-07\\
224	2.25300925528949e-07\\
225	2.21413763605311e-07\\
226	2.17600604839265e-07\\
227	2.13870236859165e-07\\
228	2.10209442563525e-07\\
229	2.06629778771019e-07\\
230	2.03121858887437e-07\\
231	1.99697101745639e-07\\
232	1.96346256760016e-07\\
233	1.93077612267592e-07\\
234	1.89879437408839e-07\\
235	1.8675588527793e-07\\
236	1.83693674394007e-07\\
237	1.80695826561331e-07\\
238	1.77751781258115e-07\\
239	1.74866927649937e-07\\
240	1.7203496821382e-07\\
241	1.69262989069495e-07\\
242	1.66546632572409e-07\\
243	1.6389204446724e-07\\
244	1.61293974969797e-07\\
245	1.5875551316249e-07\\
246	1.56269141979179e-07\\
247	1.53835117517814e-07\\
248	1.51445674513736e-07\\
249	1.49101496249028e-07\\
250	1.46797818965125e-07\\
251	1.44537569691817e-07\\
252	1.42318864455588e-07\\
253	1.40145395462209e-07\\
254	1.38015882831859e-07\\
255	1.35932601200541e-07\\
256	1.33892831305359e-07\\
257	1.31896121818247e-07\\
258	1.29938036163922e-07\\
259	1.28016648664179e-07\\
260	1.26128325903334e-07\\
261	1.24272508058062e-07\\
262	1.22448259578261e-07\\
263	1.20656759403134e-07\\
264	1.18898713955927e-07\\
265	1.17175340066319e-07\\
266	1.15487056890406e-07\\
267	1.13833413813852e-07\\
268	1.12213211425242e-07\\
269	1.10623846891004e-07\\
270	1.09063184273249e-07\\
271	1.07528397459425e-07\\
272	1.06018808171826e-07\\
273	1.04533282586286e-07\\
274	1.03073187029716e-07\\
275	1.01638476197515e-07\\
276	1.0023117731472e-07\\
277	9.88507485735132e-08\\
278	9.74983674593365e-08\\
279	9.61717763624852e-08\\
280	9.48706687385711e-08\\
281	9.35914306144385e-08\\
282	9.23336622360029e-08\\
283	9.10944859315379e-08\\
284	8.98751434685255e-08\\
285	8.86743009649587e-08\\
286	8.74944733320724e-08\\
287	8.63347943162105e-08\\
288	8.51977221238461e-08\\
289	8.40815278913704e-08\\
290	8.29875198216398e-08\\
291	8.19124368972027e-08\\
292	8.08565294804783e-08\\
293	7.98160991004336e-08\\
294	7.87920084911215e-08\\
295	7.77818013079493e-08\\
296	7.67877525898317e-08\\
297	7.5808514399778e-08\\
298	7.4846924969829e-08\\
299	7.39015935548185e-08\\
300	7.29748112052705e-08\\
301	7.20641284287889e-08\\
302	7.11706404084561e-08\\
303	7.02908106403779e-08\\
304	6.94252941556066e-08\\
305	6.85709472601778e-08\\
306	6.77294054159571e-08\\
307	6.68988038910791e-08\\
308	6.60817226179503e-08\\
309	6.52768985603525e-08\\
310	6.44869233501277e-08\\
311	6.37101230795692e-08\\
312	6.29482547583357e-08\\
313	6.21985937093378e-08\\
314	6.14619265624818e-08\\
315	6.07350732011914e-08\\
316	6.00189903024305e-08\\
317	5.93113865257715e-08\\
318	5.86141845362807e-08\\
319	5.79260806536455e-08\\
320	5.72494508065296e-08\\
321	5.65831590456527e-08\\
322	5.59292256785869e-08\\
323	5.52858898728781e-08\\
324	5.46542545854401e-08\\
325	5.40317325811371e-08\\
326	5.34188739578044e-08\\
327	5.28132407952948e-08\\
328	5.22159322215862e-08\\
329	5.16254930758258e-08\\
330	5.10437364230086e-08\\
331	5.04698115589712e-08\\
332	4.99056108221291e-08\\
333	4.93501507249485e-08\\
334	4.88047694450566e-08\\
335	4.82677994889619e-08\\
336	4.77397729019168e-08\\
337	4.7218550131255e-08\\
338	4.67045475306453e-08\\
339	4.61961465923027e-08\\
340	4.56944098298172e-08\\
341	4.51985310072199e-08\\
342	4.47099977673175e-08\\
343	4.42282886404469e-08\\
344	4.37547438593804e-08\\
345	4.32885378558294e-08\\
346	4.28303952883924e-08\\
347	4.2378872702886e-08\\
348	4.19341116547659e-08\\
349	4.14945360448553e-08\\
350	4.10604751626733e-08\\
351	4.06309900761888e-08\\
352	4.02069803162013e-08\\
353	3.9788082321588e-08\\
354	3.93754034773009e-08\\
355	3.89686543545161e-08\\
356	3.85686824821458e-08\\
357	3.81748212353549e-08\\
358	3.7787336092512e-08\\
359	3.74050711360272e-08\\
360	3.70279441234572e-08\\
361	3.66548957578859e-08\\
362	3.62861776372522e-08\\
363	3.59213254990456e-08\\
364	3.55610372038017e-08\\
365	3.52052244869547e-08\\
366	3.48546611391122e-08\\
367	3.45092130911334e-08\\
368	3.4169355188201e-08\\
369	3.3834566776569e-08\\
370	3.35048184883463e-08\\
371	3.31792320632376e-08\\
372	3.28576107090226e-08\\
373	3.25392883677737e-08\\
374	3.22244396218597e-08\\
375	3.19128991980537e-08\\
376	3.16051792447684e-08\\
377	3.13013544199482e-08\\
378	3.10019434661251e-08\\
379	3.07069166820979e-08\\
380	3.04164966908388e-08\\
381	3.01302911602305e-08\\
382	2.9848105079973e-08\\
383	2.9569289776031e-08\\
384	2.92935956785121e-08\\
385	2.90206212883237e-08\\
386	2.87504736828272e-08\\
387	2.84831580988015e-08\\
388	2.82190421003176e-08\\
389	2.79582874537154e-08\\
390	2.77012387448086e-08\\
391	2.74479334882511e-08\\
392	2.71984443997059e-08\\
393	2.69524837033533e-08\\
394	2.67097792441576e-08\\
395	2.64698513437317e-08\\
396	2.62324286242007e-08\\
397	2.59972707345516e-08\\
398	2.57644291444532e-08\\
399	2.55339964084927e-08\\
400	2.53062353989184e-08\\
401	2.5081349438976e-08\\
402	2.48595725928146e-08\\
403	2.46409869151276e-08\\
404	2.44255895699928e-08\\
405	2.42131765205684e-08\\
406	2.4003454915439e-08\\
407	2.37960607980169e-08\\
408	2.3590717481316e-08\\
409	2.33872703197653e-08\\
410	2.31857278659752e-08\\
411	2.29862151348035e-08\\
412	2.27889241333369e-08\\
413	2.25940694043768e-08\\
414	2.24018224209101e-08\\
415	2.22122935106796e-08\\
416	2.20254561399627e-08\\
417	2.18411731601079e-08\\
418	2.16591660930368e-08\\
419	2.14791468960016e-08\\
420	2.13008413412769e-08\\
421	2.112413012938e-08\\
422	2.09489866189654e-08\\
423	2.07755345124465e-08\\
424	2.06039168516899e-08\\
425	2.0434341760229e-08\\
426	2.02669482170566e-08\\
427	2.0101864872689e-08\\
428	1.99390739221765e-08\\
429	1.97784944528877e-08\\
430	1.96198859380439e-08\\
431	1.9463012124156e-08\\
432	1.93076070076076e-08\\
433	1.91535537460206e-08\\
434	1.90007948044611e-08\\
435	1.8849429419482e-08\\
436	1.86995650090015e-08\\
437	1.85513904230243e-08\\
438	1.84050299925895e-08\\
439	1.82606218700024e-08\\
440	1.81181756784166e-08\\
441	1.79776601151667e-08\\
442	1.78388890894997e-08\\
443	1.77016684076017e-08\\
444	1.75657471228076e-08\\
445	1.7430993016296e-08\\
446	1.7297317675473e-08\\
447	1.71647803220134e-08\\
448	1.70334569282151e-08\\
449	1.69035080260569e-08\\
450	1.67750498072955e-08\\
451	1.66482223237886e-08\\
452	1.65230684996125e-08\\
453	1.63996019242676e-08\\
454	1.62777028538483e-08\\
455	1.61572217185351e-08\\
456	1.60379348726296e-08\\
457	1.59196899773643e-08\\
458	1.58023679798552e-08\\
459	1.56859760567304e-08\\
460	1.55705566585214e-08\\
461	1.54562328443986e-08\\
462	1.53431116926903e-08\\
463	1.52313261959344e-08\\
464	1.5120949075186e-08\\
465	1.50120312178805e-08\\
466	1.49045239353019e-08\\
467	1.4798331330569e-08\\
468	1.46932762413827e-08\\
469	1.45891974389866e-08\\
470	1.44859499230154e-08\\
471	1.43834820855311e-08\\
472	1.4281795033604e-08\\
473	1.41809648511936e-08\\
474	1.40810819230919e-08\\
475	1.39822620868183e-08\\
476	1.38845967853995e-08\\
477	1.37881618341031e-08\\
478	1.36929732941784e-08\\
479	1.3598995567088e-08\\
480	1.35061196653622e-08\\
481	1.34142026598741e-08\\
482	1.33230926369302e-08\\
483	1.32326828975346e-08\\
484	1.31429253130932e-08\\
485	1.30538368683338e-08\\
486	1.2965478641928e-08\\
487	1.2877937555714e-08\\
488	1.27913082957178e-08\\
489	1.2705674661806e-08\\
490	1.26210993712697e-08\\
491	1.25376004175111e-08\\
492	1.24551483034994e-08\\
493	1.23736517286818e-08\\
494	1.22929889506751e-08\\
495	1.2213024818384e-08\\
496	1.21336657055218e-08\\
497	1.20548601622712e-08\\
498	1.19766224338706e-08\\
499	1.1898995203615e-08\\
500	1.18220583131973e-08\\
501	1.17458856004796e-08\\
502	1.16705610043129e-08\\
503	1.15961350593342e-08\\
504	1.15226447621935e-08\\
505	1.14500664290403e-08\\
506	1.13783478505083e-08\\
507	1.1307377637163e-08\\
508	1.12370517372702e-08\\
509	1.11672616634932e-08\\
510	1.1097962534436e-08\\
511	1.10291337837728e-08\\
512	1.09608178576227e-08\\
513	1.08930588624205e-08\\
514	1.08259369687409e-08\\
515	1.07595081651593e-08\\
516	1.06938443976631e-08\\
517	1.06289717995514e-08\\
518	1.05649130685588e-08\\
519	1.05016244395896e-08\\
520	1.04390520755105e-08\\
521	1.03770868660043e-08\\
522	1.03156485624443e-08\\
523	1.02546499658537e-08\\
524	1.01940727081916e-08\\
525	1.0133906143097e-08\\
526	1.00741996335686e-08\\
527	1.00149894209555e-08\\
528	9.95634931237708e-09\\
529	9.89832072202557e-09\\
530	9.84096711966333e-09\\
531	9.78430278990148e-09\\
532	9.72834882508443e-09\\
533	9.67306077170818e-09\\
534	9.61839674841341e-09\\
535	9.56425981443017e-09\\
536	9.51058781839984e-09\\
537	9.4573031979365e-09\\
538	9.40439487111839e-09\\
539	9.35184857758082e-09\\
540	9.29970884787268e-09\\
541	9.24800091804348e-09\\
542	9.19679004544666e-09\\
543	9.14610745959753e-09\\
544	9.09601218773472e-09\\
545	9.04651699888497e-09\\
546	8.9976511424524e-09\\
547	8.94938419910246e-09\\
548	8.90169577042346e-09\\
549	8.85450857627972e-09\\
550	8.80777496301533e-09\\
551	8.76141984662745e-09\\
552	8.71542671216321e-09\\
553	8.66976906594313e-09\\
554	8.62447849204746e-09\\
555	8.57956732671521e-09\\
556	8.53509069165024e-09\\
557	8.4910725997184e-09\\
558	8.44756920572963e-09\\
559	8.40459761116069e-09\\
560	8.36219947931568e-09\\
561	8.32036864877278e-09\\
562	8.27911614482305e-09\\
563	8.23839964292191e-09\\
564	8.19819739873047e-09\\
565	8.15844712208585e-09\\
566	8.11912916924326e-09\\
567	8.08020985633795e-09\\
568	8.04172138201318e-09\\
569	8.00369819649696e-09\\
570	7.96622724670239e-09\\
571	7.92932416809167e-09\\
572	7.89290828106937e-09\\
573	7.85674400352928e-09\\
574	7.82074982176181e-09\\
575	7.78552048945447e-09\\
576	7.75392747727925e-09\\
577	7.72809110779008e-09\\
578	7.69779003865014e-09\\
579	7.6622517645447e-09\\
580	7.62503560519558e-09\\
581	7.58723153693095e-09\\
582	7.54930037663815e-09\\
583	7.54324343152241e-09\\
584	7.51081051264872e-09\\
585	7.47244675574685e-09\\
586	7.43400101312069e-09\\
587	7.39553490678549e-09\\
588	7.35710849529428e-09\\
589	7.31878375467177e-09\\
590	7.28062202056516e-09\\
591	7.26173992818739e-09\\
592	7.24268041559784e-09\\
593	7.20500763217652e-09\\
594	7.16763777076308e-09\\
595	7.13058304394653e-09\\
596	7.09382785876087e-09\\
597	7.05732934788942e-09\\
598	7.02102814285292e-09\\
599	6.98486686571e-09\\
600	6.94880741951766e-09\\
601	6.92861750270097e-09\\
602	6.91283908170614e-09\\
603	6.87697653944451e-09\\
604	6.84125246969815e-09\\
605	6.80570959202355e-09\\
606	6.77039439292252e-09\\
607	6.73535217606648e-09\\
608	6.700622173685e-09\\
609	6.66623131496956e-09\\
610	6.63218662871704e-09\\
611	6.59846871994063e-09\\
612	6.56503151268415e-09\\
613	6.55804667401288e-09\\
614	6.53181237456659e-09\\
615	6.49875063403293e-09\\
616	6.46580509059551e-09\\
617	6.43296213714045e-09\\
618	6.40023348428296e-09\\
619	6.36764844864641e-09\\
620	6.33524599713599e-09\\
621	6.30306859473691e-09\\
622	6.27115765625427e-09\\
623	6.23954912936732e-09\\
624	6.20826775294514e-09\\
625	6.17731988838515e-09\\
626	6.15816654009471e-09\\
627	6.14668702671144e-09\\
628	6.11632510305314e-09\\
629	6.08617379200974e-09\\
630	6.05617415682987e-09\\
631	6.02628568636475e-09\\
632	5.99649418516534e-09\\
633	5.96680955828439e-09\\
634	5.93725830450613e-09\\
635	5.90787592660563e-09\\
636	5.87870118887036e-09\\
637	5.84977205967269e-09\\
638	5.8211219240138e-09\\
639	5.79277459172191e-09\\
640	5.76473774449923e-09\\
641	5.73699654020109e-09\\
642	5.73515274242425e-09\\
643	5.70951186889152e-09\\
644	5.68222792074092e-09\\
645	5.65508810400606e-09\\
646	5.62805162754955e-09\\
647	5.60110200862464e-09\\
648	5.57424568260767e-09\\
649	5.54750524178279e-09\\
650	5.52091211102263e-09\\
651	5.49450119071958e-09\\
652	5.46830720332815e-09\\
653	5.44236152417794e-09\\
654	5.41668814801781e-09\\
655	5.39129804828583e-09\\
656	5.36618288724014e-09\\
657	5.34131172167535e-09\\
658	5.31663536860604e-09\\
659	5.29427087041026e-09\\
660	5.29209962125948e-09\\
661	5.26766137754251e-09\\
662	5.2432988632785e-09\\
663	5.21901246517892e-09\\
664	5.19481910241983e-09\\
665	5.17074531652433e-09\\
666	5.14682190068524e-09\\
667	5.12308037886367e-09\\
668	5.09955038861842e-09\\
669	5.07625666682305e-09\\
670	5.05321463260519e-09\\
671	5.03042452862589e-09\\
672	5.00786657223624e-09\\
673	4.98550124897513e-09\\
674	4.96327815302286e-09\\
675	4.9411507831786e-09\\
676	4.91908922749298e-09\\
677	4.89708467721691e-09\\
678	4.87514618927871e-09\\
679	4.85329432432232e-09\\
680	4.84031215506274e-09\\
681	4.83155545512563e-09\\
682	4.80995795591981e-09\\
683	4.78852977496229e-09\\
684	4.76729627263628e-09\\
685	4.74627708300712e-09\\
686	4.72548142231107e-09\\
687	4.70490247155244e-09\\
688	4.68451392910185e-09\\
689	4.66427316363778e-09\\
690	4.6441323106417e-09\\
691	4.6240525140235e-09\\
692	4.60401325588573e-09\\
693	4.58401331226786e-09\\
694	4.56406595846086e-09\\
695	4.54419288465815e-09\\
696	4.52441965898547e-09\\
697	4.5047728623881e-09\\
698	4.48527827527183e-09\\
699	4.46595879969147e-09\\
700	4.44683117272107e-09\\
701	4.42790092481844e-09\\
702	4.40915697996442e-09\\
703	4.39056939501931e-09\\
704	4.37781473260402e-09\\
705	4.37209461260776e-09\\
706	4.35368796292022e-09\\
707	4.33531703998367e-09\\
708	4.31696844645058e-09\\
709	4.29864643739492e-09\\
710	4.28036730508003e-09\\
711	4.2621537852325e-09\\
712	4.24403126318656e-09\\
713	4.22602556635535e-09\\
714	4.20816143540471e-09\\
715	4.19046040300696e-09\\
716	4.17293717665394e-09\\
717	4.15559425886752e-09\\
718	4.13841672245749e-09\\
719	4.12137132450564e-09\\
720	4.10441387392753e-09\\
721	4.08750302766506e-09\\
722	4.07061248080794e-09\\
723	4.0537351419589e-09\\
724	4.03687980279241e-09\\
725	4.02006507343775e-09\\
726	4.00331430017928e-09\\
727	3.98665245894153e-09\\
728	3.97010447823084e-09\\
729	3.95369388734084e-09\\
730	3.93744056015058e-09\\
731	3.92135667007892e-09\\
732	3.91119181673598e-09\\
733	3.90544101174138e-09\\
734	3.88967415823218e-09\\
735	3.87401938273954e-09\\
736	3.85843235147187e-09\\
737	3.84287564524676e-09\\
738	3.8273291253641e-09\\
739	3.81179133136312e-09\\
740	3.79627458782069e-09\\
741	3.78079888168092e-09\\
742	3.76538740578852e-09\\
743	3.75006422330666e-09\\
744	3.73485299696687e-09\\
745	3.71977576045896e-09\\
746	3.70485041733677e-09\\
747	3.69008628568204e-09\\
748	3.67547815072761e-09\\
749	3.66100211972346e-09\\
750	3.64661861944636e-09\\
751	3.63228418126719e-09\\
752	3.61796597737975e-09\\
753	3.60364981756393e-09\\
754	3.58933897443161e-09\\
755	3.57504829796403e-09\\
756	3.56079831376154e-09\\
757	3.54661165549843e-09\\
758	3.53251127283119e-09\\
759	3.51851958122241e-09\\
760	3.50465717226159e-09\\
761	3.49094011814049e-09\\
762	3.47737503250468e-09\\
763	3.46395291859532e-09\\
764	3.45064594042895e-09\\
765	3.44481727290318e-09\\
766	3.43741248904477e-09\\
767	3.4242108217639e-09\\
768	3.4110129176294e-09\\
769	3.39780996470873e-09\\
770	3.38460907940166e-09\\
771	3.37142683601344e-09\\
772	3.35828393220501e-09\\
773	3.34520229835637e-09\\
774	3.33220391610559e-09\\
775	3.31931013688105e-09\\
776	3.30654044850966e-09\\
777	3.29390937191648e-09\\
778	3.28142105825716e-09\\
779	3.26906293743304e-09\\
780	3.2568035660471e-09\\
781	3.24459987551378e-09\\
782	3.23241235330534e-09\\
783	3.22021783718652e-09\\
784	3.20801239536928e-09\\
785	3.19580632888091e-09\\
786	3.1836174395698e-09\\
787	3.17146638265545e-09\\
788	3.15937445670397e-09\\
789	3.14736283024904e-09\\
790	3.13545207786543e-09\\
791	3.12366081712054e-09\\
792	3.1120021824701e-09\\
793	3.10047780194983e-09\\
794	3.08907132345489e-09\\
795	3.07774750817075e-09\\
796	3.06646213302439e-09\\
797	3.05517840702507e-09\\
798	3.04387819562938e-09\\
799	3.0325622126846e-09\\
800	3.02124347715149e-09\\
801	3.00994075559419e-09\\
802	2.99867457505164e-09\\
803	2.98746563278993e-09\\
804	2.98308825692027e-09\\
805	2.97633450729745e-09\\
806	2.96530125172775e-09\\
807	2.95438378701162e-09\\
808	2.94359380227779e-09\\
809	2.93292987411063e-09\\
810	2.92237091170222e-09\\
811	2.91187731709324e-09\\
812	2.90140401450489e-09\\
813	2.89091791018555e-09\\
814	2.88040665552633e-09\\
815	2.86987575438935e-09\\
816	2.85934085742722e-09\\
817	2.84882158246694e-09\\
818	2.83833839372563e-09\\
819	2.82791169690354e-09\\
820	2.81756183602291e-09\\
821	2.80730872433563e-09\\
822	2.79716975696757e-09\\
823	2.78715466277644e-09\\
824	2.77725747901006e-09\\
825	2.76745039139017e-09\\
826	2.75768849336004e-09\\
827	2.74792706334069e-09\\
828	2.73813894391398e-09\\
829	2.72831916357435e-09\\
830	2.71847857501004e-09\\
831	2.70863547959922e-09\\
832	2.69881028982439e-09\\
833	2.68902351941013e-09\\
834	2.6792955403049e-09\\
835	2.66964684371655e-09\\
836	2.66009742867617e-09\\
837	2.65066374563513e-09\\
838	2.64135176023098e-09\\
839	2.63214753317278e-09\\
840	2.62301336236073e-09\\
841	2.61389903160349e-09\\
842	2.60476362478846e-09\\
843	2.59558968864272e-09\\
844	2.58638136184174e-09\\
845	2.57715482400073e-09\\
846	2.56793042709626e-09\\
847	2.55872895789998e-09\\
848	2.54957087414548e-09\\
849	2.54047671963777e-09\\
850	2.53146745114966e-09\\
851	2.5304817483533e-09\\
852	2.52256303213032e-09\\
853	2.51377737227947e-09\\
854	2.50510845142738e-09\\
855	2.49652857055653e-09\\
856	2.48798733940057e-09\\
857	2.47943313390982e-09\\
858	2.47083592282028e-09\\
859	2.46219233922835e-09\\
860	2.45351624764919e-09\\
861	2.44482820905218e-09\\
862	2.4361497809433e-09\\
863	2.42750200028239e-09\\
864	2.41890577148065e-09\\
865	2.41038275399149e-09\\
866	2.40195509013014e-09\\
867	2.39364183775598e-09\\
868	2.38544958835542e-09\\
869	2.37735904737251e-09\\
870	2.36932104573145e-09\\
871	2.36127586160532e-09\\
872	2.3531832065119e-09\\
873	2.34503404768204e-09\\
874	2.33684140695407e-09\\
875	2.32862705558466e-09\\
876	2.3204139592833e-09\\
877	2.31222412281368e-09\\
878	2.30407918126624e-09\\
879	2.29600169485835e-09\\
880	2.28801561252313e-09\\
881	2.28014349448171e-09\\
882	2.27239686808369e-09\\
883	2.26475978356636e-09\\
884	2.25718017101548e-09\\
885	2.24959023558244e-09\\
886	2.24194364467777e-09\\
887	2.2342308652368e-09\\
888	2.22646735228502e-09\\
889	2.21867727423122e-09\\
890	2.21088511578619e-09\\
891	2.20311377172079e-09\\
892	2.19538568925564e-09\\
893	2.18772471238104e-09\\
894	2.18015700702982e-09\\
895	2.1727080742581e-09\\
896	2.16539063942679e-09\\
897	2.15818345929313e-09\\
898	2.15102208789635e-09\\
899	2.14382960994871e-09\\
900	2.13656211989744e-09\\
901	2.12921871651571e-09\\
902	2.12182179978869e-09\\
903	2.11439897428044e-09\\
904	2.10697607573067e-09\\
905	2.0995769121507e-09\\
906	2.09222544836808e-09\\
907	2.09161091831099e-09\\
908	2.08494839318946e-09\\
909	2.07777601689304e-09\\
910	2.07073638710096e-09\\
911	2.06383528878188e-09\\
912	2.05702789635888e-09\\
913	2.05022398882224e-09\\
914	2.04334619778433e-09\\
915	2.03637226512758e-09\\
916	2.02932115270048e-09\\
917	2.02222443885365e-09\\
918	2.01511223947307e-09\\
919	2.00801104335933e-09\\
920	2.00094624605611e-09\\
921	1.99394567939036e-09\\
922	1.98704269719709e-09\\
923	1.98027416045041e-09\\
924	1.97366244907093e-09\\
925	1.96717595984254e-09\\
926	1.96071349854422e-09\\
927	1.95417103583845e-09\\
928	1.94751304196252e-09\\
929	1.94076102476143e-09\\
930	1.9339528015322e-09\\
931	1.92712287803537e-09\\
932	1.92030017564881e-09\\
933	1.91351164393966e-09\\
934	1.90678729337129e-09\\
935	1.90016501531354e-09\\
936	1.89369021350006e-09\\
937	1.8873935051835e-09\\
938	1.88123357201887e-09\\
939	1.8750766893971e-09\\
940	1.86880156948581e-09\\
941	1.86238532601783e-09\\
942	1.85586593376314e-09\\
943	1.8492891566979e-09\\
944	1.84269205662295e-09\\
945	1.83610488746733e-09\\
946	1.82955745482853e-09\\
947	1.82308643162487e-09\\
948	1.81674256636459e-09\\
949	1.81058640300855e-09\\
950	1.80463653820029e-09\\
951	1.79877745205604e-09\\
952	1.79281623883767e-09\\
953	1.78667216014584e-09\\
954	1.78038290868583e-09\\
955	1.77401029348534e-09\\
956	1.76760267823325e-09\\
957	1.76119597488624e-09\\
958	1.75482251024707e-09\\
959	1.74852110734493e-09\\
960	1.74234928923528e-09\\
961	1.73638910589555e-09\\
962	1.73069276456632e-09\\
963	1.72512728455009e-09\\
964	1.71942362222861e-09\\
965	1.71348418813767e-09\\
966	1.70737527016242e-09\\
967	1.70117682119301e-09\\
968	1.69494271843642e-09\\
969	1.6887101658512e-09\\
970	1.68251347084154e-09\\
971	1.67639889953402e-09\\
972	1.67129135777513e-09\\
973	1.67044596495075e-09\\
974	1.66477921878437e-09\\
975	1.65944002022186e-09\\
976	1.65412220741016e-09\\
977	1.64850087106689e-09\\
978	1.64260335641254e-09\\
979	1.63655682866068e-09\\
980	1.63044445089809e-09\\
981	1.62431477067777e-09\\
982	1.61820478544e-09\\
983	1.61215968244592e-09\\
984	1.60626260634943e-09\\
985	1.60068438955972e-09\\
986	1.59559046524593e-09\\
987	1.59059436391116e-09\\
988	1.58515634644616e-09\\
989	1.57935831858032e-09\\
990	1.57339605818364e-09\\
991	1.56736968287682e-09\\
992	1.56132938738511e-09\\
993	1.55531235145424e-09\\
994	1.54937006272621e-09\\
995	1.54362183302545e-09\\
996	1.5383777684542e-09\\
997	1.53381113071938e-09\\
998	1.52893051451389e-09\\
999	1.52338072782566e-09\\
1000	1.51753868256572e-09\\
1001	1.51159354674215e-09\\
1002	1.50561542114239e-09\\
1003	1.49964099096394e-09\\
1004	1.49370980811828e-09\\
1005	1.48791628931935e-09\\
1006	1.48262220184991e-09\\
1007	1.47848835482184e-09\\
1008	1.47399258154319e-09\\
1009	1.4684623733557e-09\\
1010	1.46260770703794e-09\\
1011	1.45665756141143e-09\\
1012	1.45067576739789e-09\\
1013	1.44469283265837e-09\\
1014	1.43874580557984e-09\\
1015	1.43295750845929e-09\\
1016	1.42813373886796e-09\\
1017	1.42487146901e-09\\
1018	1.41980962639836e-09\\
1019	1.41398375742634e-09\\
1020	1.40802775524447e-09\\
1021	1.40203868114769e-09\\
1022	1.39604476439704e-09\\
1023	1.39006587385751e-09\\
1024	1.38415418389728e-09\\
1025	1.37886303614833e-09\\
1026	1.37644868930094e-09\\
1027	1.37159890041411e-09\\
1028	1.36571544678772e-09\\
1029	1.35973486738453e-09\\
1030	1.35372540138267e-09\\
1031	1.34770107974742e-09\\
1032	1.34167074866275e-09\\
1033	1.33566659265483e-09\\
1034	1.33068026598755e-09\\
1035	1.32898999169539e-09\\
1036	1.32329422609284e-09\\
1037	1.31723064672704e-09\\
1038	1.31112863128592e-09\\
1039	1.30500814941153e-09\\
1040	1.29887496337078e-09\\
1041	1.29273537939438e-09\\
1042	1.28666365023197e-09\\
1043	1.284893587736e-09\\
1044	1.28032543147447e-09\\
1045	1.27462383620498e-09\\
1046	1.27415726410293e-09\\
1047	1.26795586833586e-09\\
1048	1.26173439985219e-09\\
1049	1.25549472169115e-09\\
1050	1.24923837683439e-09\\
1051	1.24298638609141e-09\\
1052	1.2413072136879e-09\\
1053	1.23664488568309e-09\\
1054	1.23033414109187e-09\\
1055	1.22399560422147e-09\\
1056	1.2176287622529e-09\\
1057	1.21122781543657e-09\\
1058	1.20478522807455e-09\\
1059	1.19917936460424e-09\\
1060	1.19827252388325e-09\\
1061	1.19172858900108e-09\\
1062	1.1851043738431e-09\\
1063	1.17840744176211e-09\\
1064	1.17163787271897e-09\\
1065	1.16479853867858e-09\\
1066	1.15842445193972e-09\\
1067	1.15788748206128e-09\\
1068	1.1509249807542e-09\\
1069	1.14389667886823e-09\\
1070	1.13680362475922e-09\\
1071	1.12963645450279e-09\\
1072	1.12237964258428e-09\\
1073	1.11904367182438e-09\\
1074	1.11501199282723e-09\\
1075	1.10751029114763e-09\\
1076	1.09984737508473e-09\\
1077	1.09199048990156e-09\\
1078	1.08389110249195e-09\\
1079	1.08099309427987e-09\\
1080	1.0754694465718e-09\\
1081	1.06660126523579e-09\\
1082	1.05712291100442e-09\\
1083	1.04686314466696e-09\\
1084	1.0442232857634e-09\\
1085	1.03566578939053e-09\\
1086	1.02333061821045e-09\\
1087	1.00923036779796e-09\\
1088	1.00869143304008e-09\\
1089	9.74356147949135e-10\\
1090	9.41177397377521e-10\\
1091	9.09116490153809e-10\\
1092	9.07093002661229e-10\\
1093	8.78136022288041e-10\\
1094	8.48199839259184e-10\\
1095	8.19272991313883e-10\\
1096	7.91321697235153e-10\\
1097	7.64313297599244e-10\\
1098	7.38216223789167e-10\\
1099	7.12999962969421e-10\\
1100	6.88635015803185e-10\\
1101	6.65092865381191e-10\\
1102	6.42345947976625e-10\\
1103	6.20367618899836e-10\\
1104	5.99132121288808e-10\\
1105	5.78614556630461e-10\\
1106	5.74610943376329e-10\\
1107	5.58790857320058e-10\\
1108	5.39637758033969e-10\\
1109	5.21132770635571e-10\\
1110	5.03254156212497e-10\\
1111	4.8598090153233e-10\\
1112	4.69292694654171e-10\\
1113	4.5316990140838e-10\\
1114	4.37593543148341e-10\\
1115	4.22545272980244e-10\\
1116	4.0800735842111e-10\\
1117	3.93962658594571e-10\\
1118	3.80394604264728e-10\\
1119	3.67287179608875e-10\\
1120	3.54624903531193e-10\\
1121	3.4239281208419e-10\\
1122	3.30576440770125e-10\\
1123	3.1916180737983e-10\\
1124	3.0813539819183e-10\\
1125	2.97484149988026e-10\\
1126	2.87195435678972e-10\\
1127	2.82795928807375e-10\\
1128	2.77257050361863e-10\\
1129	2.67657197250522e-10\\
1130	2.5838447262896e-10\\
1131	2.49427854980376e-10\\
1132	2.40776690937256e-10\\
1133	2.32420682843451e-10\\
1134	2.24349877511293e-10\\
1135	2.16554654892825e-10\\
1136	2.09025716281282e-10\\
1137	2.01754074150851e-10\\
1138	1.94731041974702e-10\\
1139	1.87948223500906e-10\\
1140	1.81397503726832e-10\\
1141	1.75071039817366e-10\\
1142	1.68961250857609e-10\\
1143	1.6306081035562e-10\\
1144	1.57362636970575e-10\\
1145	1.51859887001232e-10\\
1146	1.46545946196618e-10\\
1147	1.41414421830697e-10\\
1148	1.36459135928766e-10\\
1149	1.31674117481166e-10\\
1150	1.27053596049306e-10\\
1151	1.22591995021934e-10\\
1152	1.18283924656707e-10\\
1153	1.1412417681974e-10\\
1154	1.10107717860467e-10\\
1155	1.06229683478566e-10\\
1156	1.02485372910312e-10\\
1157	9.8870243607247e-11\\
1158	9.5379905566105e-11\\
1159	9.2010116830073e-11\\
1160	8.87567779288179e-11\\
1161	8.56159275397503e-11\\
1162	8.25837378434796e-11\\
1163	7.96565099952816e-11\\
1164	7.68306694359823e-11\\
1165	7.41027625248755e-11\\
1166	7.14694520368328e-11\\
1167	6.89275131657119e-11\\
1168	6.64738300341013e-11\\
1169	6.41053917669314e-11\\
1170	6.18192893069534e-11\\
1171	5.96127118620729e-11\\
1172	5.74829434470618e-11\\
1173	5.54273599178711e-11\\
1174	5.34434257268197e-11\\
1175	5.15286911854849e-11\\
1176	4.96807893554093e-11\\
1177	4.78974332721294e-11\\
1178	4.61764134194774e-11\\
1179	4.45155950594267e-11\\
1180	4.29129155503609e-11\\
1181	4.13663821279136e-11\\
1182	3.98740695366746e-11\\
1183	3.84341175568823e-11\\
1184	3.70447289849691e-11\\
1185	3.57041674613607e-11\\
1186	3.44107554784399e-11\\
1187	3.31628722144959e-11\\
1188	3.19589518102734e-11\\
1189	3.07974813506187e-11\\
1190	2.96769991762641e-11\\
1191	2.85960931248509e-11\\
1192	2.75533988268064e-11\\
1193	2.6547598114782e-11\\
1194	2.5577417478004e-11\\
1195	2.46416264296456e-11\\
1196	2.37390361705475e-11\\
1197	2.2868498163027e-11\\
1198	2.202890273017e-11\\
1199	2.12191777378902e-11\\
1200	2.04382873201076e-11\\
1201	1.96852306985548e-11\\
1202	1.89590409274225e-11\\
1203	1.82587838275746e-11\\
1204	1.75835567463336e-11\\
1205	1.69324876961261e-11\\
1206	1.63047340833662e-11\\
1207	1.56994819462862e-11\\
1208	1.51159448189579e-11\\
1209	1.45533628704991e-11\\
1210	1.40110020229256e-11\\
1211	1.34881530965807e-11\\
1212	1.29841309492324e-11\\
1213	1.24982736880873e-11\\
1214	1.20299418883048e-11\\
1215	1.15785178273557e-11\\
1216	1.11434048540906e-11\\
1217	1.07240265600904e-11\\
1218	1.03198262615239e-11\\
1219	9.93026622066569e-12\\
1220	9.55482708543443e-12\\
1221	9.19300732339418e-12\\
1222	8.84432257009314e-12\\
1223	8.50830513818033e-12\\
1224	8.18450341806842e-12\\
1225	7.87248145014399e-12\\
1226	7.57181832200075e-12\\
1227	7.28210775513525e-12\\
1228	7.00295763530573e-12\\
1229	6.73398954021426e-12\\
1230	6.47483831312467e-12\\
1231	6.22515167524047e-12\\
1232	5.98458978052971e-12\\
1233	5.75282485513842e-12\\
1234	5.52954080951465e-12\\
1235	5.31443286681107e-12\\
1236	5.10720722358588e-12\\
1237	4.90758069137245e-12\\
1238	4.71528038317468e-12\\
1239	4.53004334970941e-12\\
1240	4.35161631961363e-12\\
1241	4.17975533132573e-12\\
1242	4.0142254898548e-12\\
1243	3.85480060632049e-12\\
1244	3.701262973174e-12\\
1245	3.55340302806995e-12\\
1246	3.41101911497479e-12\\
1247	3.27391718573893e-12\\
1248	3.14191053804617e-12\\
1249	3.0148195405569e-12\\
1250	2.89247139358612e-12\\
1251	2.77469983740444e-12\\
1252	2.66134492935596e-12\\
1253	2.55225277286657e-12\\
1254	2.44727527934032e-12\\
1255	2.34626992413123e-12\\
1256	2.24909950302055e-12\\
1257	2.15563191047408e-12\\
1258	2.06573989816158e-12\\
1259	1.97930087207072e-12\\
1260	1.89619665864794e-12\\
1261	1.81631330144544e-12\\
1262	1.73954086571985e-12\\
1263	1.6657732354014e-12\\
1264	1.59490793282098e-12\\
1265	1.52684593503633e-12\\
1266	1.46149151044781e-12\\
1267	1.39875204671876e-12\\
1268	1.33853791450906e-12\\
1269	1.28076231169846e-12\\
1270	1.22534112133003e-12\\
1271	1.17219280269803e-12\\
1272	1.12123825764221e-12\\
1273	1.07240073298376e-12\\
1274	1.0256057090689e-12\\
1275	9.80780795125442e-13\\
1276	9.37855675345191e-13\\
1277	8.96761982209136e-13\\
1278	8.57433253498041e-13\\
1279	8.19804846936519e-13\\
1280	7.83813874597425e-13\\
1281	7.49399144333827e-13\\
1282	7.16501091909914e-13\\
1283	6.85061734343274e-13\\
1284	6.55024591942155e-13\\
1285	6.26334650236335e-13\\
1286	5.98938288362299e-13\\
1287	5.72783214142279e-13\\
1288	5.47818388454745e-13\\
1289	5.23993957104414e-13\\
1290	5.01261133045358e-13\\
1291	4.7957210221154e-13\\
1292	4.58879881862319e-13\\
1293	4.39138135348116e-13\\
1294	4.20300968304005e-13\\
1295	4.02322667974157e-13\\
1296	3.85157306297533e-13\\
1297	3.68758379085479e-13\\
1298	3.53078138811437e-13\\
1299	3.38066991377343e-13\\
1300	3.23579698273376e-13\\
1301	3.09842959343509e-13\\
1302	2.96509776781029e-13\\
1303	2.83611289162993e-13\\
1304	2.71071542234112e-13\\
1305	2.58822579711711e-13\\
1306	2.46784674059389e-13\\
1307	2.34885980716745e-13\\
1308	2.23064823610546e-13\\
1309	2.1127358757621e-13\\
1310	1.99482428383657e-13\\
1311	1.8768115137529e-13\\
1312	1.75876487167007e-13\\
1313	1.6408942532753e-13\\
1314	1.52351306931409e-13\\
1315	1.40703148346088e-13\\
1316	1.2919185015661e-13\\
1317	1.17866923368955e-13\\
1318	1.06780672577978e-13\\
1319	9.59865855028966e-14\\
1320	8.55374617193114e-14\\
1321	7.54857191532503e-14\\
1322	6.58819753044644e-14\\
1323	5.67750229967858e-14\\
1324	4.82113689539356e-14\\
1325	4.02347735116839e-14\\
1326	3.2885915422512e-14\\
1327	2.6202369243248e-14\\
1328	2.021818777346e-14\\
1329	1.49638455885519e-14\\
1330	1.04660518731286e-14\\
1331	6.74669655328827e-15\\
1332	3.79701519820396e-15\\
1333	1.5049840313108e-15\\
1334	1.27748513780431e-15\\
};
\addlegendentry{$E=16$}
\end{axis}
\end{tikzpicture}%

%% file: residual-v1-M=100.tikz
%
%
%
\begin{tikzpicture}

\begin{axis}[%
width=0.375\textwidth,
scale only axis,
xmin=1,
xmax=10,
xlabel={Iteration},
ymode=log,
ymin=1e-14,
ymax=1e2,
yminorticks=true,
ylabel={Residual},
title style={font=\bfseries},
legend style={legend cell align=left, align=left, nodes={scale=0.6, transform shape}}
]
\addplot +[mark options={scale=0.6}]
  table[row sep=crcr]{%
1	29.445479596844\\
2	0.0240450100772832\\
3	0.00496301353907556\\
4	0.00337440722739875\\
5	0.00321462807368931\\
6	0.00311976088228501\\
7	0.00305575802590584\\
8	0.00301167740327441\\
9	0.00298097188341408\\
10	0.00295942261486531\\
};
\addlegendentry{$r=12$};

\addplot +[mark options={scale=0.6}]
  table[row sep=crcr]{%
1	33.3201810585994\\
2	0.000103572988251927\\
3	1.98851276330661e-05\\
4	2.02605545310153e-05\\
5	2.00839712592946e-05\\
6	1.99922432123682e-05\\
7	1.99338728510628e-05\\
8	1.98935665806581e-05\\
9	1.98651731786506e-05\\
10	1.98450847545271e-05\\
};
\addlegendentry{$r=18$};

\addplot +[mark options={scale=0.6}]
  table[row sep=crcr]{%
1	23.6914944857794\\
2	6.03593433717755e-06\\
3	1.93504610817013e-09\\
4	1.38670665478639e-09\\
5	1.14869956180808e-09\\
6	1.01925023588584e-09\\
7	9.43854153630493e-10\\
8	8.97964937307655e-10\\
9	8.69068974782183e-10\\
10	8.50343285084767e-10\\
};
\addlegendentry{$r=24$};

\addplot +[mark options={scale=0.6}]
  table[row sep=crcr]{%
1	22.5136651221703\\
2	1.33154557691813e-07\\
3	9.4680765848426e-13\\
4	1.12541255409108e-12\\
5	9.86656470251004e-13\\
6	1.03035298907108e-12\\
7	9.39194589826432e-13\\
8	9.72567366439629e-13\\
9	1.03395838869426e-12\\
10	9.37717772042641e-13\\
};
\addlegendentry{$r=30$};

\end{axis}
\end{tikzpicture}%

%% file: runtime-v1-3d-M=100.tikz
%
%
%
\begin{tikzpicture}

\begin{axis}[%
width=0.375\textwidth,
scale only axis,
xmin=8,
xmax=32,
xmode=log,
log basis x=2,
xlabel={Rank},
ymin=1e-1,
ymax=1e3,
ymode=log,
ylabel={Runtime (sec)},
axis background/.style={fill=white},
title style={font=\bfseries},
legend style={legend cell align=left, align=left, at={(0.02,0.98)},anchor=north west, nodes={scale=0.6, transform shape}}
]
\addplot +[mark options={scale=0.6}]
  table[row sep=crcr]{%
12	2.356\\
15	4.537\\
18	8.32\\
21	12.919\\
24	18.638\\
27	26.383\\
};
\addlegendentry{Construction operators; $\Oh{r^3}$}

\addplot +[mark options={scale=0.6}]
  table[row sep=crcr]{%
12	20.012\\
15	45.824\\
18	103.05\\
21	186.835\\
24	361.163\\
27	650.888\\
};
\addlegendentry{Solve for $\mv{U}_1$, $\mv{U}_2$ and $\mv{U}_3$; $\Oh{r^4}$}

\addplot +[mark options={scale=0.6}]
  table[row sep=crcr]{%
12	22.852\\
15	50.692\\
18	112.172\\
21	200.158\\
24	380.345\\
27	677.734\\
};
\addlegendentry{Total runtime; $\Oh{r^4}$}

\end{axis}
\end{tikzpicture}%

%% file: residual-v2-M=100.tikz
%
%
%
\begin{tikzpicture}

\begin{axis}[%
width=0.375\textwidth,
scale only axis,
xmin=1,
xmax=10,
xlabel={Iteration},
ymode=log,
ymin=1e-14,
ymax=1e2,
yminorticks=true,
ylabel={Residual},
title style={font=\bfseries},
legend style={legend cell align=left, align=left, nodes={scale=0.6, transform shape}}
]
\addplot +[mark options={scale=0.6}]
  table[row sep=crcr]{%
1	64.4912923658883\\
2	0.0830348898452449\\
3	0.0411550216650904\\
4	0.0202191466373435\\
5	0.0253348392095572\\
6	0.0179349404568083\\
7	0.0164025967636903\\
8	0.0149357709850446\\
9	0.0140429962127991\\
10	0.0140958212160474\\
};
\addlegendentry{$r=12$};

\addplot +[mark options={scale=0.6}]
  table[row sep=crcr]{%
1	54.9396033507331\\
2	2.49246597428769e-05\\
3	9.00633105631812e-06\\
4	6.82145551946668e-06\\
5	6.60351272456749e-06\\
6	6.22454316266196e-06\\
7	6.19582399344515e-06\\
8	6.75337309170577e-06\\
9	7.16751724806227e-06\\
10	7.43912284387568e-06\\
};
\addlegendentry{$r=18$};

\addplot +[mark options={scale=0.6}]
  table[row sep=crcr]{%
1	96.6176733256371\\
2	9.62192311549907e-08\\
3	4.29550282614718e-09\\
4	1.95570289371032e-09\\
5	2.02460664247373e-09\\
6	1.74434260136437e-09\\
7	1.66619966734321e-09\\
8	1.23986304272463e-09\\
9	1.480402686198e-09\\
10	1.14034624988237e-09\\
};
\addlegendentry{$r=24$};

\addplot +[mark options={scale=0.6}]
  table[row sep=crcr]{%
1	64.4723874726743\\
2	7.16113810731306e-09\\
3	4.56901730685202e-10\\
4	3.45033016344685e-10\\
5	3.32314370749094e-10\\
6	2.95208322055232e-10\\
7	1.70803128089897e-10\\
8	2.63488846050608e-10\\
9	1.269832299304e-10\\
10	1.17888840513664e-10\\
};
\addlegendentry{$r=30$};

\end{axis}
\end{tikzpicture}%

%% file: runtime-v2-3d-M=100.tikz
%
%
%
\begin{tikzpicture}

\begin{axis}[%
width=0.375\textwidth,
scale only axis,
xmin=8,
xmax=32,
xmode=log,
log basis x=2,
xlabel={Rank},
ymin=1e-1,
ymax=1e3,
ymode=log,
ylabel={Runtime (sec)},
axis background/.style={fill=white},
title style={font=\bfseries},
legend style={legend cell align=left, align=left, at={(0.02,0.98)},anchor=north west, nodes={scale=0.6, transform shape}}
]
\addplot +[mark options={scale=0.6}]
  table[row sep=crcr]{%
12	0.19\\
15	0.29\\
18	0.46\\
21	0.61\\
24	0.86\\
27	1.08\\
};
\addlegendentry{Construction operators $\mv{U}_1$, $\mv{U}_2$, $\mv{U}_3$; $\Oh{r^2}$}

\addplot +[mark options={scale=0.6}]
  table[row sep=crcr]{%
18	0.04\\
21	0.07\\
24	0.12\\
27	0.17\\
};
\addlegendentry{Construction operator $\mt{G}$}

\addplot +[mark options={scale=0.6}]
  table[row sep=crcr]{%
12	0.305\\
15	0.482\\
18	0.636\\
21	0.824\\
24	1.088\\
27	1.338\\
};
\addlegendentry{Solve for $\mv{U}_1$, $\mv{U}_2$ and $\mv{U}_3$; $\Oh{r^2}$}

\addplot +[mark options={scale=0.6}]
  table[row sep=crcr]{%
12	0.275\\
15	1.597\\
18	6.644\\
21	26.132\\
24	86.905\\
27	362.421\\
};
\addlegendentry{Solve for $\mt{G}$; $\Oh{r^9}$}

\addplot
  table[row sep=crcr]{%
12	1.086\\
15	2.738\\
18	8.137\\
21	28.040\\
24	89.474\\
27	365.918\\
};
\addlegendentry{Total runtime; $\Oh{r^7} \rightarrow \Oh{r^9}$}

\end{axis}
\end{tikzpicture}%

%% file: residual-v3-M=100.tikz
%
%
%
\begin{tikzpicture}

\begin{axis}[%
width=0.375\textwidth,
scale only axis,
xmin=1,
xmax=10,
xlabel={Iteration},
ymode=log,
ymin=1e-14,
ymax=1e2,
yminorticks=true,
ylabel={Residual},
title style={font=\bfseries},
legend style={legend cell align=left, align=left, nodes={scale=0.6, transform shape}}
]
\addplot +[mark options={scale=0.6}]
  table[row sep=crcr]{%
1	71.0846101809674\\
2	0.112866699983702\\
3	0.003514979158243\\
4	0.00580241919853146\\
5	0.00615456555490827\\
6	0.0054831096111438\\
7	0.00501576337896503\\
8	0.00476062000695987\\
9	0.00462370536577146\\
10	0.0045496203495438\\
};
\addlegendentry{$r=12$};

\addplot +[mark options={scale=0.6}]
  table[row sep=crcr]{%
1	57.9416399576333\\
2	0.00870270445954991\\
3	1.12824284017751e-05\\
4	1.04501238009951e-05\\
5	8.92560457810916e-06\\
6	6.7868998795636e-06\\
7	5.31437408121152e-06\\
8	4.50272142897365e-06\\
9	4.0770229341764e-06\\
10	3.85348439876562e-06\\
};
\addlegendentry{$r=18$};

\addplot +[mark options={scale=0.6}]
  table[row sep=crcr]{%
1	99.1326137164746\\
2	0.00124140383410696\\
3	5.61547382579972e-10\\
4	1.21133829993766e-09\\
5	1.79458851033086e-09\\
6	1.68366835680527e-09\\
7	1.56977803050825e-09\\
8	1.50728756314024e-09\\
9	1.47422096371595e-09\\
10	1.45643630085713e-09\\
};
\addlegendentry{$r=24$};

\addplot +[mark options={scale=0.6}]
  table[row sep=crcr]{%
1	56.5314106555572\\
2	2.90607724170813e-05\\
3	1.26447632772355e-12\\
4	1.23703200356548e-12\\
5	1.07376424293674e-12\\
6	1.10878792567679e-12\\
7	1.40671664593591e-12\\
8	1.16684543246428e-12\\
9	1.14491415476001e-12\\
10	9.9687514520859e-13\\
};
\addlegendentry{$r=30$};

\end{axis}
\end{tikzpicture}%

%% file: runtime-v3-3d-M=100.tikz
%
%
%
\begin{tikzpicture}

\begin{axis}[%
width=0.375\textwidth,
scale only axis,
xmin=8,
xmax=32,
xmode=log,
log basis x=2,
xlabel={Rank},
ymin=1e-1,
ymax=1e3,
ymode=log,
ylabel={Runtime (sec)},
axis background/.style={fill=white},
title style={font=\bfseries},
legend style={legend cell align=left, align=left, at={(0.02,0.98)},anchor=north west, nodes={scale=0.6, transform shape}}
]
\addplot +[mark options={scale=0.6}]
  table[row sep=crcr]{%
12	0.15\\
15	0.24\\
18	0.35\\
21	0.47\\
24	0.72\\
27	0.91\\
};
\addlegendentry{Construction operators $\mv{U}_1$, $\mv{U}_2$; $\Oh{r^2}$}

\addplot +[mark options={scale=0.6}]
  table[row sep=crcr]{%
12	0.74\\
15	1.48\\
18	2.54\\
21	4.01\\
24	6.10\\
27	9.08\\
};
\addlegendentry{Construction operators $\conj{\mv{U}_3}\mv{G}_{(3)}$; $\Oh{r^3}$}

\addplot +[mark options={scale=0.6}]
  table[row sep=crcr]{%
12	0.204\\
15	0.361\\
18	0.433\\
21	0.585\\
24	0.909\\
27	1.08\\
};
\addlegendentry{Solve for $\mv{U}_1$, $\mv{U}_2$; $\Oh{r^2}$}

\addplot +[mark options={scale=0.6}]
  table[row sep=crcr]{%
12	6.477\\
15	16.159\\
18	32.529\\
21	61.642\\
24	124.067\\
27	223.838\\
};
\addlegendentry{Solve for $\conj{\mv{U}_3}\mv{G}_{(3)}$; $\Oh{r^4}$}

\addplot +[mark options={scale=0.6}]
  table[row sep=crcr]{%
12	7.895\\
15	18.672\\
18	36.266\\
21	67.145\\
24	132.340\\
27	235.499\\
};
\addlegendentry{Total runtime; $\Oh{r^4}$}

\end{axis}
\end{tikzpicture}%

%% file: runtimes-v123-M=100-iter=4.tikz
%
%
%
\begin{tikzpicture}

\begin{axis}[%
width=0.5\textwidth,
at={(0.758in,0.481in)},
scale only axis,
xmin=8,
xmax=32,
xmode=log,
log basis x=2,
xlabel={Rank},
ymin=1e-1,
ymax=1e3,
ymode=log,
axis background/.style={fill=white},
title style={font=\bfseries},
legend style={legend cell align=left, align=left, at={(0.02,0.98)},anchor=north west}
]
\addplot
  table[row sep=crcr]{%
2	0.121032\\
4	0.298357\\
6	0.804427\\
8	2.159353\\
10	4.607214\\
12	8.528169\\
14	15.114430\\
16	26.767309\\
18	43.155919\\
20	67.556569\\
24	160.918632\\
28	340.669794\\
};
\addlegendentry{version 1};

\addplot
  table[row sep=crcr]{%
2	0.113979\\
4	0.137980\\
6	0.165415\\
8	0.219098\\
10	0.347787\\
12	0.579338\\
14	1.308169\\
16	2.298700\\
18	7.646836\\ 
20	16.559335\\
24	66.714624\\ 
};
\addlegendentry{version 2};

\addplot
  table[row sep=crcr]{%
2	0.112859\\
4	0.173383\\
6	0.428806\\
8	0.863533\\
10	1.709260\\
12	3.179376\\
14	5.427046\\
16	8.925007\\
18	14.812203\\
20	23.201373\\
24	58.669655\\
28	113.110773\\
};
\addlegendentry{version 3};

\end{axis}
\end{tikzpicture}%

%% file: error-variable-2d-r=17.tikz
%
%
%
\begin{tikzpicture}

\begin{axis}[%
width=0.375\textwidth,
scale only axis,
xmin=0,
xmax=10,
xlabel={Iteration},
ymode=log,
ymin=1e-08,
ymax=1e4,
yminorticks=true,
ylabel={Error/residual},
title style={font=\bfseries},
legend style={legend cell align=left, align=left, nodes={scale=0.6, transform shape}}
]
\addplot +[mark options={scale=0.6}]
  table[row sep=crcr]{%
0	91.7697049916572\\
1	69.1182293879963\\
2	0.0816986906568348\\
3	0.000341646810449371\\
4	3.27856246136877e-05\\
5	2.9846883171673e-05\\
6	8.77843079404105e-06\\
7	2.91294113848747e-06\\
8	2.09596106252717e-06\\
9	1.91681540678325e-06\\
10	1.86651261228684e-06\\
};
\addlegendentry{Error w.r.t. semi-exact solution}

\addplot +[mark options={scale=0.6}]
  table[row sep=crcr]{%
1	1642.1469962723\\
2	0.539133340498152\\
3	0.000811160365871902\\
4	0.000164155462217153\\
5	0.00014568976518238\\
6	0.000101700059645779\\
7	0.000111104315552941\\
8	0.000120132698033311\\
9	0.00012372694496398\\
10	0.000124958850838427\\
};
\addlegendentry{Residual}

\addplot +[mark options={scale=0.6}]
  table[row sep=crcr]{%
0	9.81801715216901e-07\\
1	9.81801715216901e-07\\
2	9.81801715216901e-07\\
3	9.81801715216901e-07\\
4	9.81801715216901e-07\\
5	9.81801715216901e-07\\
6	9.81801715216901e-07\\
7	9.81801715216901e-07\\
8	9.81801715216901e-07\\
9	9.81801715216901e-07\\
10	9.81801715216901e-07\\
};
\addlegendentry{Expected min. error}

\end{axis}
\end{tikzpicture}%

%% file: singularvalues-variable-2d.tikz
%
%
%
\begin{tikzpicture}

\begin{axis}[%
width=0.375\textwidth,
scale only axis,
xmin=0,
xmax=18,
xlabel={Rank},
ymode=log,
ymin=1e-15,
ymax=100000,
yminorticks=true,
ylabel={Singular value},
title style={font=\bfseries},
legend style={legend cell align=left, align=left, nodes={scale=0.6, transform shape}}
]
\addplot +[line width=0.75pt, mark options={scale=0.75}]
  table[row sep=crcr]{%
1	91.5604707115939\\
2	74.6389384752908\\
3	29.2972299146646\\
4	6.84195765401993\\
5	1.43372131933781\\
6	0.313646603222012\\
7	0.0721465200993463\\
8	0.0182599098236162\\
9	0.00579283478411671\\
10	0.00173820610043303\\
11	0.000459145695822899\\
12	0.000162865463588722\\
13	6.17008711807617e-05\\
14	1.58088958602978e-05\\
15	7.9372264391688e-06\\
16	4.06171839552595e-06\\
17	9.81801715216901e-07\\
};
\addlegendentry{Semi-exact (solution)}

\addplot +[mark options={scale=0.6}]
  table[row sep=crcr]{%
1	104.94499873105\\
2	73.8971898422717\\
3	23.4630834885581\\
4	3.11580754233616\\
5	0.0757433921198959\\
6	0.00431922168464984\\
7	0.000151172821897063\\
8	8.24146473218871e-05\\
9	4.78131342640318e-06\\
10	5.01150486952746e-07\\
11	4.96206089718419e-08\\
12	8.9779514066572e-09\\
13	2.38187283176429e-09\\
14	1.82957041168799e-10\\
15	4.09996682756953e-11\\
16	9.76568677305862e-12\\
17	6.86218867460603e-13\\
};
\addlegendentry{Iter = 1}

\addplot +[mark options={scale=0.6}]
  table[row sep=crcr]{%
1	91.5557053322345\\
2	74.6401038501054\\
3	29.2929248019765\\
4	6.84231965717022\\
5	1.433524367262\\
6	0.314015961178774\\
7	0.0725501783204869\\
8	0.0183042652295159\\
9	0.00621869334954903\\
10	0.00252608257553629\\
11	0.000325257685876546\\
12	6.24421270997746e-05\\
13	6.51789015257302e-06\\
14	8.10760761207674e-07\\
15	1.04160113551023e-07\\
16	1.11413793910024e-08\\
17	8.43697585119063e-10\\
};
\addlegendentry{Iter = 2}

\addplot +[mark options={scale=0.6}]
  table[row sep=crcr]{%
1	91.5604143551281\\
2	74.6389234175081\\
3	29.2972355490188\\
4	6.84194541468434\\
5	1.43372076389966\\
6	0.313654363068101\\
7	0.072137691176779\\
8	0.0182664588040593\\
9	0.00578745104107383\\
10	0.00174819351050666\\
11	0.000455281062643901\\
12	0.000183704695552869\\
13	6.02615489837228e-05\\
14	8.37991280642365e-06\\
15	1.28531119188795e-06\\
16	1.47385622614736e-07\\
17	1.08619238131996e-08\\
};
\addlegendentry{Iter = 3}

\addplot +[mark options={scale=0.6}]
  table[row sep=crcr]{%
1	91.5604609689959\\
2	74.6389360096314\\
3	29.2972258567257\\
4	6.84195649594181\\
5	1.4337200144241\\
6	0.31364741466218\\
7	0.0721468665648639\\
8	0.0182590439832759\\
9	0.00579337281806888\\
10	0.00173782411066976\\
11	0.000460130892135346\\
12	0.000162326005827751\\
13	6.23917083632737e-05\\
14	1.90772110584587e-05\\
15	5.17285298489361e-06\\
16	7.43452535249698e-07\\
17	7.32586425548494e-08\\
};
\addlegendentry{Iter = 4}

\addplot +[mark options={scale=0.6}]
  table[row sep=crcr]{%
1	91.5604648157331\\
2	74.6389385943661\\
3	29.297231078388\\
4	6.841956131037\\
5	1.4337223917996\\
6	0.313647457923915\\
7	0.0721449021353797\\
8	0.0182614574650077\\
9	0.00579136701816684\\
10	0.00173969146367384\\
11	0.000457659607364912\\
12	0.000165060814096801\\
13	6.12433455737195e-05\\
14	1.69528672683657e-05\\
15	8.01831675389378e-06\\
16	2.12455667815171e-06\\
17	2.68836907924902e-07\\
};
\addlegendentry{Iter = 5}

\end{axis}
\end{tikzpicture}%

%% file: error-variable-2d-M=1000.tikz
%
%
%
\begin{tikzpicture}

\begin{axis}[%
width=0.375\textwidth,
scale only axis,
xmin=0,
xmax=10,
xlabel={Iteration},
ymode=log,
ymin=1e-12,
ymax=1e4,
yminorticks=true,
ylabel={Error w.r.t. semi-exact solution},
axis background/.style={fill=white},
title style={font=\bfseries},
legend style={legend cell align=left, align=left, nodes={scale=0.6, transform shape}}
]
\addplot +[mark options={scale=0.6}]
  table[row sep=crcr]{%
0	91.769399089003\\
1	142.503179767938\\
2	0.67234654692724\\
3	0.00606996786823398\\
4	0.00110015438609123\\
5	0.000494055822559004\\
6	0.000370998764989986\\
7	0.00034997374113609\\
8	0.000344640042499494\\
9	0.000343198839089851\\
10	0.000342804323823468\\
};
\addlegendentry{$r=12$}

\addplot +[mark options={scale=0.6}]
  table[row sep=crcr]{%
0	91.7626704966655\\
1	63.3406993105331\\
2	0.0680667938870474\\
3	0.000200040124516136\\
4	3.34708911096567e-05\\
5	5.3375032758725e-06\\
6	1.25086485891803e-06\\
7	7.35526654745818e-07\\
8	5.64209089638106e-07\\
9	4.92470010913578e-07\\
10	4.60953575977283e-07\\
};
\addlegendentry{$r=18$}

\addplot +[mark options={scale=0.6}]
  table[row sep=crcr]{%
0	91.7745433240276\\
1	54.1085922527951\\
2	0.00912561078654957\\
3	1.29681039667139e-06\\
4	7.23461826763191e-09\\
5	1.7315470006327e-09\\
6	1.04102879910399e-09\\
7	7.81927039244e-10\\
8	6.74522434568001e-10\\
9	6.29062686360479e-10\\
10	6.0929936185099e-10\\
};
\addlegendentry{$r=24$}

\addplot +[mark options={scale=0.6}]
  table[row sep=crcr]{%
0	91.7668599890111\\
1	48.1304422901816\\
2	0.000262319545765502\\
3	7.29049809534655e-11\\
4	4.22548700996891e-11\\
5	8.29302068208756e-11\\
6	4.63198345464395e-11\\
7	5.04428102678417e-11\\
8	7.52695698489469e-11\\
9	6.48048965295417e-11\\
10	4.61200974187957e-11\\
};
\addlegendentry{$r=36$}

\end{axis}
\end{tikzpicture}%

%% file: residuals-variable-2d-M=1000.tikz
%
%
%
\begin{tikzpicture}

\begin{axis}[%
width=0.375\textwidth,
scale only axis,
xmin=0,
xmax=10,
xlabel={Iteration},
ymode=log,
ymin=1e-12,
ymax=1e4,
yminorticks=true,
ylabel={Residual},
title style={font=\bfseries},
legend style={legend cell align=left, align=left, nodes={scale=0.6, transform shape}}
]
\addplot +[mark options={scale=0.6}]
  table[row sep=crcr]{%
1	1805.83506487425\\
2	2.48278261819334\\
3	0.0139447558894335\\
4	0.00421705054587041\\
5	0.00307327601499764\\
6	0.00343119893046187\\
7	0.00357089281969172\\
8	0.00361177555784844\\
9	0.00362320015195076\\
10	0.00362635694205762\\
};
\addlegendentry{$r=12$}

\addplot +[mark options={scale=0.6}]
  table[row sep=crcr]{%
1	1687.83619529764\\
2	0.301006053951278\\
3	0.000902615617343562\\
4	0.000164670139951428\\
5	4.68939679483661e-05\\
6	4.46997553370085e-05\\
7	5.24591439795395e-05\\
8	5.76195011782948e-05\\
9	6.03421245283033e-05\\
10	6.16645280914799e-05\\
};
\addlegendentry{$r=18$}

\addplot +[mark options={scale=0.6}]
  table[row sep=crcr]{%
1	1312.48912081465\\
2	0.187418405960824\\
3	1.18286004812246e-05\\
4	1.91892990786301e-07\\
5	1.38926048331364e-07\\
6	1.69006544004044e-07\\
7	1.94980521117327e-07\\
8	2.09902437114329e-07\\
9	2.17356189790476e-07\\
10	2.20856712966105e-07\\
};
\addlegendentry{$r=24$}

\addplot +[mark options={scale=0.6}]
  table[row sep=crcr]{%
1	1157.50444452196\\
2	0.013967750671266\\
3	4.98483703646609e-10\\
4	5.07848113990708e-10\\
5	5.46770217610707e-10\\
6	5.10405398498919e-10\\
7	5.37242731799656e-10\\
8	5.06250776075895e-10\\
9	5.0551546363387e-10\\
10	5.11993842674513e-10\\
};
\addlegendentry{$r=36$}

\end{axis}
\end{tikzpicture}%

%% file: runtime-variable-2d-M=1000.tikz
%
%
%
\begin{tikzpicture}

\begin{axis}[%
width=0.375\textwidth,
at={(0.758in,0.481in)},
scale only axis,
xmin=0,
xmax=50,
xlabel={Rank (r)},
ymin=0,
ymax=120,
ylabel={Runtime (sec)},
axis background/.style={fill=white},
title style={font=\bfseries},
legend style={legend cell align=left, align=left, draw=white!15!black, nodes={scale=0.6, transform shape}}
]
\addplot +[mark options={scale=0.6}]
  table[row sep=crcr]{%
1	0.0888668\\
3	0.4259461\\
6	1.5748525\\
9	3.8105173\\
12	6.4289317\\
18	13.9983622\\
24	25.6150505\\
36	57.870017\\
48	104.849505\\
};
\addlegendentry{Update U}

\addplot +[mark options={scale=0.6}]
  table[row sep=crcr]{%
1	0.0740352\\
3	0.1838032\\
6	0.6372322\\
9	1.6990159\\
12	2.8384972\\
18	5.9882744\\
24	10.3674736\\
36	22.6386341\\
48	39.637221\\
};
\addlegendentry{Update V}

\addplot +[mark options={scale=0.6}]
  table[row sep=crcr]{%
1	0.162902\\
3	0.6097493\\
6	2.2120847\\
9	5.5095332\\
12	9.2674288\\
18	19.9866366\\
24	35.9825241\\
36	80.5086511\\
48	144.486726\\
};
\addlegendentry{Runtime}

\end{axis}
\end{tikzpicture}%

%% file: runtime-variable-2d-M=1000+loglog.tikz
%
%
%
\begin{tikzpicture}

\begin{axis}[%
width=0.375\textwidth,
at={(0.758in,0.481in)},
scale only axis,
xmin=1e0,
xmax=1e2,
xmode=log,
xlabel={Rank (r)},
ymin=1e-2,
ymax=1e3,
ymode=log,
ylabel={Runtime (sec)},
axis background/.style={fill=white},
title style={font=\bfseries},
legend style={legend cell align=left, align=left, draw=white!15!black, nodes={scale=0.6, transform shape}}
]
\addplot +[mark options={scale=0.6}]
  table[row sep=crcr]{%
1	0.0888668\\
3	0.4259461\\
6	1.5748525\\
9	3.8105173\\
12	6.4289317\\
18	13.9983622\\
24	25.6150505\\
36	57.870017\\
48	104.849505\\
};
\addlegendentry{Update U}

\addplot +[mark options={scale=0.6}]
  table[row sep=crcr]{%
1	0.0740352\\
3	0.1838032\\
6	0.6372322\\
9	1.6990159\\
12	2.8384972\\
18	5.9882744\\
24	10.3674736\\
36	22.6386341\\
48	39.637221\\
};
\addlegendentry{Update V}

\addplot +[mark options={scale=0.6}]
  table[row sep=crcr]{%
1	0.162902\\
3	0.6097493\\
6	2.2120847\\
9	5.5095332\\
12	9.2674288\\
18	19.9866366\\
24	35.9825241\\
36	80.5086511\\
48	144.486726\\
};
\addlegendentry{Runtime}

\end{axis}
\end{tikzpicture}%

%% file: cprank-orderDvr=7-wavenumber.tikz
%
%
%
\begin{tikzpicture}

\begin{axis}[%
width=0.8\textwidth,
scale only axis,
xmode=log,
xmin=1,
xmax=100,
xminorticks=true,
xlabel={CP-rank wave number ($s$)},
ymode=log,
ymin=1e-05,
ymax=10,
yminorticks=true,
]
\addplot +[mark options={scale=0.6}]
  table[row sep=crcr]{%
2	1.05227345286093\\
3	0.244146486443476\\
4	0.133671324696763\\
5	0.078697749639831\\
6	0.044359489377007\\
7	0.031641486582529\\
8	0.037281118362448\\
9	0.014015530708451\\
10	0.008830546281349\\
12	0.004624701678211\\
14	0.002818499729945\\
16	0.002016876121034\\
18	0.001416447606645\\
20	0.000535917640435368\\
24	0.000323434951709611\\
28	0.000353898905469884\\
32	0.000140229106573297\\
36	0.00017608556280857\\
40	8.38711719413076e-05\\
50	0.000228187681150959\\
60	2.15329643311612e-05\\
70	3.93978992820898e-05\\
80	1.81745231477847e-05\\
90	1.69387324393532e-05\\
100	2.73835693520574e-05\\
125	1.24834798862494e-05\\
150	3.12113923124357e-05\\
175	1.57655860029137e-05\\
200	6.16405807729592e-06\\
};
\end{axis}
\end{tikzpicture}%

%% file: residual-all-orderDvr=7-solrank=16-waverank=32.tikz
%
%
%
\begin{tikzpicture}

\begin{axis}[%
width=0.8\textwidth,
scale only axis,
xmin=0,
xmax=10,
xlabel={Iteration},
ymode=log,
ymin=0.0001,
ymax=10,
yminorticks=true,
ylabel={Residual},
legend style={legend cell align=left, align=left, nodes={scale=0.6, transform shape}}
]
\addplot +[mark options={scale=0.6}]
  table[row sep=crcr]{%
1	0.630212863576271\\
2	0.00236939652617741\\
3	0.000183696106291433\\
4	0.000205125867192498\\
5	0.000217073908342878\\
6	0.000232252713804881\\
7	0.000240434684979846\\
8	0.000246202738610989\\
9	0.000250613309670309\\
10	0.000253840851242934\\
};
\addlegendentry{version 1}

\addplot +[mark options={scale=0.6}]
  table[row sep=crcr]{%
1	5.76193390263778\\
2	0.00199165947589355\\
3	0.00480116735319347\\
4	0.0074527179610751\\
5	0.00497272483988253\\
6	0.00349635641595582\\
7	0.00353391600616445\\
8	0.00448447435444537\\
9	0.00391651289897862\\
10	0.00449446446658272\\
};
\addlegendentry{version 2}

\addplot +[mark options={scale=0.6}]
  table[row sep=crcr]{%
1	5.99885915687037\\
2	0.00162912852608232\\
3	0.000156412253631898\\
4	0.000234316940658095\\
5	0.000238189211380308\\
6	0.000246058475960696\\
7	0.000248372833675671\\
8	0.0002494817995274\\
9	0.000249323193314131\\
10	0.000249377753508427\\
};
\addlegendentry{version 3}

\end{axis}
\end{tikzpicture}%

%% file: residual-all-orderDvr=7-waverank=32.tikz
%
%
%
\begin{tikzpicture}

\begin{axis}[%
width=0.8\textwidth,
scale only axis,
xmin=0,
xmax=24,
xlabel={Maximal low-rank $r$ of approximation},
ymode=log,
ymin=1e-07,
ymax=10,
yminorticks=true,
ylabel={Residual},
legend style={legend cell align=left, align=left, nodes={scale=0.6, transform shape}}
]
\addplot +[mark options={scale=0.6}]
  table[row sep=crcr]{%
2	0.416579533795744\\
3	0.422447282349009\\
4	0.344852710759694\\
5	0.216751409956652\\
6	0.193475913766274\\
8	0.115671986660269\\
10	0.0441083721944136\\
12	0.00625681163438005\\
16	0.000254168229684401\\
20	1.96606765920935e-05\\
24	3.38112405230756e-06\\
};
\addlegendentry{version 1}

\addplot +[mark options={scale=0.6}]
  table[row sep=crcr]{%
2	0.408104357851707\\
3	0.416593716596844\\
4	0.381425321591709\\
5	0.355503942249896\\
6	0.612902907009588\\
8	0.172461327072626\\
10	0.317638418643652\\
12	0.058630324218368\\
16	0.00253457993217098\\
20	0.000854124485647378\\
24	2.41804360041505e-05\\
};
\addlegendentry{version 2}

\addplot +[mark options={scale=0.6}]
  table[row sep=crcr]{%
2	0.262723216096621\\
3	0.33735200592287\\
4	0.235274682627076\\
5	0.208621518498484\\
6	0.215974441770047\\
8	0.139189188437221\\
10	0.0506823170626152\\
12	0.00871632332524819\\
16	0.000249359149380671\\
20	6.09567793406564e-06\\
24	6.19513061823379e-07\\
};
\addlegendentry{version 3}

\end{axis}
\end{tikzpicture}%

%% file: residual-all-orderDvr=14-waverank=32.tikz
%
%
%
\begin{tikzpicture}

\begin{axis}[%
width=0.8\textwidth,
scale only axis,
xmin=0,
xmax=24,
xlabel={Maximal low-rank $r$ of approximation},
ymode=log,
ymin=1e-07,
ymax=10,
yminorticks=true,
ylabel={Residual},
legend style={legend cell align=left, align=left, nodes={scale=0.6, transform shape}}
]
\addplot +[mark options={scale=0.6}]
  table[row sep=crcr]{%
2	0.529736425863644\\
3	0.41185292265979\\
4	0.797751005461586\\
5	0.383963910501724\\
6	0.235218216916467\\
8	0.177712693954026\\
10	0.192198532968015\\
12	0.107536267749343\\
16	0.00661420701429977\\
20	0.000375658790639547\\
};
\addlegendentry{version 1}

\addplot +[mark options={scale=0.6}]
  table[row sep=crcr]{%
2	0.419350644231988\\
3	0.418288955229564\\
4	0.619579465639644\\
5	1.30874082625856\\
6	0.43753775897463\\
8	1.92229404964787\\
10	0.125836027945509\\
12	0.0818702332994378\\
16	0.0497005473708251\\
20	0.0262923273291113\\
};
\addlegendentry{version 2}

\addplot +[mark options={scale=0.6}]
  table[row sep=crcr]{%
2	0.40992504816427\\
3	0.298816882560173\\
4	0.273374874299591\\
5	0.254216237458379\\
6	0.245731322259959\\
8	0.183896368606948\\
10	0.162804215189307\\
12	0.0683209417430981\\
16	0.0118254838121869\\
20	0.000301607605454202\\
};
\addlegendentry{version 3}

\end{axis}
\end{tikzpicture}%

%% file: runtimes-all-orderDvr=7-waverank=32+loglog.tikz
%
%
%
\begin{tikzpicture}

\begin{axis}[%
width=0.8\textwidth,
scale only axis,
xmode=log,
xmin=2,
xmax=32,
log basis x = 2,
xminorticks=true,
xlabel={Maximal low-rank $r$ of approximation},
ymode=log,
ymin=0.1,
ymax=100000,
yminorticks=true,
ylabel={Runtime (sec)},
legend style={legend cell align=left, align=left, nodes={scale=0.6, transform shape}}
]
\addplot +[mark options={scale=0.6}]
  table[row sep=crcr]{%
2	0.3874572\\
3	0.7061261\\
4	1.7231405\\
5	4.4983776\\
6	10.4394877\\
8	30.6514726\\
10	74.8277133\\
12	162.0073253\\
16	535.7783054\\
20	1355.4675901\\
24	2947.3223489\\
};
\addlegendentry{version 1}

\addplot +[mark options={scale=0.6}]
  table[row sep=crcr]{%
2	0.3251689\\
3	0.249121\\
4	0.2236293\\
5	0.3117947\\
6	0.374879\\
8	0.8782378\\
10	1.6516168\\
12	3.4270538\\
16	14.3528478\\
20	61.2257691\\
24	244.1689815\\
};
\addlegendentry{version 2}

\addplot +[mark options={scale=0.6}]
  table[row sep=crcr]{%
2	0.3097978\\
3	0.5053149\\
4	0.7786363\\
5	1.7957841\\
6	3.6903191\\
8	10.5277331\\
10	31.1629924\\
12	54.2125362\\
16	180.8078978\\
20	455.2348111\\
24	994.7267573\\
};
\addlegendentry{version 3}

\end{axis}
\end{tikzpicture}%

%% file: runtimes-all-orderDvr=14-waverank=32+loglog.tikz
%
%
%
\begin{tikzpicture}

\begin{axis}[%
width=0.8\textwidth,
scale only axis,
xmode=log,
xmin=2,
xmax=32,
log basis x = 2,
xminorticks=true,
xlabel={Maximal low-rank $r$ of approximation},
ymode=log,
ymin=0.1,
ymax=100000,
yminorticks=true,
ylabel={Runtime (sec)},
legend style={legend cell align=left, align=left, nodes={scale=0.6, transform shape}}
]
\addplot +[mark options={scale=0.6}]
  table[row sep=crcr]{%
2	0.4081247\\
3	1.4887474\\
4	3.6488767\\
5	9.2777145\\
6	20.4467515\\
8	74.4404381\\
10	186.7417244\\
12	405.5468823\\
16	1409.5032916\\
20	3822.2569589\\
};
\addlegendentry{version 1}

\addplot +[mark options={scale=0.6}]
  table[row sep=crcr]{%
2	0.6022992\\
3	0.316849\\
4	0.404419\\
5	0.6588588\\
6	0.7896445\\
8	1.3425819\\
10	2.5071714\\
12	4.3659661\\
16	16.7118949\\
20	64.7096582\\
};
\addlegendentry{version 2}

\addplot +[mark options={scale=0.6}]
  table[row sep=crcr]{%
2	0.3096721\\
3	0.6641755\\
4	1.5161428\\
5	3.4818666\\
6	7.3406607\\
8	25.5986031\\
10	63.3223526\\
12	139.4288137\\
16	475.323175\\
20	1282.0702298\\
};
\addlegendentry{version 3}

\end{axis}
\end{tikzpicture}%